\documentclass[10pt]{amsart} 
\usepackage[mathscr]{eucal}
\usepackage{amsmath,amsfonts}
\parskip=\smallskipamount

\newtheorem{theorem}{Theorem}[section]

\newtheorem{proposition}[theorem]{Proposition}

\newtheorem{corollary}[theorem]{Corollary}
\newtheorem{lemma}[theorem]{Lemma}
\newtheorem{example}[theorem]{Example}

\newtheorem{remark}[theorem]{Remark}

\makeatletter
\@addtoreset{equation}{section}
\makeatother

\newcommand{\BB}{{\mathbb B}}
\newcommand{\CC}{{\mathbb C}}
\newcommand{\NN}{{\mathbb N}}

\newcommand{\ZZ}{{\mathbb Z}}
\newcommand{\DD}{{\mathbb D}}
\newcommand{\RR}{{\mathbb R}}
\newcommand{\FF}{{\mathbb F}}
\newcommand{\TT}{{\mathbb T}}

\newcommand{\cA}{{\mathcal A}}
\newcommand{\cB}{{\mathcal B}}
\newcommand{\cC}{{\mathcal C}}
\newcommand{\cD}{{\mathcal D}}
\newcommand{\cE}{{\mathcal E}}
\newcommand{\cF}{{\mathcal F}}
\newcommand{\cG}{{\mathcal G}}
\newcommand{\cH}{{\mathcal H}}
\newcommand{\cK}{{\mathcal K}}

\newcommand{\cM}{{\mathcal M}}
\newcommand{\cN}{{\mathcal N}}
\newcommand{\cP}{{\mathcal P}}
\newcommand{\cR}{{\mathcal R}}
\newcommand{\cS}{{\mathcal S}}

\newcommand{\cX}{{\mathcal X}}
\newcommand{\cW}{{\mathcal W}}

\newdimen\expt
\def\boxit#1{\setbox0\hbox{$\displaystyle{#1}$}
      \hbox{\lower.4\expt
 \hbox{\lower3\expt\hbox{\lower\dp0
      \hbox{\vbox{\hrule height.4\expt
 \hbox{\vrule width.4\expt\hskip3\expt
      \vbox{\vskip3\expt\box0\vskip2\expt}%
 \hskip3\expt\vrule width.4\expt}\hrule height.4\expt}}}}}}
 
\begin{document}

 \pagestyle{myheadings}
\markboth{ Gelu Popescu}{Unitary invariants in multivariable operator theory}
%\pagestyle{plain}

%\begin{flushright}
 % {\it October 6, 2003}
%\end{flushright}
%\bigskip

% \title {Unitary invariants for $n$-tuples of operators}

\title [  Unitary invariants in multivariable operator theory  
  ] 
{  Unitary invariants in multivariable operator theory } 
  \author{Gelu Popescu}
\date{ April 12, 2004}
\thanks{Research supported in part by an NSF grant}

\address{Department of Mathematics, The University of Texas 
at San Antonio \\ San Antonio, TX 78249, USA}
\email{\tt gpopescu@math.utsa.edu}

%\begin{abstract}
% 

%\end{abstract}

\maketitle

%\tableofcontents

 %\clearpage

\section*{Contents}
{\it 
 Introduction
 
 Part I. Unitary invariants for $n$-tuples of operators
 \begin{enumerate}
 \item[ 1.]  Joint numerical radius
 \item[ 2.]  Euclidean operator radius
 \item[ 3.]  Joint numerical range and spectrum
 \item[ 4.]  $\rho$-operator radius
 \end{enumerate}
 
 Part II. Joint operator radii, inequalities, and applications
 
 \begin{enumerate}
 \item[ 5.]  von Neumann inequalities
 \item[ 6.]  Constrained von Neumann inequalities
 \item[ 7.]  Multivariable Haagerup-de la Harpe inequalities
 \item[ 8.]  Multivariable Fej\' er inequalities
 \end{enumerate}

 References
}

  \bigskip
  
\section*{Introduction}
 The problems considered in this paper come as a natural continuation
 of our program to develop
a {\it free} analogue of Sz.-Nagy--Foia\c s theory, 
 for row contractions. An $n$-tuple $(T_1,\ldots, T_n)$ of operators 
 acting on a Hilbert space  is called row contraction if 
 $$T_1 T_1^*+\cdots +T_nT_n^*\leq I.
 $$
 In this study, the role of the unilateral shift is played by 
 the left creation operators on the full
 Fock space with $n$ generators, $F^2(H_n)$, and the Hardy algebra $H^\infty(\DD)$
 is replaced by the noncommutative analytic Toeplitz algebra $F_n^\infty$.
   
 The    algebra   $F_n^\infty$ 
  and  its norm closed version,
  the noncommutative disc
 algebra  $\cA_n$,  were introduced by the author   \cite{Po-von} in connection
   with a multivariable noncommutative von Neumann inequality.
$F_n^\infty$  is the algebra of left multipliers of  
$F^2(H_n)$  and  can be identified with 
 the
  weakly closed  (or $w^*$-closed) algebra generated by the left creation operators
   $S_1,\dots, S_n$  acting on   $F^2(H_n)$,
    and the identity.
     The noncommutative disc algebra $\cA_n$ is 
    the  norm closed algebra generated by  
   $S_1,\dots, S_n$,  
    and the identity.
     When $n=1$, $F_1^\infty$ 
   can be identified
   with $H^\infty(\DD)$, the algebra of bounded analytic functions
    on the open unit disc. The algebra $F_n^\infty$ can be viewed as a
     multivariable noncommutative 
    analogue of $H^\infty(\DD)$.
 There are many analogies with the invariant
  subspaces of the unilateral 
 shift on $H^2(\DD)$, inner-outer factorizations,
  analytic operators, Toeplitz operators, $H^\infty(\DD)$--functional
   calculus, bounded (resp.~spectral) interpolation, etc.
The  noncommutative analytic Toeplitz algebra   $F_n^\infty$ 
  has  been studied
    in several papers
\cite{Po-charact},  \cite{Po-multi},  \cite{Po-funct}, \cite{Po-analytic},
\cite{Po-disc}, \cite{Po-poisson}, 
 \cite{ArPo}, and recently in
  \cite{DP1}, \cite{DP2},   \cite{DP},
  \cite{ArPo2}, \cite{Po-curvature},  \cite{DKP},  \cite{PPoS}, and \cite{Po-similarity}.

  In  \cite{Po-poisson}, \cite{Arv1}, \cite{Arv2},  \cite{ArPo2}, \cite{ArPo1}, \cite{DP},  
  \cite{Po-tensor}, \cite{Po-curvature}
  and  \cite{PPoS}, a good case is made 
  that the appropriate
        multivariable commutative analogue
   of $H^\infty(\DD)$  is the algebra 
   $W_n^\infty:=P_{F_s^2(H_n)} F_n^\infty|_{F_s^2(H_n)}$,
  the compression  of $ F_n^\infty$ to the symmetric Fock space
   ${F_s^2(H_n)} \subset F^2(H_n)$, which 
    was   also proved to be  the $w^*$-closed algebra generated by  the operators
   $B_i:=P_{F_s^2(H_n)} S_i|_{F_s^2(H_n)}, \ i=1,\dots, n$, and the identity. 
   The multivariable commutative disc algebra $\cA_n^c$ is the norm closed algebra 
   generated by the creation operators $B_1,\ldots, B_n$  acting on 
   the symmetric Fock space, 
    and the identity.
   Arveson showed in  \cite{Arv1}
   that the algebra $W_n^\infty$ can be seen as the   multiplier algebra of the 
reproducing kernel Hilbert space 
with reproducing kernel $K_n: \BB_n\times \BB_n\to \CC$ defined by
 $$
 K_n(z,w):= {\frac {1}
{1-\langle z, w\rangle_{\CC^n}}}, \qquad z,w\in \BB_n,
$$ 
where $\BB_n$ is the open unit ball of $\CC^n$. 

 In recent  years, there has been exciting progress in multivariable operator theory,
 especially in connection with unitary invariants for $n$-tuples of operators
 and interpolation in several variables.
 In \cite{Po-charact}, we defined the characteristic function of a row contraction 
 (a multi-analytic operator acting on Fock spaces) which, as in the classical case
 \cite{SzF-book}, turned out to be a complete unitary invariant for completely 
 non-coisometric
 row contractions. In 2000, Arveson \cite{Arv2} introduced and studied the curvature 
 and Euler characteristic associated with a row contraction with commuting entries. 
 Noncommutative
 analogues of these numerical invariants were defined and studied  by the
 author \cite{Po-curvature} and, independently, by D.~Kribs \cite{Kr}.
 Refinements of these unitary invariants were considered in 
 \cite{Po-similarity}.

This paper concerns unitary invariants for $n$-tuples $ (T_1,\ldots, T_n)$ 
 of (not necessarily commuting) bounded linear operators  acting 
 on a  Hilbert space $\cH$.
First, we introduce a  notion of 
 {\it joint numerical radius} for  $ (T_1,\ldots, T_n)$ 
    by setting
 \begin{equation*} 
  w(T_1,\ldots, T_n):=\sup\left|\sum_{\alpha\in \FF_n^+}\sum_{j=1}^n
   \left< h_\alpha, T_j h_{g_j\alpha} \right>\right|,
 \end{equation*} 
 where the supremum is taken over all families of vectors
  $\{h_\alpha\}_{\alpha\in \FF_n^+} \subset \cH$ with 
  $\sum\limits_{\alpha\in \FF_n^+}\|h_\alpha\|^2=1$, and $\FF_n^+$ is the free 
  semigroup with generators $g_1,\ldots, g_n$ and neutral element $g_0$.
  In Section 1, 
  we work out the basic properties
  of the joint numerical radius and show that it is a natural multivariable generalization 
  of the
  classical numerical radius of a bounded linear  operator $T$, i.e.,
   $$
   \omega(T):=\sup\{ |\left< Th,h\right>|:\ h\in \cH, \|h\|=1\}.$$When $n=1$, they coincide.
The joint numerical radius
  turns out to be a norm equivalent to the operator norm  on $B(\cH)^{(n)}$. Characterizations  of the unit balls
  $$\left\{(T_1,\ldots, T_n)\in B(\cH)^{(n)}:\ w(T_1,\ldots, T_n)\leq 1\right\}$$
  and
  $$\left\{(T_1,\ldots, T_n)\in B(\cH)^{(n)}:\ 
  \left\|\sum_{i=1}^n T_iT_i^*\right\|\leq 1\right\} $$
in terms of certain completely  positive maps on operator systems generated by the left
creation operators $S_1,\ldots, S_n$ on the full Fock space, lead to a multivariable version of Berger's dilation theorem (\cite{B}) and  an appropriate generalization of Berger-Kato-Stampfli mapping theorem (\cite{BS}, \cite{K}), to $n$-tuples of operators. More precisely, we prove that if $w(T_1,\ldots, T_n)\leq 1$ and $f_1,\ldots, f_k$ are in the noncommutative disc algebra $\cA_n$, then
$$
  w(f_1(T_1,\ldots, T_n),\ldots, f_k(T_1,\ldots, T_n))\leq
  \|[f_1,\ldots, f_k]\|+2\left(\sum_{j=1}^n |f_j(0)|^2\right)^{1/2}.$$Actually, a matrix-valued version of this inequality is obtained. When $n=k=1$, we obtain the classical result \cite{BS}. 
  On the other hand, we provide a multivariable operatorial generalization of Schwarz's lemma.
   In particular, if $\|[T_1,\ldots, T_n]\|\leq 1$ and 
  $f_1,\ldots, f_k\in \cA_n$ with   $f_j(0)=0$, then
  \begin{equation*}r(f_1(T_1,\ldots, T_n), \ldots, f_k(T_1,\ldots, T_n))\leq r(T_1,\ldots, T_n)
  \|[f_1,\ldots, f_k\|,
  \end{equation*}
  where $r(\cdot)$ is the joint spectral radius of  a tuple of operators. A similar inequality holds true
  if we replace $r(\cdot)$ by the operator norm. Moreover, if $\|[T_1,\ldots, T_n]\|<1$,
  then $\cA_n$ can be replaced  by the noncommutative analytic Toeplitz algebra $F_n^\infty$.
  
  We should add that if $T_1,\ldots, T_n$ are mutually commuting operators, 
  we find commutative versions for all the results of this section. In this case, 
  the noncommutative disc algebra $\cA_n$ (resp. $F_n^\infty$)
   is replaced by its commutative version $\cA_n^c$ (resp. $W_n^\infty$).

In Section 2, we present basic properties of the {\it euclidean operator radius} 
of an $n$-tuple of operators $(T_1,\ldots, T_n)$, defined by\begin{equation*}
 w_e(T_1,\ldots, T_n):= \sup_{\|h\|=1}\left(\sum_{i=1}^n\left|\left<
 T_ih,h\right>\right|^2\right)^{1/2},
 \end{equation*}in connection with the joint numerical radius and several other operator radii. In general, when $n\geq 2$,  $w_e(T_1,\ldots, T_n)\leq w(T_1,\ldots, T_n)$. However,  the two notions coincide
with the classical numerical radius if $n=1$. We show that the euclidean operator radius is a norm equivalent to the operator norm and the joint numerical radius on $B(\cH)^{(n)}$, and provide sharp inequalities.

In Section 3, we study the {\it joint (spatial) numerical range} of $(T_1,\ldots, T_n)$, defined by
$$W(T_1,\ldots, T_n):=\left\{ (\left< T_1h,h\right>,\ldots, \left<T_nh,h\right>):\ h\in \cH, \|h\|=1\right\},$$ in connection with the {\it right spectrum} $\sigma_r(T_1,\ldots, T_n)$, the joint numerical radius $w(T_1,\ldots, T_n)$, the euclidean operator radius $w_e(T_1,\ldots, T_n)$,
and the {\it joint spectral radius} $r(T_1,\ldots, T_n)$. In particular, we show that

$$ \sigma_r(T_1,\ldots, T_n)\subseteq \overline{W(T_1,\ldots, T_n)} \subseteq 
 \overline{(\CC^n)}_{w_e(T_1,\ldots, T_n)}\subseteq 
 \overline{(\CC^n)}_{w(T_1,\ldots, T_n)},
  $$
  where $\left(\CC^n\right)_r$ denotes the open unit ball of   radius $r>0$. If $T_1,\ldots, T_n$ are commuting operators and $\sigma(T_1,\ldots, T_n)$ is the Harte spectrum of $(T_1,\ldots, T_n)$ with respect to any  commutative closed subalgebra of $B(\cH)$ containing $T_1,\ldots, T_n$, and the identity, then we deduce that
  $$\sigma(T_1,\ldots, T_n)\subseteq\overline{(\CC^n)}_{r(T_1,\ldots, T_n)}\subseteq 
 \overline{(\CC^n)}_{w(T_1,\ldots, T_n)}. $$

It is well known that in general the joint numerical range 
$W(T_1,\ldots, T_n)$ is not a convex subset of $\CC^n$ if $n\geq 2$ (see \cite{GR}).
%There are a few exceptions for example when the operators $T_1,\ldots, T_n$
% are commuting normal operators (\cite{De}) or double commuting (\cite{}).
 We prove in Section 3 an analogue of Toeplitz-Hausdorff theorem (\cite{T}, \cite{Hau}) on the convexity 
 of the spatial numerical range of an operator on a Hilbert space, for 
 the joint numerical range
 of operators in the noncommutative analytic Toeplitz algebra $F_n^\infty$.
 We show that, if $f_1,\ldots, f_k\in F_n^\infty$, then
 \begin{enumerate}
 \item[1.] $ W(P_{\cP_{m}}f_1|\cP_{ m},\ldots,
   P_{\cP_{m}}f_k|\cP_{m})$ is a convex compact subset of $\CC^k$,
   where $\cP_{m}$ is the set of all polynomials in $F^2(H_n)$ of degree
   $\leq m$, and $P_{\cP_m}$ is the  projection of $F^2(H_n)$
    onto $\cP_m$;
  \item[2.]
  $\overline{W(f_1,\ldots, f_k)}$ is a convex compact subset of $\CC^k$.
 \end{enumerate}
 If $g_1,\ldots, g_k$ are elements of  Arveson's algebra
  $W_n^\infty $, we obtain a corona type
   result  which provides a characterization of the 
   Harte  spectrum  $\sigma(g_1,\ldots, g_k)$  relative to $W_n^\infty$.
   Moreover, we show that 
   $$ \sigma(g_1,\ldots, g_k)\subseteq
   \overline{W(g_1,\ldots, g_k)}\subseteq 
 \overline{(\CC^n)}_{w_e(g_1,\ldots, g_k)}\ \text{ and }$$     
   $$
 \sigma(g_1,\ldots, g_k)\subseteq
   \overline{\left(\CC^n\right)}_{r(g_1,\ldots, g_k)}
   \subseteq
   \overline{\left(\CC^n\right)}_{w_e(g_1,\ldots, g_k)}\subseteq 
 \overline{(\CC^n)}_{w(g_1,\ldots, g_k)}.
   $$
   
   In Section 4, we provide an appropriate
    generalization of the Sz.-Nagy--Foia\c s theory of $\rho$-contractions
    (\cite{Sz1},
    \cite{SzF}, \cite{SzF-book}), 
    to our multivariable setting. We say that an $n$-tuple 
    $(T_1,\ldots, T_n)$ of operators acting on a Hilbert space $\cH$  belongs to the class $\cC_\rho$, 
    $\rho>0$,
    if there is a Hilbert space $\cK\supseteq \cH$ and isometries
     $V_i\in B(\cK)$, 
    \ $i=1,\ldots, n$, with orthogonal ranges such that
    $$T_\alpha =\rho P_\cH V_\alpha |\cH\quad \text{ for any }   
    \alpha\in \FF_n^+\backslash \{g_0\},
    $$
    where 
     $T_\alpha :=  T_{i_1}T_{i_2}\cdots T_{i_k}$
if $\alpha=g_{i_1}g_{i_2}\cdots g_{i_k}$,  and $P_\cH$ is the orthogonal projection of $\cK$ on $\cH$.
If $n=1$, one can easily see that this definition is equivalent to the classical
one.
    The results of this section can be seen as the unification of
     the theory of isometric dilations for row contractions \cite{Sz1},
      \cite{SzF-book},
      \cite{Fr}, \cite{Bu}, \cite{Po-models}, \cite{Po-isometric}, 
      \cite{Po-charact}
      (which corresponds to the case $\rho=1$) 
     and Berger type dilations  of Section 1 for $n$-tuples $(T_1,\ldots, T_n)$ with the joint numerical
     radius $w(T_1,\ldots, T_n)\leq 1$ (which corresponds to the case $\rho=2$). 
    When $T_1,\ldots, T_n$ are commuting operators, we prove that 
    $(T_1,\ldots, T_n)\in \cC_\rho$ if and only if  there is a Hilbert space $\cG\supseteq \cH$ and a $*$-representation
    $\pi:C^*(B_1,\ldots, B_n)\to B(\cG)$ such that 
    $$T_\alpha =\rho P_\cH \pi(B_\alpha)| \cH \quad \text{ for any }  
    \alpha\in \FF_n^+\backslash \{g_0\},
    $$
    where $B_1,\ldots, B_n$ are the left creation operators 
    on the symmetric Fock space.
    We obtain  in this section  several intrinsic characterizations for $n$-tuples of operators of class
     $\cC_\rho$ (noncommutative and commutative case).
     
    Following the classical case (\cite{Ho1}, \cite{W}), we define the operator radius
     $\omega_\rho: B(\cH)^{(n)}\to [0,\infty)$, \ $\rho>0$, by setting
     $$
     \omega_\rho(T_1,\ldots, T_n):=\inf\left\{t>0:\ 
     \left(\frac{1} {t}T_1,\ldots, \frac{1} {t}T_n\right)\in \cC_\rho\right\}
     $$
     and 
     $\omega_\infty:=\lim\limits_{\rho\to\infty}\omega_\rho(T_1,\ldots, T_n)$.
     In particular,  $\omega_1(T_1,\ldots, T_n)$ coincides with the norm of the row
     operator $[T_1,\ldots, T_n]$, and $\omega_2(T_1,\ldots, T_n)$ coincides with
      the joint numerical radius of $(T_1,\ldots, T_n)$.
     We present in this section basic properties of the joint
     $\rho$-operator radius  $\omega_\rho(\cdot )$ and extend
       to our multivariable setting (noncommutative and commutative)
     several classical results obtained  by Sz.-Nagy and Foia\c s, Halmos,
       Berger and Stampfli, Holbrook, Paulsen, and  others   
     (\cite{B}, \cite{BS}, \cite{BS1},  \cite{Ha}, \cite{Ha1}, 
     \cite{Ho1}, \cite{Ho2}, \cite{K}, \cite{Pa-book}, \cite{Pe}, \cite{SzF},
       and \cite{W}).
     
     In Part II  of this paper, we prove
     von Neumann type inequalities \cite{vN}
      (see also \cite{Pi-book} and \cite{Pa-book})
     for arbitrary admissible (resp. strongly admissible) operator radii
     $\omega:B(\cH)^{(k)}\to [0,\infty)$. We mention that the joint operator radii
      considered in 
     Part I of this paper, namely, $\|\cdot\|$, \ $\|\cdot\|_e$, \ $w(\cdot)$,
     \ $ w_e(\cdot)$, \ $r(\cdot)$,\ 
      $r_e(\cdot)$, and $\omega_\rho(\cdot)$ ($0<\rho  \leq 2$) are strongly
      admissible, while $\omega_\rho (\cdot )$ ($\rho>2$) is admissible.
      We show  that, in general,  given a row contraction
      $[T_1,\ldots, T_n]$,  
       the inequality
      $$\omega(f_1(T_1,\ldots, T_n),\ldots, f_k(T_1,\ldots, T_n))\leq
      \omega(f_1(S_1,\ldots, S_n),\ldots, f_k(S_1,\ldots, S_n))
      $$
      holds if $f_1(S_1,\ldots, S_n),\ldots, f_k(S_1,\ldots, S_n)$ belong to  operator algebras 
      (resp. operator systems) generated by the left creation operators
       $S_1,\ldots, S_n$ and the identity such as
      the noncommutative disc algebra $\cA_n$,  the Toeplitz $C^*$-algebra
       $C^*(S_1,\ldots, S_n)$,  the noncommutative analytic Toeplitz algebra
        $F_n^\infty$, 
        the noncommutative Douglas type algebra ${\cD}_n$, 
        and the operator system
         \ $F_n^\infty (F_n^\infty)^*$. 
        The operator $f_j(T_1,\ldots, T_n)$ 
       is  defined by an appropriate functional calculus for row contractions.
       Actually, we obtain matrix-valued generalizations of the above inequality.
       An important role in this investigation is played by
         the noncommutative Poisson transforms
     associated with row contractions (see \cite{Po-poisson}, \cite{ArPo2}, 
     \cite{Po-tensor}, \cite{Po-curvature},
    and  \cite{Po-similarity}).

       If $[T_1,\ldots, T_n]$ is a row contraction  satisfying certain   constrains,
        we prove that there is a suitable invariant subspace 
        $\cE\subset F^2(H_n)$ under each operator $S_1^*, \ldots, S_n^*$, such that
        $$\omega(f_1(T_1,\ldots, T_n),\ldots, f_k(T_1,\ldots, T_n))\leq
      \omega(P_\cE f_1(S_1,\ldots, S_n)|\cE,\ldots, P_\cE f_k(S_1,\ldots, S_n)|
      \cE).
      $$
     This type  of constrained von Neumann inequalities as well as matrix-valued
     versions are 
     considered in Section 6.
     In particular, we obtain appropriate  multivariable generalizations
     of several von Neumann type inequalities obtained in
     \cite{vN}, \cite{PY}, \cite{Sz2}, \cite{Dr}, \cite{Po-von}, \cite{Po-funct}, 
     \cite{Po-disc}, \cite{Po-poisson},  \cite{Arv1},   \cite{ArPo2},
     and \cite{Po-tensor},
     to joint operator
     radii. 
     For example, 
     if $[T_1,\ldots, T_n]$  is a row contraction with  commuting entries, then we prove that
     $$\omega(f_1(T_1,\ldots, T_n),\ldots, f_k(T_1,\ldots, T_n))\leq
      \omega(f_1(B_1,\ldots, B_n),\ldots, f_k(B_1,\ldots, B_n))
      $$
        if $f_1(B_1,\ldots, B_n),\ldots, f_k(B_1,\ldots, B_n)$  are elements of 
        operator algebras generated by the creation operators 
        $B_1,\ldots, B_n$  acting on the symmetric Fock space,  such as
      the commutative disc algebra $\cA_n^c$,  the Toeplitz $C^*$-algebra
       $C^*(B_1,\ldots, B_n)$,   Arveson's algebra
        $W_n^\infty$, 
        the   Douglas type algebra ${\cD_n^c}$, and  the operator system
         \ $W_n^\infty (W_n^\infty)^*$. 
        The operator $f_j(T_1,\ldots, T_n)$ 
       is defined by an appropriate functional calculus for row contractions with commuting entries.
     Applications of these inequalities are considered in the next
      sections.

In 1992, Haagerup and de la Harpe (\cite{HD}) proved that any bounded linear
 operator of norm $1$ on a Hilbert space $\cH$ such that $T^m=0$, $m\geq 2$,
 satisfies the inequality $$ \omega(T)\leq \cos \frac {\pi} {m+1}, $$ where
 $\omega(T)$  is the numerical radius of $T$. 
 In Section 7, we apply the results of the previous sections
 to obtain several multivariable generalizations of the 
 Haagerup--de la Harpe  inequality. In particular, we show that
if  the operators $T_1,\ldots, T_n$  are such that $T_\alpha=0$ for any 
$\alpha\in \FF_n^+$ with lenght  $|\alpha|=m$
$(m\geq 2)$, then 
\begin{equation*}
  w(T_\alpha:\ |\alpha|=k)
\leq 
 \left\|\sum_{i=1}^n T_i T_i^*\right\|^{k/2} \cos\frac {\pi} 
 {\left[\frac{m-1} {k}\right]+2}
\end{equation*}
for  $1\leq k\leq m-1$, where $[x]$ is the integer part of $x$.
A similar inequality holds for the euclidean operator radius. 
Using   a result of Boas and Kac
  (\cite{BoK}) as generalized by Janssen (\cite{J}), we also obtain an
   epsilonized  version of the above inequality, when one gives up 
   the condition $T_\alpha=0$. This  extends a  recent result of
    Badea and Cassier (\cite{BC}), to our setting.
   
   Haagerup and de la Harpe   showed in \cite{HD} that their inequality 
   is equivalent to Fej\' er's inequality 
\cite {Fe} for positive trigonometric polynomials of the form
$$f(e^{i\theta}):= \sum_{k=-m+1}^{m-1} a_k e^{ik\theta}, \quad a_k\in \CC,
$$
which asserts that
$$
|a_1|\leq a_0 \cos \frac {\pi} {m+1}.
$$
   In the last section of this paper, we provide 
    multivariable noncommutative
 (resp.~commutative)
  analogues of classical inequalities  (Fej\' er \cite{Fe},
   Egerv\' ary-Sz\' azs \cite{ES}) for the coefficients of positive
   trigonometric  polynomials, 
   and  of recent extensions to positive rational functions, 
   obtained by Badea and Cassier
   \cite{BC}.

    We  obtain the following  multivariable operator-valued generalization of
     the classical
inequalities
  of Fej\' er and Egerv\'ary-Sz\'azs, to the tensor product $C^*(S_1,\ldots, S_n)\otimes B(\cH)$. 
   Let $m\geq 2$ be a nonnegative 
integer and let  $\left\{A_{(\alpha)}\right\}_{|\alpha|\leq m-1}$  be 
a sequence of operators in $B(\cH)$ such that
  the  operator
 \begin{equation*}
 \sum_{1\leq k\leq m-1}  S_\alpha^*\otimes A_{(\alpha)}+  
 I\otimes A_{0}+
 \sum_{1\leq k\leq m-1}  S_\alpha \otimes A_{(\alpha)}^*
 \end{equation*}
is positive.  
 Then, 
\begin{equation*}
w(A_{(\alpha)}:\ |\alpha|=k)
\leq 
 \|A_{0}\|\cos\frac {\pi} {\left[\frac{m-1} {k}\right]+2}
\end{equation*}
 for  $1\leq k\leq m-1$. Similar 
  inequalities are obtained    for  the reduced group 
  $C^*$-algebra $C^*_{red}(\FF_n)$ and the full group $C^*$-algebra $C^*(\FF_n)$.
  In  particular, we deduce that
   if
   $$
   f=\sum_{1\leq|\alpha|\leq m-1} \overline{a}_\alpha S_\alpha^* +a_0 I+
   \sum_{1\leq|\alpha|\leq m-1} a_\alpha S_\alpha
   $$
   is a positive polynomial in $C^*(S_1,\ldots, S_n)$, then
   $$
   \left(\sum_{|\alpha|=k} |a_\alpha|^2\right)^{1/2}\leq a_0
   \cos\frac{\pi} {\left[\frac{m-1}{k}\right]+2}
   $$
    for  $1\leq k\leq m-1$. 
    Similar results hold for $C^*_{red}(\FF_n)$ and $C^*(\FF_n)$.
  
 There is a large literature relating to the classical numerical 
 range and no attempt has been made to compile a comprehensive 
 list of reference here. For general proprieties, we refer to \cite{Ha},
 \cite{BD1}, \cite{BD2},  and \cite{GR}.
 Finally, we want to acknowledge that we were strongly influenced in writing
 this paper by the work of Haagerup and de la Harpe \cite{HD},
  Arveson \cite{Ar},  Paulsen \cite{Pa-book}, 
  Sz.-Nagy and Foia\c s \cite{SzF-book},  and Badea and Cassier \cite{BC}.

\smallskip

 \section*{Part I. Unitary invariants for $n$-tuples of operators}
  \section{Joint numerical radius}
  \label{JNR}
 
We introduce a  notion of 
 joint numerical radius for  $n$-tuples $ (T_1,\ldots, T_n)$ 
 of   operators acting  on a  Hilbert space $\cH$, 
    by setting
 \begin{equation}\label{def1}
  w(T_1,\ldots, T_n):=\sup\left|\sum_{\alpha\in \FF_n^+}\sum_{j=1}^n
   \left< h_\alpha, T_j h_{g_j\alpha} \right>\right|,
 \end{equation} 
 where the supremum is taken over all families of vectors
  $\{h_\alpha\}_{\alpha\in \FF_n^+} \subset \cH$ with 
  $\sum\limits_{\alpha\in \FF_n^+}\|h_\alpha\|^2=1$.
  In this section, we present basic properties of the joint numerical range, which 
   lead to a multivariable version of Berger's dilation theorem 
   (\cite{B}) and  an appropriate generalization of Berger-Kato-Stampfli
    mapping theorem (\cite{BS}, \cite{K}), to $n$-tuples of operators
    with $w(T_1,\ldots, T_n)\leq 1$. We also obtain a multivariable operatorial  generalization of Schwarz's lemma.
    If $T_1,\ldots, T_n$ are mutually commuting operators, 
    we find commutative versions
    for all the results of this section.

    First notice that, when $n=1$, our  numerical radius coincides
   with the classical numerical radius of an operator $T\in B(\cH)$ 
   (the algebra of all bounded linear operators on $\cH$), i.e.,
   $$
   \omega(T):=\sup\{ |\left< Th,h\right>|:\ h\in \cH, \|h\|=1\}.
   $$
  Indeed,   if $\sum\limits_{k=0}^\infty \|h_k\|^2<\infty$ and
   $f_m(\theta):=\sum\limits_{k=0}^m e^{ik\theta}h_k$, $h_k\in \cH$, 
    $m=1,2,\ldots
$,  then
  \begin{equation*}
  \begin{split}
  \left|\sum _{k=0}^m \left< h_k, Th_{k+1}\right>\right|&=
  \left|\frac{1} {2\pi}\int_0^{2\pi}\left<  f_m(\theta), 
  Tf_m(\theta)\right>e^{i\theta} d\theta\right|\\
  &\leq \omega(T) \frac{1} {2\pi}\int_0^{2\pi} \|f_m(\theta)\|^2 d\theta\\
  &=\omega(T) \sum_{k=0}^m \|h_k\|^2.
  \end{split}
  \end{equation*}
  Taking $m\to\infty$, we obtain
    $w(T)\leq \omega(T)$.
  On the other hand, let $h\in \cH$, $\|h\|=1$, and define
  $h_k:=\lambda^k\sqrt{1-|\lambda|^2} h$, $k=0,1,\ldots$, 
  where $\lambda\in \DD:=\{z\in \CC:\ |z|<1\}$. 
  It is easy to see that $\sum\limits_{k=0}^\infty \|h_k\|^2 =1$ and 
  \begin{equation*}
  \left|\sum _{k=0}^\infty \left< h_k, Th_{k+1}\right>\right|=
  |\lambda| |\left<Th,h\right>|
  \end{equation*} 
  for any  $\lambda\in \DD$. Hence, we deduce the inequality
   $w(T)\geq \omega(T)$, which proves
   our assertion, i.e., $w(T)=\omega(T)$.

Let $H_n$ be an $n$-dimensional complex  Hilbert space with orthonormal basis
$e_1$, $e_2$, $\dots,e_n$, where $n\in \{1,2,\dots\}$ or $n=\infty$.
  We consider the full Fock space  of $H_n$ defined by
$$F^2(H_n):=\bigoplus_{k\geq 0} H_n^{\otimes k},$$ 
where $H_n^{\otimes 0}:=\CC 1$ and $H_n^{\otimes k}$ is the (Hilbert)
tensor product of $k$ copies of $H_n$.
Define the left creation 
operators $S_i:F^2(H_n)\to F^2(H_n), \  i=1,\dots, n$,  by
$$
 S_i\psi:=e_i\otimes\psi, \quad  \psi\in F^2(H_n).
$$
Let $\FF_n^+$ be the unital free semigroup on $n$ generators 
$g_1,\dots,g_n$, and the identity $g_0$.
The length of $\alpha\in\FF_n^+$ is defined by
$|\alpha|:=k$, if $\alpha=g_{i_1}g_{i_2}\cdots g_{i_k}$, and
$|\alpha|:=0$, if $\alpha=g_0$.
We also define
$e_\alpha :=  e_{i_1}\otimes e_{i_2}\otimes \cdots \otimes e_{i_k}$  
 and $e_{g_0}= 1$.
It is  clear that 
$\{e_\alpha:\alpha\in\FF_n^+\}$ is an orthonormal basis of $F^2(H_n)$.
 We denote by $\cP$ the set of all polynomials in $F^2(H_n)$, 
 i.e., all the elements of the form 
 $$p=\sum_{|\alpha|\leq m } a_\alpha e_\alpha, 
 \quad a_\alpha\in \CC,  \ m=0,1,2,\ldots.
 $$
 If $T_1,\dots,T_n\in B(\cH)$, define 
$T_\alpha :=  T_{i_1}T_{i_2}\cdots T_{i_k}$
if $\alpha=g_{i_1}g_{i_2}\cdots g_{i_k}$, and 
$T_{g_0}:=I_\cH$.  If $p$ is a polynomial as the one above, we set
$p(T_1,\ldots, T_n):=
\sum_{|\alpha|\leq m } a_\alpha T_\alpha$.

  The  joint spectral radius associated
 with an  $n$-tuple of operators 
$ (T_1,\ldots, T_n)$  is given by
$$
r(T_1,\ldots, T_n):=\lim_{k\to \infty}\left\|\sum_{|\alpha|=k} T_\alpha 
T_\alpha^*\right\|^{1/2k}.
$$
  In general,   $r(T_1,\ldots, T_n)\neq r(T_1^*,\ldots, T_n^*)$ if $n\geq 2$. Note
   for example that 
   $$r(S_1,\ldots, S_n)=1\  \text{  and  }\ r(S_1^*,\ldots, S_n^*)=\sqrt{n},
   $$
   where $S_1,\ldots, S_n$ are the left creation operators on the full Fock space. 

  For simplicity, throughout
 this paper, $[T_1,\ldots, T_n]$ denotes either the $n$-tuple
  $(T_1,\ldots, T_n)$ or the operator row matrix $[T_1\ \cdots \ T_n]$.

In what follows we present basic properties of the joint numerical radius.
  \begin{theorem}
  \label{propri}
 The joint numerical radius $w:B(\cH)^{(n)}\to [0,\infty)$   
 for $n$-tuples of operators  satisfies the following properties:
 \begin{enumerate}
 \item[(i)]
 $w(T_1,\ldots, T_n)=0$ if and only if $T_1=\cdots=T_n=0$;
 \item[(ii)]
 $w(\lambda T_1,\ldots, \lambda T_n)=|\lambda|~w(T_1,\ldots, T_n) $
 for any $\lambda\in \CC$;
 \item[(iii)]
 $w(T_1+T_1',\ldots, T_n+T_n')\leq w(T_1,\ldots, T_n)+w(T_1',\ldots, T_n');$
 \item[(iv)]
 $w(U^*T_1U,\ldots, U^*T_nU)= w(T_1,\ldots, T_n)$ for any unitary operator
  $U:\cK\to \cH$;
  \item [(v)]
  $w(X^*T_1X,\ldots, X^*T_n X)\leq \|X\|^2 w(T_1,\ldots, T_n)$ for any 
    operator $X:\cK\to \cH$;
 \item[(vi)]
 $\frac{1}{2} \left\|\sum\limits_{i=1}^n T_iT_i^*\right\|^{1/2}\leq
 w(T_1,\ldots, T_n)\leq  \left\|\sum\limits_{i=1}^n T_iT_i^*\right\|^{1/2};$
 \item[(vii)]
 $r(T_1,\ldots, T_n)\leq w(T_1,\ldots, T_n);$
 \item [(viii)]
 $w(I_\cE\otimes T_1,\ldots, I_\cE\otimes T_n)=w(T_1,\ldots, T_n)$ for 
 any separable
 Hilbert space $\cE$;
 \item[(ix)]
 $w$ is a continuous map in the norm topology.
 \end{enumerate}
  \end{theorem}
  \begin{proof}
  The first three properties follow easily using the definition
  \eqref{def1}, and  show that  the joint numerical 
  radius is a norm on $B(\cH)^{(n)}$. 
  To prove relations (iv) and (v), let 
  $\{k_\alpha\}_{\alpha\in \FF_n^+}$   be an arbitrary sequence  of vectors in 
  $ \cK$
   such that
  $\sum\limits_{\alpha\in \FF_n^+} \|k_\alpha\|^2=1$. Fix an operator $X:\cK\to \cH$ and
  define the vectors $h_\alpha:= \frac {1} {M} Xk_\alpha$, \ $\alpha\in \FF_n^+$,
  where $M:=\left(\sum\limits_{\alpha\in \FF_n^+} \|Xk_\alpha \|^2\right)^{1/2}$.
  Notice that 
  $\sum\limits_{\alpha\in \FF_n^+} \|h_\alpha\|^2=1$ and $M\leq \|X\|$.
  On the other hand, we have
  \begin{equation*}
  \begin{split}
  \left| \sum_{ \alpha\in \FF_n^+}\sum_{j=1}^n\left< k_\alpha, X^*T_jXk_{g_j\alpha}
  \right>\right|&\leq
  \left| \sum_{ \alpha\in \FF_n^+}\sum_{j=1}^n\left< h_\alpha, T_jh_{g_j\alpha}
  \right>\right|\|X\|^2\\
  &\leq\|X\|^2w(T_1,\ldots, T_n).
  \end{split}
  \end{equation*}
  Taking the supremum over  all sequences
  $\{k_\alpha\}_{\alpha\in \FF_n^+}\subset \cK$   
    with 
  $\sum\limits_{\alpha\in \FF_n^+} \|k_\alpha\|^2=1$, we deduce inequality (v).
  A closer look     reveals that, 
  when $X=U$ is a unitary operator,
   we have equality  in the above inequality. Therefore, relation (iv) holds true.
   
  Now, let us prove (vi).
  Any vector $x\in F^2(H_n)$ with $\|x\|=1$ has the form 
  $x=\sum\limits_{\alpha\in \FF_n^+} e_\alpha\otimes h_\alpha$, where
  the sequence $\{h_\alpha\}_{\alpha\in \FF_n^+}\subset \cH$  is such that
  $\sum\limits_{\alpha\in \FF_n^+} \|h_\alpha\|^2=1$.
  If $S_1,\ldots, S_n$ are the left creation operators on the full Fock space
 $F^2(H_n)$,
  note that
  \begin{equation*}
  \begin{split}
  \left|\left<\left(\sum_{j=1}^n S_j\otimes T_j^*\right) x,x\right>\right| 
  &=\left|\left< \sum_{ \alpha\in \FF_n^+}\sum_{j=1}^n e_{g_j\alpha}\otimes T_j^*h_\alpha,
   \sum_{\beta\in
  \FF_n^+} e_\beta\otimes h_\beta\right> \right|\\
  &=\left| \sum_{ \alpha,\beta\in \FF_n^+}\sum_{j=1}^n 
   \left<e_{g_i\alpha}, e_\beta \right> 
  \left<T_j^*h_\alpha, h_\beta\right>\right|\\
  &=\left| \sum_{\alpha \in \FF_n^+}\sum_{j=1}^n
  \left<h_{\alpha}, T_jh_{g_i \alpha} \right>\right|.
  \end{split}
  \end{equation*} 
  Hence, we infer that 
 \begin{equation}\label{reduct}
 w(T_1,\ldots, T_n)=w(S_1\otimes T_1^*+\cdots + S_n\otimes T_n^*).
 \end{equation}
 On the other hand, 
  it is well known that the classical numerical radius of an operator $X$ 
  satisfies the inequalities  $\frac{1} {2}\|X\|\leq w(X)\leq \|X\|$.
  Using relation \eqref{reduct} and taking into account that
 the left creation operators
  have orthogonal ranges,
   %
  % $S_i^*S_j=\delta_{ij}$, ~$i,j=1,\ldots, n$,
    we have
   \begin{equation*}\begin{split}
   w(T_1,\ldots, T_n)&=w\left(\sum_{i=1}^n S_i\otimes T_i^*\right)\leq \left\|
   \sum_{i=1}^n S_i\otimes T_i^*\right\|\\
   &=\left\| \left(\sum_{i=1}^n S_i^*\otimes T_i\right)
   \left(\sum_{i=1}^n S_i\otimes T_i^*\right)\right\|^{1/2}\\
   &= \left\| \sum_{i=1}^n T_iT_i^*\right\|^{1/2}=\|[T_1,\ldots, T_n]\|.
   \end{split}
   \end{equation*}
 Similarly, one can prove the first inequality in (vi).
  
  To prove (vii), notice that,
  since $S_i^* S_j=\delta_{ij}I$, \ $i,j=1,\ldots, n$, we have
 \begin{equation*}
 \begin{split}
  r\left(\sum_{i=1}^n  S_i\otimes T_i^*\right) =&\lim_{k\to \infty}
 \left\|\sum_{|\alpha|=k} 
    S_\alpha\otimes T_{\tilde \alpha}^*\right\|^{1/k}\\
    &=
    \lim_{k\to \infty}
 \left\| 
    \left(\sum_{|\alpha|=k}S_\alpha^*\otimes T_{\tilde \alpha}\right)
    \left(\sum_{|\alpha|=k}S_\alpha\otimes T_{\tilde \alpha}^*
    \right)\right\|^{1/2k}\\
    &=\lim_{k\to \infty}
 \left\|I\otimes \sum_{|\alpha|=k} T_{\tilde \alpha}T_{\tilde \alpha}^*
 \right\|^{1/2k}\\
 &= r(T_1,\ldots, T_n).
 \end{split}
 \end{equation*}
  Consequently, we deduce that
 \begin{equation*}\begin{split}
 r(T_1,\ldots, T_n)&=r(S_1\otimes T^*_1+\cdots +S_n\otimes T^*_n)\\
 &\leq w(S_1\otimes T^*_1+\cdots +S_n\otimes T^*_n)\\
 &=w(T_1,\ldots, T_n)
 \end{split}
 \end{equation*}
   Therefore, the inequality (vii) is established. 
  
  To prove (viii),  we use again  relation \eqref{reduct} and the classical result
  $w(I_\cE\otimes X)=w(X)$.  Notice that
  \begin{equation*}
  \begin{split}
  w(I_\cE\otimes T_1,\ldots, I_\cE\otimes T_n)&=
  w\left(\sum_{i=1}^n S_i\otimes (I_\cE\otimes T_i^*)\right)\\
  &=w\left( I_\cE\otimes \left( \sum_{i=1}^n S_i\otimes T_i^*\right) \right)\\
  &=w(T_1,\ldots, T_n).
  \end{split}
  \end{equation*}
  The property (ix) follows imediately from (vi).
  The proof is complete.
  \end{proof}

  In general,  $w(T_1,\ldots, T_n)\neq w(T_1^*,\ldots, T_n^*)$ if $n\geq 2$.
  Indeed, we have
  $$w(S_1^*,\ldots, S_n^*)\geq r(S_1^*,\ldots, S_n^*)=\sqrt{n}
  $$  and
 $w(S_1,\ldots, S_n)=1$.

  \begin{corollary}\label{pro}
  If $(T_1,\ldots, T_n)\in B(\cH)^{(n)}$, then 
  \begin{equation*}
  \begin{split}
  \|[T_1,\ldots, T_n]\|&=\|S_1\otimes T^*_1+\cdots +S_n\otimes T^*_n\|,\\
  r(T_1,\ldots, T_n)&=r(S_1\otimes T^*_1+\cdots +S_n\otimes T^*_n), \text{ and }\\
  w(T_1,\ldots, T_n)&=w(S_1\otimes T^*_1+\cdots +S_n\otimes T^*_n),
  \end{split}
  \end{equation*}
  where $S_1,\ldots, S_n$ are the left creation operators on the full Fock space
  $F^2(H_n)$.
  \end{corollary}

We need to recall from  \cite{Po-charact},
\cite{Po-multi}, \cite{Po-von},  \cite{Po-funct}, and  \cite{Po-analytic} 
 a few facts
 concerning multi-analytic  (resp. multi-Toeplitz) operators on Fock spaces.
   We say that 
 a bounded linear
  operator 
$M$ acting from $F^2(H_n)\otimes \cK$ to $ F^2(H_n)\otimes \cK'$ is 
 multi-analytic
if 
\begin{equation}
M(S_i\otimes I_\cK)= (S_i\otimes I_{\cK'}) M\quad 
\text{\rm for any }\ i=1,\dots, n.
\end{equation}
Notice that $M$ is uniquely determined by the operator
$\theta:\cK\to F^2(H_n)\otimes \cK'$, which is   defined by ~$\theta k:=M(1\otimes k)$, \ $k\in \cK$, 
 and
is called the  symbol  of  $M$. We denote $M=M_\theta$. Moreover, 
$M_\theta$ is 
 uniquely determined by the ``coefficients'' 
  $\theta_{(\alpha)}\in B(\cK, \cK')$, which are given by
 $$
\left< \theta_{(\tilde\alpha)}k,k'\right>:= \left< \theta k, e_\alpha 
\otimes k'\right>=\left< M_\theta(1\otimes k), e_\alpha 
\otimes k'\right>,\quad 
k\in \cK,\ k'\in \cK',\ \alpha\in \FF_n^+,
$$
where $\tilde\alpha$ is the reverse of $\alpha$, i.e., $\tilde\alpha= g_{i_k}\cdots g_{i_1}$ if
$\alpha= g_{i_1}\cdots g_{i_k}$.
Note that 
$$
\sum\limits_{\alpha \in \FF_n^+} \theta_{(\alpha)}^*   \theta_{(\alpha)}\leq 
\|M_\theta\|^2 I_\cK.
$$
 We can associate with $M_\theta$ a unique formal Fourier expansion 
\begin{equation}\label{four}
M_\theta\sim \sum_{\alpha \in \FF_n^+} R_\alpha \otimes \theta_{(\alpha)},
\end{equation}
where $R_i:=U^* S_i U$, \ $i=1,\ldots, n$, are the right creation operators
on $F^2(H_n)$ and 
$U$ is the (flipping) unitary operator on $F^2(H_n)$ mapping 
  $e_{i_1}\otimes e_{i_2}\otimes\cdots\otimes e_{i_k}$ into
 $e_{i_k}\otimes\cdots\otimes e_{i_2}\otimes e_{i_1} $.
 Since $M_\theta$ acts like its Fourier representation on ``polynomials'', 
  we will identify them for simplicity.
Based on the noncommutative von Neumann
   inequality (\cite{Po-funct}, \cite{Po-tensor})
  we proved that
  $$M_\theta=\text{\rm SOT}-\lim_{r\to 1}\sum_{k=0}^\infty \sum_{|\alpha|=k}
   r^{|\alpha|} R_\alpha\otimes \theta_{(\alpha)},
   $$
   where, for each $r\in (0,1)$ the series converges in the uniform norm.
The set of  all multi-analytic operators in 
$B(F^2(H_n)\otimes \cK,
F^2(H_n)\otimes \cK')$  coincides  with   
$R_n^\infty\bar \otimes B(\cK,\cK')$, 
the WOT closed algebra generated by the spatial tensor product, where $R_n^\infty=U^* F_n^\infty U$
(see \cite{Po-analytic}  and \cite{Po-central}).
We also denote $\theta(R_1,\ldots, R_n):=M_\theta$.

 An operator  $T$  acting on $F^2(H_n)\otimes \cE$, where $\cE$ is 
  a Hilbert space, 
  is called multi-Toeplitz  if 
 $$
 (S_i\otimes I_\cE)^* T (S_j\otimes I_\cE)=\delta_{ij} I\quad 
  \text{ for any }\ i,j=1,\ldots, n. 
 $$ 
Basic properties for multi-Toeplitz operators on Fock spaces can be found in 
 \cite{Po-multi}, \cite{Po-analytic}.   

 Given $\alpha, \beta\in \FF_n^+$, we say that $\alpha>\beta$ if
  $\alpha=\beta \omega$
  for some  $\omega\in \FF_n^+\backslash\{g_0\}$. We denote 
  $\omega:=\alpha\backslash \beta$.  
  A kernel $K:\FF_n^+\times \FF_n^+\to B(\cH)$ is called multi-Toeplitz if 
  $K(g_0,g_0)=I_\cH$ and 
  $$
  K(\alpha, \beta)= 
  \begin{cases}
    K(\alpha\backslash \beta, g_0) &\text{ if } \alpha>\beta\\
    K( g_0, \beta\backslash \alpha) &\text{ if } \alpha<\beta\\
    0\quad &\text{ otherwise}.
   \end{cases}$$
  It is said to be positive definite provided that
  $$
  \sum_{\alpha,\beta\in \FF_n^+} \left<K(\alpha,\beta) h(\beta), 
  h(\alpha)\right>\geq 0
  $$
  for all finitely supported functions $h$ from $\FF_n^+$ into $\cH$. 
  
  The next result, which will play an 
  important role in this section,  provides a characterization for the $n$-tuples of operators 
  with the joint numerical radius 
$w(T_1,\ldots, T_n)\leq 1$.
  
  \begin{theorem}\label{C2}
  Let $T_1,\ldots, T_n\in B(\cH)$ and let $\cS\subset C^*(S_1,\ldots, S_n)$
  be the operator system defined by
  \begin{equation}\label{Sys}
  \cS:=\{p(S_1,\ldots, S_n)+q(S_1,\ldots, S_n)^*:\ p,q\in \cP\}.
  \end{equation}
  Then the following statements are equivalent:
  \begin{enumerate}
 \item[(i)] 
  $w(T_1,\ldots, T_n)\leq 1$;
 \item[(ii)] The map $\Psi:\cS\to B(\cH)$ defined by
 \begin{equation*} 
 \Psi(p(S_1,\ldots, S_n)+q(S_1,\ldots, S_n)^*):=
 p(T_1,\ldots, T_n)+q(T_1,\ldots,T_n)^*
 + (p(0)+\overline{q(0)})I
 \end{equation*}
 is completely positive.
  \end{enumerate}
  \end{theorem}
  
  \begin{proof} First,  we prove that (i)$\implies$(ii).
  For each
    $q=0,1,\ldots$, 
   define  the operator
   \begin{equation}\label{mq}
   M_q:=
   P_{\cP_q\otimes \cH}\left(\sum_{1\leq |\alpha|\leq q} 
   R_\alpha^*\otimes T_{\tilde \alpha}+ 2I+\sum_{1\leq |\alpha|\leq q} 
   R_\alpha\otimes T_{\tilde \alpha}^*
   \right)|\cP_q\otimes \cH,
   \end{equation}
  where   $P_{\cP_q\otimes \cH}$ is the orthogonal 
  projection of $F^2(H_n)\otimes \cH$ onto  $\cP_q\otimes \cH$, and 
  $\cP_q$ is the set of all polynomials
  of degree $\leq q$ in $F^2(H_n)$.
   Let $\cR$ be the operator system obtained from $\cS$ by replacing the left
   creation operators $S_1,\ldots, S_n$ with the right creation operators
   $R_1,\ldots, R_n$.
   Any element $G\in \cR\otimes B(\CC^m)$, $m=1,2\ldots$,  has the form
   $$
   G=G(R_1,\ldots, R_n):=\sum_{|\alpha|\leq q} R_\alpha^*\otimes B_{(\alpha)}+
   \sum_{|\alpha|\leq q} R_\alpha\otimes A_{(\alpha)}
   $$
   for some operators $A_{(\alpha)}, B_{(\alpha)}\in B(\CC^m)$, 
   $|\alpha|\leq q$.
   Note that $G$ is a multi-Toeplitz operator acting on the Hilbert space 
   $F^2(H_n)\otimes \CC^m$, i.e.,
   $$
   (S_i^*\otimes I_{\CC^m}) G(S_j\otimes I_{\CC^m})=\delta_{ij}G,\qquad 
   i,j=1,\ldots, n.
   $$
    Assume that $w(T_1,\ldots, T_n)\leq 1$ and $G\geq 0$. Then, according 
    to the Fej\' er type 
   factorization theorem of \cite{Po-multi}, there 
   exists a multi-analytic operator 
   $\Theta\in B(F^2(H_n)\otimes \CC^m)$ such that 
   $G= \Theta^* \Theta$ and 
   $$
   (S_\alpha^*\otimes I_{\CC^m}) \Theta (1\otimes h)=0
   $$
   for any $h\in \CC^m$ and $|\alpha|>q$.
   Therefore, 
   $\Theta=\sum\limits_{|\alpha|\leq q} R_\alpha \otimes \Theta_{(\alpha)}$ 
   for some operators $\Theta_{(\alpha)}\in B(\CC^m)$ and 
   $$
   G=\Theta^*\Theta= \sum_{|\alpha|, |\beta|\leq q}
    R_\alpha^* R_\beta \otimes \Theta_{(\alpha)}^*\Theta_{(\beta)}.
    $$
    Since $R_i^* R_j=\delta_{ij}$, $i,j=1,\ldots, n$, we have
    \begin{equation}\label{G=}
    G=\sum_{ \alpha>\beta, ~   |\alpha|, |\beta|\leq q
     }R_{\alpha\backslash \beta}^*\otimes 
    \Theta_{(\alpha)}^*\Theta_{(\beta)}
    +
    \sum_{ \beta\geq \alpha, ~ |\alpha|, |\beta|\leq q
     }R_{\beta\backslash \alpha}\otimes 
    \Theta_{(\alpha)}^*\Theta_{(\beta)}.
   \end{equation} 
   To prove  that (i)$\implies$(ii), it is enough to show that  the operator
   $G(T_1,\ldots, T_n)+G(0)I$ is positive whenever $G\geq 0$.
    To this end, define the multi-Toeplitz  kernel 
   $K_{T,2}:\FF_n^+\times \FF_n^+\to B(\cH)$
   by setting
   $$
   K_{T,2}(\alpha, \beta):= 
   \begin{cases}
   T_{\beta\backslash \alpha} &\text{ if } \beta>\alpha\\
   2I  &\text{ if } \alpha=\beta\\
    T_{\alpha\backslash \beta}^*  &\text{ if } \alpha>\beta\\
    0\quad &\text{ otherwise}.
   \end{cases} 
   $$
   Note that 
   $$
   G(T_1,\ldots, T_n)+G(0)I=\sum_{|\alpha|,|\beta|\leq q} K_{T,2}(\alpha, \beta)
   \otimes \Theta_{(\alpha)}^*\Theta_{(\beta)}.
   $$
   Hence, for any element 
   $\sum\limits_{i=1}^p h_i\otimes z_i\in \cH\otimes \CC^m $, 
   $p=1,2,\ldots$, we
   obtain
   \begin{equation*}
   \begin{split}
   \Biggl<[G(T_1,\ldots, T_n)\Biggr.&\Biggl.+G(0)I]\left(\sum_{i=1}^p
    h_i\otimes z_i\right),
   \sum_{i=1}^p h_i\otimes z_i\Biggr>\\
   &=
   \sum_{i, j=1}^p 
   \sum_{|\alpha|,|\beta|\leq q}\left< K_{T,2}(\alpha, \beta) h_i\otimes 
   \Theta_{(\alpha)}^*\Theta_{(\beta)}z_i, h_j\otimes z_j
   \right>\\
   &=
   \sum_{i, j=1}^p 
   \sum_{|\alpha|,|\beta|\leq q}\left< K_{T,2}(\alpha, \beta) h_i\otimes 
   \Theta_{(\beta)}z_i, h_j\otimes \Theta_{(\alpha)}z_j
   \right>\\
   &=
   \sum_{|\alpha|,|\beta|\leq q}\left< 
   (K_{T,2}(\alpha, \beta)\otimes I_{\CC^m}) \left(\sum_{i=1}^p h_i\otimes 
   \Theta_{(\beta)}z_i\right), \sum_{i=1}^ph_i\otimes \Theta_{(\alpha)}z_i
   \right>.
   \end{split}
   \end{equation*}
   We need to prove that the operator matrix
   $[K_{T,2}(\alpha,\beta)]_{|\alpha|,|\beta|\leq q}$ is positive.
   Since $[K_{T, 2}(\alpha,\beta)]_{|\alpha|,|\beta|\leq q}$ is 
   the   matrix representation of the operator $M_q$ with respect to
    the decomposition
   $\cP_q\otimes \cH= \oplus_{i=1}^N\cH$, where $N:=1+n+\cdots + n^q$,
    it is enough
   to prove that $M_q\geq 0$.
   Let the operator  $A_q: \cP_q\otimes \cH\to \cP_q\otimes \cH$ be defined by 
   \begin{equation} \label{aq}
   A_q:=P_{\cP_q\otimes \cH}(R_1^*\otimes T_1+
   \cdots +R_n^*\otimes T_n)|\cP_q\otimes \cH.
   \end{equation}
    Since $A_q^{q+1}=0$, we have
   $$M_q=\sum_{k=0}^q A_q^k +\sum_{k=0}^q {A_q^*}^k=
   (I-A_q)^{-1}+ (I-A_q^*)^{-1}.
   $$
   Let $h:=(I-A_q) k$,  $k\in \cP_q\otimes \cH$, and note that
   \begin{equation*}
   \begin{split}
   \left<M_qh,h\right>&=\left<k+(I-A_q^*)^{-1}(I-A_q)k, (I-A_q)k\right>\\
   &=\left<k,(I-A_q)k\right>+\left<(I-A_q)k,k \right>\\
   &=2\|k\|^2-2\,\text{\rm Re} \left<A_qk, k\right>.
   \end{split}
   \end{equation*}
   Therefore, $M_q\geq 0$ if and only if 
   $ \text{\rm Re}\left<A_qk,k\right>\leq 1$ for all
  $k\in \cP_q\otimes \cH$ with $\|k\|=1$.
  We prove now the latter condition is equivalent to $\omega(A_q)\leq 1$.
  Since one implication is clear, assume that  
    $\text{\rm Re} \left<A_qk,k \right>\leq 1$ for all
  $k\in \cP_q\otimes \cH$ with $\|k\|=1$. If $\omega(A_q) > 1$, then 
  there exists
  $x=\sum\limits_{|\alpha|\leq q} e_\alpha \otimes h_\alpha
  \in \cP_q\otimes \cH$ such that $\|x\|=1$ and $|\left<A_qx,x \right>|>1$.
  A short calculation reveals that
  \begin{equation*}
  \begin{split}
  \left|\left<A_qx,x\right>\right|&=
  \left|\left<\left(\sum_{j=1}^n R_i^*\otimes T_j\right) x,x\right>\right|\\
  &=
  \left|\left< \sum_{|\beta|\leq q}e_\beta\otimes h_\beta, 
  \sum_{j=1}^n \sum_{|\alpha|\leq q} e_{\alpha g_j}\otimes T_j^* h_\alpha
  \right>\right|\\
  &=\left|\sum_{j=1}^n \sum_{|\alpha|,|\beta|\leq q}\left<e_\beta, 
  e_{\alpha g_j} \right>
  \left<T_jh_\beta, h_\alpha \right>\right|\\
  &=\left|\sum_{j=1}^n \sum_{|\alpha|\leq q}
  \left<T_jh_{\alpha g_j}, h_\alpha \right>\right|.
  \end{split}
  \end{equation*} 
 Now, let $\theta_\alpha\in [0,2\pi]$, \ $|\alpha|\geq 1$ ,  be such that
   $\left<T_jh_{g_j}, h_{g_0}\right>e^{i \theta_{g_j}}\geq 0$ for any $j=1,\ldots, n$, and 
  $$\left<T_jh_{\alpha g_j}, 
  h_{\alpha}\right>e^{i \theta_{\alpha g_j}} e^{-i \theta_{\alpha}}\geq 0$$
  for any 
   $j=1,\ldots, n$, and $\alpha\in \FF_n^+$ with $|\alpha|\leq q-1$.
   Let $h_\alpha':=e^{i\theta_\alpha} h_\alpha$ if $0<|\alpha|\leq q$ 
   and $h_{g_0}':= h_{g_0}$, and note that 
   $$
   \left<T_jh_{\alpha g_j}', h_\alpha'\right>=
   \left<T_jh_{\alpha g_j}, 
  h_{\alpha}\right>e^{i \theta_{\alpha g_j}} e^{-i \theta_{\alpha}}\geq 0.
  $$
  Therefore, $\sum\limits_{|\alpha|\leq q} \|h_\alpha'\|^2=
  \sum\limits_{|\alpha|\leq q} \|h_\alpha\|^2=1$
  and, if $y:= \sum\limits_{|\alpha|\leq q} e_\alpha\otimes h_\alpha'$, then we obtain
  \begin{equation*}\begin{split}
  \left<A_qy,y\right>&=\sum_{j=1}^n\sum_{|\alpha|\leq q}
   e^{i \theta_{\alpha g_j}} e^{-i \theta_{\alpha}}\left<T_jh_{\alpha g_j}, 
  h_{\alpha}\right>\\
  &\geq \left|\sum_{j=1}^n\sum_{|\alpha|\leq q}
  \left<T_jh_{\alpha g_j}, 
  h_{\alpha}\right>\right|>1.
  \end{split}
  \end{equation*}
  Hence,
  $1<\text{\rm Re}\left<A_qk',k'\right>$, which is a contradiction. 
   This proves that we must have 
  $\omega(A_q)\leq 1$.
   Therefore, 
   \begin{equation}\label{Mo}
   M_q\geq 0 \ \text{ if and only if } \  \omega(A_q)\leq 1.
   \end{equation}
    According to Corollary \ref{pro}, we have   
   \begin{equation}\label{def2}
   w(T_1,\ldots, T_n)=w(S_1\otimes T_1^*+\ldots + S_n\otimes T_n^*).
   \end{equation}
   Since $U^* S_iU=R_i$, \ $i=1,\ldots, n$, and $U$ is unitary,  
   the unitary invariance of the joint numerical radius and relation \eqref{def2}
   imply
   $$ w(T_1,\ldots, T_n)=w(R_1\otimes T_1^*+\ldots + R_n\otimes T_n^*).
   $$
  Therefore,  if  $w(T_1,\ldots, T_n)\leq1$, then
   $\omega(A_q)\leq w(R_1\otimes T_1^*+\ldots + R_n\otimes T_n^*)\leq 1$.
   Now, relation \eqref{Mo} implies condition   (ii) of the theorem.

  Conversely, assume that   (ii) holds.
   Define the multi-Toeplitz kernel $K_{S,1}:\FF_n^+\times \FF_n^+\to B(\cH)$
   by setting
   $$
   K_{S,1}(\alpha, \beta):= 
   \begin{cases}
   S_{\beta\backslash \alpha}  &\text{ if } \beta>\alpha\\
   I &\text{ if } \alpha=\beta\\
    S_{\alpha\backslash \beta}^*  &\text{ if } \alpha>\beta\\
    0 &\text{ othewise}.
   \end{cases} 
   $$
  Since $S_i^*S_j=\delta_{ij} I$, $i,j=1,\ldots, n$, we have
  \begin{equation*}
  \begin{split}
  \sum_{|\alpha|, |\beta|\leq q}
  \left<K_{S,1}(\alpha, \beta)h_\beta, h_\alpha\right>&=
  \sum_{|\alpha|, |\beta|\leq q}
  \left<S_\alpha^* S_\beta h_\beta, h_\alpha\right>\\
  &=\left\| \sum_{|\alpha|\leq q}S_\alpha h_\alpha
  \right\|^2\geq 0
  \end{split}
  \end{equation*}
  for any $\{h_\alpha\}_{|\alpha|\leq q}\subset \cH$ and $q=1,2,\ldots, $.
  Therefore, the operator matrix $[K_{S,1}(\alpha, \beta)]_{|\alpha|,|\beta|\leq q}$
  is positive. Since $\Psi$ is a completely positive map and 
  $$
  [\Psi(K_{S,1}(\alpha, \beta))]_{|\alpha|,|\beta|\leq q}=
  [K_{T,2}(\alpha, \beta)]_{|\alpha|,|\beta|\leq q},
  $$
  we deduce that $M_q\geq 0$. Hence and using \eqref{Mo},  we get $\omega(A_q)\leq 1$ for any 
  $q=1,2,\ldots.$
  According to relations \eqref{aq} and \eqref{def2}, we obtain
  $w(T_1,\ldots, T_n)\leq 1$. The proof is complete.

  \end{proof}

 \begin{lemma}\label{spec-ra}
 If 
 $(T_1,\ldots, T_n)$ is an $n$-tuple of operators on a Hilbert space
   $\cH$, then 
   $$
   r(T_1,\ldots, T_n)= r(R_1\otimes T_1^*+\cdots +R_n\otimes T_n^*).
   $$
   Moreover, if the spectral radius  $r(T_1,\ldots, T_n)\leq 1$ and $z_i\in \DD$, for any $i=1,\ldots, n$,
    then 
   the operator $I-\sum\limits_{i=1}^n z_iR_i\otimes 
T_i^* $ is invertible and 
   \begin{equation}\label{inverse}
   \left(I-\sum_{i=1}^n z_iR_i\otimes 
   T_i^*\right)^{-1}=I+\sum_{k=1}\sum_{|\alpha|=k} 
   z_\alpha R_\alpha\otimes T_{\tilde \alpha}^*,
   \end{equation}
   where the series converges in the operator norm.
 \end{lemma}
 \begin{proof}
The first part of the theorem follows from Corollary \ref{pro}
  and the fact 
 that $R_i=U^* S_iU$, \ $i=1,\ldots, n$, where $U$ is the (flipping)
  unitary operator.
  Assume  now that 
 $ r(T_1,\ldots, T_n)\leq 1$ and  $z_i\in \DD$,\ $i=1,\ldots, n$. Then we have
 $$
 r\left(\sum_{i=1}^n  z_iR_i\otimes T_i^*\right)=
  r(\bar z_1T_1,\ldots, \bar z_nT_n)\leq \gamma\, r(T_1,\ldots, T_n)\leq \gamma,
 $$
 where $\gamma:=\max \{|z_1|, \ldots, |z_n|\}<1$.
 This proves that  the operator 
 $I-\sum_{i=1}^n z_iR_i\otimes T_i^*$ is invertible.
  On the other hand, the root test implies the convergence of the
   series in the operator
  norm. The equality in \eqref{inverse} is now obvious.
 \end{proof}

   The following  result  provides  a characterization for row contractions
    as well as
     a new proof of the noncommutative von Neumann
  inequality obtained in 
  \cite{Po-von}.

  \begin{theorem}\label{vN}
  Let $T_1,\ldots, T_n\in B(\cH)$ and let $\cS\subset C^*(S_1,\ldots, S_n)$
  be the operator system defined by \eqref{Sys}.
  Then the following statements are equivalent:
  \begin{enumerate}
 \item[(i)] 
  $T_1T_1^*+\cdots + T_n T_n^*\leq 1$;
 \item[(ii)] The map $\Psi:\cS\to B(\cH)$ defined by
 \begin{equation*} 
 \Psi(p(S_1,\ldots, S_n)+q(S_1,\ldots, S_n)^*):=
 p(T_1,\ldots, T_n)+q(T_1,\ldots,T_n)^*
 \end{equation*}
 is completely positive.
 \item[(iii)] The spectral radius $r(T_1\ldots, T_n)\leq 1$ and  the multi-Toeplitz operator
  $A_r $ defined by
  $$
  A_r:= \sum_{k=1}^\infty \sum_{|\alpha|=k}r^k 
  R_\alpha^*\otimes T_{\tilde \alpha} +I+
  \sum_{k=1}^\infty \sum_{|\alpha|=k}r^k 
  R_\alpha\otimes T_{\tilde \alpha}^*
  $$
  is positive for any $0<r<1$, where the convergence is in the operator norm.
  \end{enumerate}
  \end{theorem}
  
  \begin{proof} 
  First we prove that (i)$\implies$(ii).
  For each $q=0,1,\ldots$, 
  define  the operator
   \begin{equation}\label{nq}
   N_q:=
   P_{\cP_q\otimes \cH}\left(\sum_{1\leq |\alpha|\leq q} 
   R_\alpha^*\otimes T_{\tilde \alpha}+ I+\sum_{1\leq |\alpha|\leq q} 
   R_\alpha\otimes T_{\tilde \alpha}^*
   \right)|\cP_q\otimes \cH,
   \end{equation}
  where   $P_{\cP_q\otimes \cH}$ is the orthogonal 
  projection of $F^2(H_n)\otimes \cH$ onto  $\cP_q\otimes \cH$, and 
  $\cP_q$ is the set of all polynomials
  of degree $\leq q$ in $F^2(H_n)$.
  The multi-Toeplitz  kernel 
   $K_{T,1}:\FF_n^+\times \FF_n^+\to B(\cH)$ is given 
   by  
   $$
   K_{T,1}(\alpha, \beta):= 
   \begin{cases}
   T_{\beta\backslash \alpha} &\text{ if } \beta>\alpha\\
   I  &\text{ if } \alpha=\beta\\
    T_{\alpha\backslash \beta}^*  &\text{ if } \alpha>\beta\\
    0\quad &\text{ otherwise}.
   \end{cases} 
   $$
    As in the proof of Theorem \ref{C2}, the Fej\' er type factorization theorem
    of \cite{Po-analytic} shows that any positive  operator
    $G\in \cR\otimes B(\CC^m)$, $m=1,2,\ldots$, has the form \eqref{G=}, and
    similar calculations show that
   $$
   G(T_1,\ldots, T_n)=\sum_{|\alpha|,|\beta|\leq q} K_{T,1}(\alpha, \beta)
   \otimes \Theta_{(\alpha)}^*\Theta_{(\beta)}.
   $$
   Hence, for any element 
   $\sum\limits_{i=1}^p h_i\otimes z_i\in \cH\otimes \CC^m $, 
   $p=1,2,\ldots$, one 
   obtains
   \begin{equation*}
   \begin{split}
   \Biggl<[G(T_1,\ldots, T_n)]&\left(\sum_{i=1}^p 
    h_i\otimes z_i\right),
   \sum_{i=1}^p h_i\otimes z_i\Biggr>\\
   &=\sum_{|\alpha|,|\beta|\leq q}\left< 
   (K_{T,1}(\alpha, \beta)\otimes I_{\CC^m}) \left(\sum_{i=1}^p h_i\otimes 
   \Theta_{(\beta)}z_i\right), \sum_{i=1}^ph_i\otimes \Theta_{(\alpha)}z_i
   \right>.
   \end{split}
   \end{equation*}
   Therefore, to show that $G(T_1,\ldots, T_n)\geq 0$, 
    we need to prove that the  matrix
   $[K_{T,1}(\alpha,\beta)]_{|\alpha|,|\beta|\leq q}$ is positive.
   Since $[K_{T, 1}(\alpha,\beta)]_{|\alpha|,|\beta|\leq q}$ is 
   the   matrix representation of the operator $N_q$ with respect to
    the decomposition
   $\cP_q\otimes \cH= \oplus_{i=1}^N\cH$, where $N:=1+n+\cdots + n^q$,
    it is enough
   to prove that $N_q\geq 0$.
   Let the operator  $A_q: \cP_q\otimes \cH\to \cP_q\otimes \cH$ be defined by 
   \begin{equation*} \
   A_q:=P_{\cP_q\otimes \cH}(R_1^*\otimes T_1+
   \cdots +R_n^*\otimes T_n)|\cP_q\otimes \cH.
   \end{equation*}
    Since $A_q^{q+1}=0$, we have
   \begin{equation*}
   \begin{split}
   N_q&=\sum_{k=1}^q A_q^k +I+\sum_{k=1}^q {A_q^*}^k=
   (I-A_q)^{-1}+ (I-A_q^*)^{-1}-I\\
   &=(I-A_q)^{-1}(I-A_qA_q^*)(I-A_q^*)^{-1}.
   \end{split}
   \end{equation*}
  Since $[T_1,\ldots, T_n]$ is a row contraction and
    $R_i^* R_j=\delta_{i j} I$,
\ $i,j=1,\ldots, n$,
   we have
   
  \begin{equation*}
   \begin{split}
   A_q A_q^*&= \left(\sum_{i=1}^n R_i^*\otimes T_i\right) P_{\cP_q\otimes \cH}
   \left(\sum_{i=1}^n R_i\otimes T_i^*\right)\\
   &\leq \sum_{i,j=1}^n R_i^*R_j\otimes T_iT_j^*\leq I.
   \end{split}
   \end{equation*}
  Therefore $N_q\geq 0$, which shows that $G(T_1,\ldots, T_n)\geq 0$. 
  This proves that
  (i)$\implies$ (ii).
  
   Assume now that (ii) holds.
  According to Arveson's extension theorem \cite{Ar}, there exists a completely positive map
   $\Phi:C^*(R_1,\ldots, R_n)\to B(\cH)$ extending $\Psi$. 
   Using  Stinespring's dilation theorem \cite{St}, there exists a $*$-representation 
   $\pi:C^*(R_1,\ldots, R_n)\to B(\cK)$ such that 
   $$
   \Phi(x)=P_\cH \pi(x)|\cH, \quad x\in C^*(R_1,\ldots, R_n).
  $$
  Since  $R_i^* R_j=\delta_{i j} I$, \ $i,j=1,\ldots, n$,  it is clear that
  $$\|[T_1,\ldots, T_n]\|\leq \|[\pi(R_1),\ldots, \pi(R_n)]\|\leq 1.
  $$
  This completes the proof of  (ii)$\implies$(i).
  
  Now we prove the implication (i) $\implies$ (iii).
  Assume that condition (i) holds.
  Since
   $$
   r\left(\sum_{i=1}^n R_i\otimes T_i^*\right)=r(T_1,\ldots, T_n)\leq 1,
   $$
    Lemma \ref{spec-ra} 
   shows that
   the operator $I-\sum\limits_{i=1}^n rR_i\otimes T_i^*$ is invertible and 
   $$
   \left(I-\sum_{i=1}^n rR_i\otimes T_i^*\right)^{-1}= 
   \sum_{k=0}^\infty \sum_{|\alpha|=k} r^k R_\alpha \otimes T_{\tilde \alpha}^*.
   $$
   Moreover, we have
   $$
   A_r=\left(I-\sum_{i=1}^n rR_i\otimes T_i^*\right)^{-1}
   \left[ I\otimes \left( I-r^2\sum_{i=1}^n T_i T_i^*\right) \right]
   \left(I-\sum_{i=1}^n rR_i^*\otimes T_i\right)^{-1}\geq 0
   $$
   for any $0<r<1$. 
   
    To prove  the implication
   (iii) $\implies$ (i), assume that $A_r\geq 0$  for any $0<r<1$.
  %%%%%%%%% 
   %
  First we show that
  \begin{equation}\label{Ar}
  \left< A_r\left( \sum_{|\beta|\leq q} e_\beta\otimes h_\beta \right), 
   \sum_{|\gamma|\leq q} e_\gamma\otimes h_\gamma\right> 
   = \sum_{|\beta|, |\gamma|\leq q}\left< K_{T,1,r}(\gamma, \beta)
    h_\beta, h_\gamma\right>,
  \end{equation} 
  where
  the multi-Toeplitz kernel  $K_{T,1,r}:\FF_n^+\times \FF_n^+\to B(\cH)$
  is defined  by  
   $$
   K_{T,1,r}(\alpha, \beta):= 
   \begin{cases}
  r^{|\beta\backslash \alpha|} T_{\beta\backslash \alpha}
   &\text{ if } \beta>\alpha\\
   I  &\text{ if } \alpha=\beta\\
    r^{|\alpha\backslash \beta|}(T_{\alpha\backslash \beta})^* 
     &\text{ if } \alpha>\beta\\
    0\quad &\text{ otherwise}.
   \end{cases} 
   $$
   Since $r(T_1,\ldots, T_n)\leq 1$, the series
    $\sum\limits_{k=0}^\infty \sum\limits_{|\alpha|=k} r^k
     R_\alpha\otimes T_{\tilde \alpha}^*$
    is convergent in norm for any $0<r<1$.
    Note that if $\{h_\beta\}_{|\beta|\leq q}\subset \cH$, then 
    \begin{equation*}
    \begin{split}
    \left< \left(\sum_{k=0}^\infty \sum_{|\alpha|=k} r^k
     R_\alpha\otimes T_{\tilde \alpha}^*\right)
     \right.& \left.\left(\sum_{|\beta|\leq q} e_\beta\otimes h_\beta \right), 
   \sum_{|\gamma|\leq q} e_\gamma\otimes h_\gamma\right>\\
   &= 
     \sum_{k=0}^\infty \sum_{|\alpha|=k}\left<\sum_{|\beta|\leq q} r^k
     R_\alpha e_\beta\otimes T_{\tilde \alpha}^*h_\beta,
    \sum_{|\gamma|\leq q} e_\gamma\otimes h_\gamma\right> \\
    &=\sum_{\alpha\in \FF_n^+}\sum_{|\beta|, |\gamma|\leq q}
    r^{|\alpha|} \left< e_{\beta \tilde{\alpha}}, e_\gamma\right>
    \left<T_{\tilde\alpha}^* h_\beta, h_\gamma\right>\\
    &=
    \sum_{ \gamma\geq\beta; ~|\beta|, |\gamma|\leq q}
    r^{|\gamma\backslash \beta|}
     \left<T_{\gamma\backslash \beta}^* h_\beta, h_\gamma\right>\\
     &=
     \sum_{\gamma\geq\beta; ~|\beta|, |\gamma|\leq q}\left<K_{T,1,r}
      (\gamma, \beta)h_\beta, h_\gamma\right>.
    \end{split}
    \end{equation*}
   Hence, taking into account that $K_{T,1,r}
      (\gamma, \beta)=K_{T,1,r}^*
      ( \beta, \gamma)$,  relation \eqref{Ar} follows.
  Now, since $A_r\geq 0$, we must have
  $[K_{T,1,r}(\alpha,\beta)]_{|\alpha|, |\beta|\leq q}\geq 0$ for any $0<r<1$
  and $q=0,1,\ldots$.
  Taking $r\to 1$, we  get
   $[K_{T,1}(\alpha,\beta)]_{|\alpha|, |\beta|\leq q}\geq 0$ for any 
    $q=0,1,\ldots$.
    Using the Naimark type dilation theorem of   \cite{Po-posdef},
     we deduce that $[T_1,\ldots, T_n]$ is a row
     contraction. The proof is complete.
  \end{proof}

  \begin{corollary}$($\cite{Po-von}$)$\label{VN}
  If $[T_1,\ldots, T_n]$ is a row contraction, then 
  $$
  \|p(T_1,\ldots, T_n)\|\leq \|p(S_1,\ldots, S_n)\|
  $$
  for any polynomial $p\in \cP$. Moreover the map $\Psi$ of Theorem $\ref{vN}$
  is completely contractive and can be extended to $\cA_n+ \cA_n^*$.
  \end{corollary}

A  consequence of the noncommutative von Neumann inequality and the noncommutative commutant lifting
theorem \cite{Po-isometric} is the following 
multivariable matrix-valued  generalization of Schwarz's lemma. We denote by $M_m$ the set of all $m\times m$ complex matrices.
 
\begin{theorem}\label{Schw} Let $[T_1,\ldots, T_n]$ be a row contraction and $m,k\geq 1$.
 
\begin{enumerate}
 
\item[(i)]
 
If $F_j:=F_j(S_1,\ldots, S_n)\in M_m\otimes \cA_n$ such that   $F_j(0)=0$,  $j=1,\ldots, k$, then \begin{equation*} 
\omega(F_1(T_1,\ldots, T_n), \ldots, F_k(T_1,\ldots, T_n))\leq \omega(T_1,\ldots, T_n) \|[F_1,\ldots, F_k]\|,
\end{equation*} where $\omega$ is the norm or the joint spectral radius.

\item[(ii)] If 
  $\cA_n$ is replaced by $F_n^\infty$ in (i),  then  the above inequality   holds true
  for any
 completely non-coisometric row contraction, in particular if \ $\|[T_1,\ldots, T_n]\|<1$.
 \item[(iii)] If $T_iT_j=T_jT_i$  for any  $i,j=1,\ldots, n$, then $(i)$ and $(ii)$ remain true if we replace 
 $\cA_n$ (resp. $F_n^\infty$)   by the commutative  algebra $\cA_n^c$ (resp. $W_n^\infty$).
 \end{enumerate}
 
\end{theorem}
 
\begin{proof}
 
Let $X_q:=X_q(S_1,\ldots, S_n)$ be the row matrix with entries $F_\alpha$, $\alpha\in \FF_n^+$ with $|\alpha|=k$, and denote $X_q:=[F_\alpha:\ |\alpha|=q]$. Note that $X_q\in M_{m,mk^q}\otimes \cA_n$, i.e., a $m\times mk^q$ matrix with entries in $\cA_n$.
Taking into account the structure  of the elements in $ F_n^\infty$, the condition $F_j(0)=0$, $j=1,\ldots, k$,  implies $$X_q(S_1,\ldots, S_n)=\sum_{|\beta|=q} (I_m\otimes S_\beta) \Phi_\beta(S_1,\ldots, S_n)$$
for some operator matrices $\Phi_\beta(S_1,\ldots, S_n)\in M_{m,mk^q}\otimes \cA_n$. Since the operators  $S_\beta$, $|\beta|=q$,  are isometries with orthogonal ranges, we have $$X_q(S_1,\ldots, S_n)^*X_q(S_1,\ldots, S_n)= \sum_{|\beta|=q}\Phi_\beta(S_1,\ldots, S_n)^*\Phi_\beta(S_1,\ldots, S_n).$$
 
Now, using the  
the noncommutative von Neumann inequality, we obtain
 \begin{equation*} \begin{split} 
 \|X_q(T_1,\ldots, T_n)\|&=\left\|\sum_{|\beta|=q} (I_m\otimes T_\alpha) \Phi_\beta(T_1,\ldots, T_n)\right\| \\ 
 &\leq \left\|[I_m\otimes T_\beta: |\beta|=q] \left[\begin{matrix}\Phi_\beta(T_1,\ldots, T_n)\\:\\  |\beta|=q\end{matrix}\right]\right\|\\ 
 &\leq \|[I_m\otimes T_\beta:\ |\beta|=q]\|\left\|\sum_{|\beta|=q}^n \Phi_\beta(S_1,\ldots, S_n)^*\Phi_\beta(S_1,\ldots, S_n)\right\|^{1/2}\\ 
 &=\|[T_\beta: \ |\beta|=q]\| \|X_q(S_1,\ldots, S_n)\|\\
 &=\|[T_\beta: \ |\beta|=q]\|\|[F_1,\ldots, F_k]\|^q.
 \end{split} \end{equation*}
 Therefore, we have
 $$
 \left\| \sum_{|\beta|=q}F_\beta(T_1,\ldots, T_n) F_\beta(T_1,\ldots, T_n)^*\right\|^{1/2q}\leq 
 \left\| \sum_{|\beta|=q}T_\beta T_\beta^*\right\|^{1/2q}\|\|[F_1,\ldots, F_k]\|.
 $$
 Taking $q=1$ or $q\to\infty$, we obtain the inequality (i).
 Part (ii) follows in a similar manner if one uses the $F_n^\infty$-functional calculus
 for completely non-coisometric row contractions (see \cite{Po-funct}).
 
  To prove  (iii), let $[T_1,\ldots, T_n]$ be  a completely non-coisometric
   row contraction  with commuting  entries. Let 
   $G_j:= G_j(B_1,\ldots, B_n)\in M_m\otimes W_n^\infty$, $j=1,\ldots, k$,  be  such that $[G_1,\ldots, G_k]$ is a row   contraction  with $G_j(0)=0$, $j=1,\ldots,k$.
   Since $G:=[G_1,\ldots, G_k]\in M_{m,mk}\otimes W_n^\infty$, according to \cite{Po-tensor},  we have
   $$
   G\left(I_{\CC^{mk}}\otimes B_i\right)= \left(I_{\CC^{m}}\otimes B_i\right)G
   $$
   for any $i=1,\ldots, n$.
     Using the noncommutative commutant 
   lifting theorem \cite{Po-isometric} (see  also \cite{ArPo2}), we find
    $F(S_1,\ldots, S_n)\in M_{m,mk}\otimes F_n^\infty$  with the properties
    $$ P_{\CC^m\otimes F^2_s(H_n)} F(S_1,\ldots, S_n)|(\CC^{mk} \otimes F^2_s(H_n))
    =G(B_1,\ldots, B_n)
    $$
    and $\|F(S_1,\ldots, S_n)\|=\|G(B_1,\ldots, B_n)\|$.
    Since $G(0)=0$, we also have $F(0)=0$. On the other hand, the commutativity of the operators  $T_1,\ldots, T_n$
    implies $F(T_1,\ldots, T_n)=G(T_1,\ldots, T_n)$. Applying now part (ii) of the theorem
    to $F(S_1,\ldots, S_n)$, we obtain
    $$\|G(T_1,\ldots, T_n)\|=\|F(T_1,\ldots, T_n)\|\leq \|[T_1,\ldots, T_n]\|.
    $$
   Assume now that $G(B_1,\ldots, B_n)\in M_m\otimes \cA_n^c$ and let  $[T_1,\ldots, T_n]$ be  an  
   arbitrary row contraction with commuting entries.
   Note that $[rT_1,\ldots, rT_n]$ is completely non-coisometric for $0<r<1$. Therefore, we have
   $$
   \|G(rT_1,\ldots, rT_n)\| \leq \|[rT_1,\ldots, rT_n]\|.
   $$
   Since $G(rT_1,\ldots, rT_n)$ converges in norm to  $G(T_1,\ldots, T_n)$ as $r\to 1$, 
    the above inequality implies
   $\|G(T_1,\ldots, T_n)\|\leq \|[T_1,\ldots, T_n]\|$.
 The proof is complete.\end{proof}
  
  Note that in the particular case when $n=m=k=1$, $\|T\|<1$,  and $f\in H^\infty(\DD)$ with $f(0)=0$ and $\|f\|_\infty\leq 1$,
  we get 
  $$\|f(T)\|\leq \|T\|\quad \text{ and } \quad r(f(T))\leq r(T).$$

    Now we can obtain the following multivariable  version of 
     Berger's dilation 
    theorem (see  \cite{B}).
  \begin{theorem} \label{berger} Let $T_1,\ldots, T_n\in B(\cH)$. Then 
  $w(T_1,\ldots, T_n)\leq 1$ if and only if
   there exists a Hilbert
   space $\cK\supseteq \cH$ and isometries with orthogonal ranges 
   $V_1,\ldots, V_n \in B(\cK)$ such that
   $$
   T_\sigma=2P_\cH V_\sigma |\cH \quad  \text{ for any }  
   \sigma\in \FF_n^+\backslash \{g_0\}. 
   $$
  \end{theorem}
 \begin{proof} According to the proof of Theorem \ref{C2}, the multi-Toeplitz
  kernel
 $\frac{1} {2} K_{T,2}: \FF_n^+\times \FF_n^+\to B(\cH)$ is positive definite
 if and only if the numerical radius $w(T_1,\ldots, T_n)\leq 1$.
 Using the Naimark type dilation theorem of \cite{Po-posdef},
  the result follows.
 \end{proof}
Let us point out another proof of Theorem \ref{berger}. According to
 Arveson's extension theorem \cite{Ar},  the map $\Psi$ of Theorem \ref{C2} has
 a completely positive extension $\tilde\Psi:C^*(S_1,\ldots, S_n)\to B(\cH)$.
 By Stinespring's theorem \cite{St}, there exists a representation 
 $\pi:C^*(S_1,\ldots, S_n)\to B(\cK)$ and an operator $V\in B(\cH,\cK)$ such that
 $\frac {1}{2}\tilde\Psi(x)=V^* \pi(x)V$ for any $x\in C^*(S_1,\ldots, S_n)$.
 Since $\Psi(I)=2I$, we have $V^*V=I$. Identifying $\cH$ with $V\cH$, we get
 $\frac {1}{2}\tilde\Psi(x)=P_\cH \pi(x)|\cH$. In particular,  this shows that
 $T_\sigma= 2P_\cH \pi(S_\sigma)|\cH$, for any $\sigma\in \FF_n^+$,
  $\sigma\neq g_0$. Since $V_i:=\pi(S_i)$, $i=1,\ldots, n$, are isometries with 
  orthogonal ranges, the result follows.

 \begin{corollary}\label{simi}
  Let
$(T_1,\ldots, T_n)\in B(\cH)^{(n)}$ be an  $n$-tuple of operators  with
 the numerical 
radius
$
w(T_1,\ldots, T_n)\leq1.
$
Then there is a completely bounded map $\Phi:\cA_n+\cA_n^*\to B(\cH)$ such that 
$$
\Phi(p(S_1,\ldots, S_n)+q(S_1,\ldots, S_n)^* )=p(T_1,\ldots, T_n)+
q(T_1,\ldots, T_n)^*
$$
for any polynomials $p, q\in \cP$, 
and $\|\Phi\|_{cp}\leq 3$.
Moreover,  the $n$-tuple $(T_1,\ldots, T_n)$ is simultaneously similar 
to a row contraction.
\end{corollary}
\begin{proof}
Note that Theorem \ref{berger} implies
$$
p(T_1,\ldots, T_n)+q(T_1,\ldots, T_n)^*
=P_\cH[-(p(0)+\overline{q(0)})I_\cK+2(p(V_1,\ldots, V_n)+p(V_1,\ldots, V_n)^*)]|\cH
$$
for any polynomials $p$ and $q$. Using the noncommutative von Neumann
 inequality (see
Corollary \ref{VN}), we deduce that the map $\Phi$ is completely bounded
and  $\|\Phi\|_{cp}\leq 3$.
According to 
Theorem 2.4 of \cite{Po-disc}   (which uses Paulsen's similarity
 result \cite{Pa1}) applied to $\Phi|\cA_n$,
 there exists a row contraction
$[A_1,\ldots, A_n]$ and an invertible operator $X$ such that
$T_i=X^{-1} A_i X$, $i=1,\ldots, n$.
This completes the proof.
\end{proof}
  
 According to Corollary \ref{simi}, if $f(S_1,\ldots, S_n)\in \cA_n$, then
 $$
 f(T_1,\ldots, T_n):=\Phi(f(S_1,\ldots, S_n))=\lim_{m\to\infty} p_m(T_1,\ldots, T_n),
 $$
 where $\{p_m(S_1,\ldots, S_n)\}_{m=1}^\infty$ is a sequence of polynomials in $\cA_n$
 which converges to $f(S_1,\ldots, S_n)$ in the operator norm. Using Theorem
 \ref{berger}, we have
 \begin{equation}\label{ffpf}  
  f(T_1,\ldots, T_n)+f(0)I=2P_\cH f(V_1,\ldots, V_n)|\cH.
  \end{equation}
  Hence, and using  the fact that $f(rV_1,\ldots, rV_n)$  converges to $f(V_1,\ldots, V_n)$
  in norm as $r\to 1$ (see \cite{Po-poisson}), we deduce that
  $$
   f(T_1,\ldots, T_n)=\lim_{r\to 1} f(rT_1,\ldots, rT_n),
   $$
   where the limit exists in the norm topology and $$
f(rT_1,\ldots, rT_n):=\sum_{k=0}^\infty \sum_{|\alpha|=k} r^k a_\alpha T_\alpha
$$ 
for  $f(S_1,\ldots, S_n)=\sum_{k=0}^\infty \sum_{|\alpha|=k} a_\alpha S_\alpha$.
   
  The next result is a    multivariable generalization of Berger-Kato-Stampfli 
  mapping theorem  (see \cite{BS}),  to $n$-tuples of operators.
  
  \begin{theorem}\label{BKS1}
  If $(T_1,\ldots, T_n)\in B(\cH)^{(n)}$ is an $n$-tuple of operators with 
   joint numerical radius  $w(T_1,\ldots, T_n)\leq 1$ and $F_1,\ldots, F_k\in M_m\otimes\cA_n$, then
  $$
  w(F_1(T_1,\ldots, T_n),\ldots, F_k(T_1,\ldots, T_n))\leq
  \|[F_1,\ldots, F_k]\|+2\left(\sum_{j=1}^n \|F_j(0)\|^2\right)^{1/2}.
  $$
  \end{theorem}
  \begin{proof} First, we prove the result when $m=1$,  $g_j\in \cA_n$,
    $\|[g_1,\ldots, g_k] \|\leq 1$, and
  $g_j(0)=0$ for 
  any $j=1,\ldots, k$.
  Let $\cS\subset C^*(S_1,\ldots, S_n)$
  be the operator system defined by
  \begin{equation*} 
  \cS:=\{p(S_1,\ldots, S_n)+q(S_1,\ldots, S_n)^*:\ p,q\in \cP\}.
  \end{equation*}
   According to Theorem \ref{C2}, 
   the map $\Psi:\cS\to B(\cH)$ defined by
 \begin{equation*} 
 \Psi(p(S_1,\ldots, S_n)+q(S_1,\ldots, S_n)^*):=
 p(T_1,\ldots, T_n)+q(T_1,\ldots,T_n)^*
 + (p(0)+\overline{q(0)})I
 \end{equation*}
 is completely positive.
  Let $L_1,\ldots, L_k$ be the  be the left creation operators on the full
  Fock space $F^2(H_k)$ and let $\cX\subset C^*(L_1,\ldots, L_k)$ be 
  the operator system defined by
  $$\cX:=\{r(L_1,\ldots, L_k)+s(L_1,\ldots, L_k)^*:\  r,s \text{ are polynomials}\}.
  $$
  Given $g_j=g_j(S_1,\ldots, S_n)\in \cA_n$, $j=1,\ldots, k$,
   with $\|[g_1,\ldots, g_k] \|\leq 1$,
  define the mapping
  $\Phi_{g_1,\ldots, g_k}:\cX\to C^*(S_1, \ldots, S_n)$ by setting
  $$
  \Phi_{g_1,\ldots, g_k}(r(L_1,\ldots, L_k)+s(L_1,\ldots, L_k)^*):=
   r(g_1 ,\ldots,
   g_k )
  +s(g_1 ,\ldots, g_k )^*.
  $$
  Using Theorem \ref{vN}, we  conclude that 
  the map $\Phi_{g_1,\ldots, g_k}$ is completely positive.
  Assume for the moment that $g_j(S_1,\ldots, S_n)$, $j=1,\ldots, k$, are
   polynomials in $S_1,\ldots, S_n$.
  Notice that $\Psi\circ \Phi_{g_1,\ldots, g_k}:\cX\to B(\cH)$  is completely
  positive and, since $g_j(0)=0$ for any $j=1,\ldots, k$,  we have
   \begin{equation*}\begin{split}
  \Psi\circ \Phi_{g_1,\ldots, g_k}(r(L_1,\ldots, L_k)&+s(L_1,\ldots, L_k)^*) \\
  &=
 r( Y_1,\ldots, Y_k)+s(Y_1,\ldots,Y_k))^*
 + (r(0)+\overline{s(0)})I,
  \end{split}
  \end{equation*}
 where $Y_j:= g_j(T_1,\ldots, T_n)$, $j=1,\ldots, k$.
Since $\Psi\circ \Phi_{g_1,\ldots, g_k}$ is completely positive,
  we can apply Theorem \ref{C2}
   to deduce that $w(Y_1,\ldots, Y_k)\leq 1$.
 Since the numerical radius is norm continuous, 
 the conclusion remains true for any $g_j(S_1,\ldots, S_n)$ in the noncommutative disc
 algebra $\cA_n$. 
 
 Now, for each $j=1,\ldots, k$, define
 $$
 g_j(S_1,\ldots, S_n):=\frac{1} {K}\left(f_j(S_1,\ldots, S_n)- f_j(0) I\right),
 $$
 where $K:=\|[f_1,\ldots, f_k]\| +\left(\sum_{j=1}^n |f_j(0)|^2\right)^{1/2}$.
 Notice that $\|[g_1,\ldots, g_k]\|\leq 1$ and $g_j(0)=0$ for any $j=1,\ldots, k$.
 Applying the first part of the proof to the $n$-tuple $(g_1,\ldots, g_k)$,
 we obtain
 $$w(f_1(T_1,\ldots, T_n)-f_1(0) I,\ldots, f_k(T_1,\ldots, T_n)-f_k(0) I)\leq K.
 $$
 Taking into account that the joint numerical radius is a norm
 on $B(\cH)^{(k)}$, the result follows.
 The matrix-valued extension of this result follows in a similar manner.
 The proof is complete.
 \end{proof}

  \begin{corollary} \label{BKS}If $w(T_1,\ldots, T_n)\leq 1$ and $F(S_1,\ldots, S_n)\in M_m\otimes\cA_n$, 
  $F(0)=0$, then
  $$
w(F(T_1,\ldots, T_n))\leq \|F(S_1,\ldots, S_n)\|.$$
  
\end{corollary}

 Note that if $n=m=1$, we find again the cassical result of  Berger-Kato-Stampfli
 \cite{BS}, \cite{K}. 
  We recall that the  numerical range  of an operator $X\in B(\cH)$ is defined by
 $$W(X):=\{\left<Xh,h\right>:\ h\in \cH, \|h\|=1\}.
 $$
 \begin{theorem}\label{Real-ine}
 If $w(T_1,\ldots, T_n)\leq 1$ and $f\in \cA_n$ with $\text{\rm Re}\, f\geq 0$, then
 $$
 \text{\rm Re} \, W(f(T_1,\ldots, T_n))\geq -\text{\rm Re}\, f(0).
 $$
 \end{theorem}
 \begin{proof} According to  the remarks following Corollary \ref{simi}, we have 
 \begin{equation}\label{ffpf2}
 f(T_1,\ldots, T_n)+f(0)I= 2P_\cH f(V_1, \ldots, V_n)|\cH
 \end{equation}
 for some isometries with orthogonal ranges acting on a Hilbert space $\cK\supseteq \cH$.
 If $0<r<1$, then  $\|[rV_1,\ldots, rV_n]\|<1$, and according to \cite{Po-isometric}, there
 is a Hilbert space $\cG$ such that $\cK$ can be identified with a subspace of 
 the Hilbert space
 $\cG\otimes F^2(H_n)$ and 
 $$ 
 rV_i^*=(I_\cG\otimes S_i^*)|\cK, \quad i=1,\ldots, n.
 $$
 Hence, 
 \begin{equation}\label{F*F}
 f(rV_1,\ldots, rV_n)^*=(I_\cG\otimes f(S_1,\ldots, S_n)^*)|\cK.
 \end{equation}
 Since $\text{Re}\, f(S_1,\ldots, S_n)\geq 0$, relation \eqref{F*F} implies 
 $$
 f(rV_1,\ldots, rV_n)^*+f(rV_1,\ldots, rV_n)\geq 0
 $$
 for $0<r<1$. Since 
 $$
 \lim_{r\to 1} f(rV_1,\ldots, rV_n)=f(V_1,\ldots, V_n)
 $$ 
 in the operator norm, we infer that
 $\text{Re}\, f(V_1,\ldots, V_n)\geq 0$.
 Hence and using  relation \eqref{ffpf2},  we obtain
 \begin{equation*}
 \text{Re}\,\left< f(T_1,\ldots, T_n)h,h\right> \geq  
 -\text{Re}\,\left< f(0)h,h\right>.
 \end{equation*}
 for any $h\in \cH$. The proof is complete.
 \end{proof}

 Another consequence of Theorem \ref{BKS1}, is the following multivariable
 power inequality.
 \begin{corollary}\label{power-ine}
 If $T_1,\ldots, T_n\in B(\cH)$, 
 then
 \begin{equation} \label{power-ineq}
 w(T_\alpha:\ |\alpha|=k)\leq w(T_1,\ldots, T_n)^k
 \end{equation}
 for any $k=1,2,\ldots $.
 \end{corollary}
 \begin{proof}
 Since the joint numerical radius is homogeneous, we can assume that
  $w(T_1,\ldots, T_n)=1$. Applying Theorem \ref{BKS1} to $(S_\alpha:\ |\alpha|=k)$,
  we obtain
  $$
  w(T_\alpha:\ |\alpha|=k)\leq \|[S_\alpha:\ |\alpha|=k]\|=1,
  $$
  which completes the proof.
 \end{proof}
 
We remark that if $n=1$ we obtain the classical 
power inequality \cite{B} (originally a conjecture of Halmos). Notice also that if $(A_1,\ldots, A_n)$ and $(B_1,\ldots, B_n)$ are arbitrary
 $n$-tuples
of operators then 
$$
w(A_iB_j:\ i,j=1,\ldots, n)\leq 4w(A_1,\ldots, A_n) w(B_1,\ldots, B_n).
$$
Indeed, using Theorem \ref{propri}, we deduce that
\begin{equation}
\begin{split}
w(A_iB_j:\ i,j=1,\ldots, n)&\leq \|[A_iB_j:\ i,j=1,\ldots, n]\|\\
&\leq 
\|[A_1,\ldots, A_n]\|\|[B_1,\ldots, B_n]\|\\
&\leq 2\|[A_1,\ldots, A_n]\| w(B_1,\ldots, B_n)\\
&\leq 4w(A_1,\ldots, A_n) 
w(B_1,\ldots, B_n).
\end{split}
\end{equation}

According to Theorem \ref{propri}, we have \begin{equation}\label{rwn}
r(T_1,\ldots, T_n)\leq w(T_1,\ldots, T_n)\leq\|[T_1,\ldots, T_n]\|.
\end{equation}In what follows we  characterize the $n$-tuples of operators for which 
equalities in \eqref{rwn} occur.

\begin{proposition}

If $(T_1,\ldots, T_n)\in B(\cH)^{(n)}$, then the following statements are equivalent:

\begin{enumerate}

\item[(i)] $r(T_1,\ldots, T_n)=\|[T_1,\ldots, T_n]\|;$

\item[(ii)]  $\|[T_\alpha:\ |\alpha|=k]\|=\|[T_1,\ldots, T_n]\|^k$ 
for any $k=1,2,\ldots;$

\item[(iii)] $w(T_1,\ldots, T_n)=\|[T_1,\ldots, T_n]\|.$

\end{enumerate}

\end{proposition}

\begin{proof}

Using Corollary \ref{pro}, it is easy to deduce that \begin{equation}\label{rr}
r(T_1,\ldots, T_n)^k=r(T_\alpha:\ |\alpha|=k)\quad \text{ for } k=1,2,\ldots.
\end{equation}
Consequently, if condition (i) holds, then $$\|[T_1,\ldots, T_n]\|^k=r(T_1,\ldots, T_n)^k=r(T_\alpha:\ |\alpha|=k)\leq 
\|[T_\alpha:\ |\alpha|=k]\|.
$$Since the reverse inequality is always true, we deduce (ii). The implication (ii)$\implies$(i) follows immediately  using the definition of the joint spectral radius. The inequality \eqref{rwn} shows that (i)$\implies$(iii). It remains to prove that (i)$\implies$(iii). Assume that $w(T_1,\ldots, T_n)=\|[T_1,\ldots, T_n]\|=1.$ Using Corollary \ref{pro}, we get

$$ w\left(\sum_{i=1}^n S_i\otimes T_i^*\right)=\left\|\sum_{i=1}^n S_i\otimes T_i^*\right\|=1.$$Hence, there is  a sequence $\{y_m\}\subset F^2(H_n)\otimes \cH$ such that $\|y_m\|=1$  and  $\left<\left(\sum_{i=1}^n S_i\otimes T_i^*\right)y_m,y_m\right>$ converges to 1, as $m\to\infty$. A standard argument shows that $1$ is in the approximate spectrum of $\sum_{i=1}^n S_i\otimes T_i^*$. Consequently, $r\left(\sum_{i=1}^n S_i\otimes T_i^*\right)=1$.
Using again Corollary \ref{pro}, the result follows.\end{proof}

\begin{proposition}

If $(T_1,\ldots, T_n)\in B(\cH)^{(n)}$, then   $w(T_1,\ldots, T_n)= r(T_1,\ldots, T_n)$ if and only if  $$w(T_\alpha:\ |\alpha|=k)=w(T_1,\ldots, T_n)^k\quad  
\text{ for  any } \ k=1,2,\ldots.$$

\end{proposition}

\begin{proof} Assume that $w(T_1,\ldots, T_n)= r(T_1,\ldots, T_n)$. Using relations \eqref{rr} and \eqref{rwn}, we have $$ w(T_1,\ldots, T_n)^k= r(T_1,\ldots, T_n)^k=r(T_\alpha:\ |\alpha|=k)\leq w(T_\alpha:\ |\alpha|=k).$$
Hence and using the power inequality \eqref{power-ineq}, we deduce that $w(T_1,\ldots, T_n)^k=w(T_\alpha:\ |\alpha|=k)$ for any $k=1,2,\ldots$. Conversely, assume that the latter equality holds. Using again \eqref{rwn}, we have$$w(T_1,\ldots, T_n)=w(T_\alpha:\ |\alpha|=k)^{1/k}\leq \|[T_\alpha:\ |\alpha|=k]\|^{1/k}$$ for any $k=1, 2, \ldots $. Taking $k\to\infty$, we obtain  $w(T_1,\ldots, T_n)\leq r(T_1,\ldots, T_n)$. The reverse inequality is always true (see \eqref{rwn}). This completes the proof.\end{proof}

 We consider now the case when the $n$-tuple $(T_1,\ldots, T_n)$ has commuting entries.

  \begin{theorem}\label{C2-commut}
  Let $T_1,\ldots, T_n\in B(\cH)$ be commuting operators and let $\cS_c\subset C^*(B_1,\ldots, B_n)$
  be the operator system defined by
  \begin{equation}\label{Sys-c}
  \cS_c:=\{p(B_1,\ldots, B_n)+q(B_1,\ldots, B_n)^*:\ p,q\in \cP\},
  \end{equation}
  where $B_1,\ldots, B_n$ are the creation operators on the symmetric Fock space.
  Then the following statements are equivalent:
  \begin{enumerate}
 \item[(i)] 
  $w(T_1,\ldots, T_n)\leq 1$;
 \item[(ii)] The map $\Psi_c:\cS_c\to B(\cH)$ defined by
 \begin{equation*} 
 \Psi_c(p(B_1,\ldots, B_n)+q(B_1,\ldots, B_n)^*):=
 p(T_1,\ldots, T_n)+q(T_1,\ldots,T_n)^*
 + (p(0)+\overline{q(0)})I
 \end{equation*}
 is completely positive.
 \item[(iii)] There is a Hilbert space $\cG\supseteq \cH$ and a $*$-representation
    $\pi:C^*(B_1,\ldots, B_n)\to B(\cG)$ such that 
    $$T_\alpha =2 P_\cH \pi(B_\alpha)| \cH\quad \text{ for any }  
    \alpha\in \FF_n^+\backslash \{g_0\}.
    $$
  \end{enumerate}
  \end{theorem}
  \begin{proof}
 Assume that $w(T_1,\ldots, T_n)\leq 1$. According the the proof of Theorem \ref{C2}, 
  the multi-Toeplitz  kernel 
   $\frac{1} {2}K_{T,2}:\FF_n^+\times \FF_n^+\to B(\cH)$
 defined
   by  
   $$
   \frac{1} {2}K_{T,2}(\alpha, \beta):= 
   \begin{cases}
   \frac{1} {2}T_{\beta\backslash \alpha} &\text{ if } \beta>\alpha\\
   I  &\text{ if } \alpha=\beta\\
    \frac{1} {2}T_{\alpha\backslash \beta}^*  &\text{ if } \alpha>\beta\\
    0\quad &\text{ otherwise}.
   \end{cases} 
   $$
  is positive definite.  Since $T_1,\ldots, T_n$ are commuting operators, 
  according to Theorem 3.3 of \cite{Po-moment}, there exists a completely positive linear map
  $\Phi:C^*(B,\ldots, B_n)\to B(\cH)$ such that $\Phi(I)=I$ and
  $$
  \Phi(B_\sigma)=\frac {1} {2} T_\sigma, \quad \sigma\in \FF_n^+\backslash\{g_0\}.
  $$
  According to Stinespring's theorem (see \cite{St}), there exists a Hilbert space $\cG\supseteq \cH$ 
  and a $*$-representation $\pi:C^*(B_1,\ldots, B_n)\to B(\cH)$ such  that (iii) holds.
  
  Assume now (iii) and let $W_i:=\pi(B_i)$, $i=1,\ldots, n$. Since  $[W_1,\ldots, W_n]$
  is a row contraction, there exists a Hilbert space $\cK\supseteq \cG$ and $V_1,\ldots, V_n\in B(\cK)$
  isometries with orthogonal ranges such that
  $W_\alpha^*=V_\alpha^*|\cK$, $\alpha\in \FF_n^+$. Since $\cH\subseteq \cG$, (iii) implies
  $T_\alpha=2P_\cH V_\alpha|\cH$, $\alpha\in \FF_n^+$. Using Theorem \ref{berger}, we deduce that
   $w(T_1,\ldots, T_n)\leq 1$,
  therefore (iii)$\implies$(i). 
  Now we show that 
(iii)$\implies$(ii).
  Notice  that (iii)  implies
  \begin{equation*}
  \begin{split}
  p(T_1,\ldots, T_n)+ q(T_1,\ldots, T_n)^*&+(p(0)+
  \overline{q(0)})I\\
  &=2P_\cH\pi[p(B_1, \ldots, B_n)+
  q(B_1, \ldots, B_n)^*]|\cH.
  \end{split}
  \end{equation*}
 This clearly implies (ii).
 To prove that  (ii)$\implies$(i), note that there is a completely positive linear map
 $\gamma:C^*(S_1,\ldots, S_n)\to C^*(B_1,\ldots, B_n)$ such that 
 $\gamma(S_\alpha S_\beta^*)=B_\alpha B_\beta^*$
 (see \cite{Po-poisson}).
 Since the map $\Psi_c\circ \gamma$ is completely positive and 
 $(\Psi_c\circ \gamma)|\cS=\Psi$,  where $\cS$ and $\Psi$ are defined in Theorem \ref{C2},
 the latter theorem  implies $w(T_1,\ldots, T_n)\leq 1$.
 The proof is complete.
  \end{proof}

 \begin{remark} If $T_1,\ldots, T_n$ are commuting operators, then  ``commutative''
  versions
 of  Corollary $\ref{simi}$, Theorem $\ref{BKS1}$, Corollary $\ref{BKS}$, and 
 Theorem $\ref{Real-ine}$  hold true, if we replace the noncommutative disc algebra
  $\cA_n$ with its commutative 
 version $\cA_n^c$. The proofs are exactly the same but one uses Theorem $\ref{C2-commut}$.
 \end{remark}

 \smallskip

 \section{Euclidean  operator radius}
 \label{euclid}

In this section  we present basic properties of the {\it euclidean operator radius} 
of an $n$-tuple of operators $(T_1,\ldots, T_n)$, defined by\begin{equation*}
 w_e(T_1,\ldots, T_n):= \sup_{\|h\|=1}\left(\sum_{i=1}^n\left|\left<
 T_ih,h\right>\right|^2\right)^{1/2},
 \end{equation*}in connection with the joint numerical radius and several other operator radii.   
We define a new norm  and   ``spectral radius''  on $B(\cH)^{(n)}$
   by setting
 \begin{equation}
 \label{e-norm}
 \|(T_1,\ldots, T_n)\|_e:=\sup\limits_{(\lambda_1,\ldots, \lambda_n)\in \BB_n} 
 \|\lambda_1  T_1+\cdots + \lambda_n  T_n\|
 \end{equation}
  and
 \begin{equation}
 \label{e-radius}
 r_e(T_1,\ldots, T_n):=\sup\limits_{(\lambda_1,\ldots, \lambda_n)\in \BB_n} 
 r(\lambda_1  T_1+\cdots + \lambda_n  T_n),
 \end{equation}
 where $r(X)$ denotes the  usual spectral radius of an operator $X\in B(\cH)$.
Notice that $\|\cdot\|_e$ is a norm on   $B(\cH)^{(n)}$,    
 $$
 \|(T_1,\ldots, T_n)\|_e=\|(T_1^*,\ldots, T_n^*)\|_e,\quad \text{ and }\quad
 r_e(T_1,\ldots, T_n)=r_e(T_1^*,\ldots, T_n^*).
 $$

  In what follows we show that $\|\cdot\|_e$ is equivalent to the operator norm on 
$B(\cH)^{(n)}$.
 \begin{theorem}\label{norm-ine}    
 If $(T_1,\ldots, T_n)\in B(\cH)^{(n)}$, then 
  \begin{equation}\label{radical}
 \frac {1} {\sqrt{n}}\|[T_1,\ldots, T_n]\|\leq 
 \|(T_1,\ldots, T_n)\|_e\leq \|[T_1,\ldots, T_n]\|,
 \end{equation}
 where the  constants $\frac {1} {\sqrt{n}}$ and $1$ are the best possible, and
 \begin{equation}\label{re}
 r_e(T_1,\ldots, T_n)\leq r(T_1,\ldots, T_n).
 \end{equation}
 \end{theorem}
 \begin{proof} 
Let $\sigma$ be the rotation-invariant normalized positive 
 Borel measure on the unit sphere $\partial \BB_n$.   Using  the  relations 
 (see \cite{R}) 
 $$
 \int_{\partial \BB_n}|\lambda_i|^2d\sigma(\lambda)=\frac{1} {n}
  \quad \text{\rm and }\quad \int_{\partial \BB_n}\lambda_i 
  \overline{\lambda}_j d\sigma(\lambda)=0 \quad  \text{\rm  if }\  i\neq j,\ i,j=1,\ldots, n,
  $$
 we  deduce that
 \begin{equation*}
 \begin{split}
 \|(T_1,\ldots, T_n)\|^2_e &=
 \sup_{(\lambda_1,\ldots, \lambda_n)\in \BB_n} \sup_{\|h\|=1}\left<
 \left(\sum_{i=1}^n \lambda_i T_i\right)
 \left(\sum_{j=1}^n \overline{\lambda}_j T_j^*\right)h,h\right>\\
 &\geq
  \sup_{\|h\|=1} \int_{\partial \BB_n}\sum_{i,j=1}^n \lambda_i 
  \overline{\lambda}_j\left<T_iT_j^*h,h\right>d\sigma(\lambda)
  \\
  &=\sup_{\|h\|=1}\frac{1} {n} \left<\sum_{i=1}^nT_iT_i^*h,h\right>
  \\
  &=\frac{1} {n}\|[T_1,\ldots, T_n]\|^2,
 \end{split}
 \end{equation*}
On the other hand,  we have
 \begin{equation*}
 \begin{split}
 \|(T_1,\ldots, T_n)\|_e &=
 \sup_{(\lambda_1,\ldots, \lambda_n)\in \BB_n}\|\lambda_1T_1+\cdots 
 \lambda_n T_n\|\\
 &\leq \sup_{(\lambda_1,\ldots, \lambda_n)\in \BB_n}
 \left( \sum_{i=1}^n |\lambda_i|^2\right)^{1/2}\left\|\sum_{i=1}^n
  T_i T_i^*\right\|^{1/2}\\
 &=\|[T_1,\ldots, T_n]\|.
 \end{split}
 \end{equation*} 
  Now, notice that if $S_1,\ldots, S_n$ are the left creation operators, then
 \begin{equation}\begin{split}
 1&\leq\frac{1}{\sqrt{n}}\|[S_1^*,\ldots, S_n^*]\|\leq\|[S_1^*,\ldots, S_n^*]\|_e\\
 &=\|[S_1,\ldots, S_n]\|_e\leq \|[S_1,\ldots, S_n]\|=1.
 \end{split}
 \end{equation}
 This shows that the inequalities \eqref{radical} are best possible.
 
  To prove inequality \eqref{re},
  notice that
 \begin{equation*}
 \begin{split}
 r_e(T_1,\ldots, T_n)&=
 \sup_{(\lambda_1,\ldots, \lambda_n)\in \BB_n}\inf_{m\in \NN}
 \left\| \left( \sum_{i=1}^n \lambda_i T_i\right)^m\right\|^{1/m}\\
 &\leq \inf_{m\in \NN}\sup_{(\lambda_1,\ldots, \lambda_n)\in \BB_n}
 \left\|\left( \sum_{i=1}^n \overline{\lambda}_i T_i^*\right)^m\right\|^{1/m}\\
 &\leq 
 \inf_{m\in \NN}\sup_{(\lambda_1,\ldots, \lambda_n)\in \BB_n}
 \left[\left(\sum_{|\alpha|=m} |\lambda_\alpha|^2\right)^{1/2m}
 \left\|\sum_{|\alpha|=m}T_\alpha T_\alpha^*\right\|^{1/2m}\right]\\
 &\leq \inf_{m\in \NN} \left\|\sum_{|\alpha|=m}T_\alpha 
 T_\alpha^*\right\|^{1/2m}\\
 &=r(T_1,\ldots, T_n).
 \end{split}
 \end{equation*}
  Notice that the inequality \eqref{re} is strict in general.
 For example, we have
$r(S_1^*,\ldots, S_n^*)=\sqrt{n}$ and 
 $r_e(S_1^*,\ldots, S_n^*)=r_e(S_1,\ldots, S_n)=1$.
 Therefore, $r_e(S_1^*,\ldots, S_n^*)<r(S_1^*,\ldots, S_n^*)$ if $n\geq 2$.
 On the other hand, we have equality in \eqref{re} if $T_i=S_i$, $i=1,\ldots, n$.
 Indeed,  
  $
  r_e(S_1,\ldots, S_n)=r(S_1,\ldots, S_n)=1.
  $
 This completes the proof.
 \end{proof}

The next result summarizes some of the basic properties of 
the euclidean operator radius of an $n$-tuple of operators.

  \begin{theorem}
  \label{propri2}
 The  euclidean operator  radius $w_e:B(\cH)^{(n)}\to [0,\infty)$   
 for $n$-tuples of operators  satisfies the following properties:
 \begin{enumerate}
 \item[(i)]
 $w_e(T_1,\ldots, T_n)=0$ if and only if $T_1=\cdots=T_n=0$;
 \item[(ii)]
 $w_e(\lambda T_1,\ldots, \lambda T_n)=|\lambda|~w_e(T_1,\ldots, T_n) $
 for any $\lambda\in \CC$;
 \item[(iii)]
 $w_e(T_1+T_1',\ldots, T_n+T_n')\leq w_e(T_1,\ldots, T_n)+w_e(T_1',\ldots, T_n');$
 \item[(iv)]
 $w_e(U^*T_1U,\ldots, U^*T_nU)= w_e(T_1,\ldots, T_n)$ for any unitary operator
  $U:\cK\to \cH$;
  \item [(v)]
  $w_e(X^*T_1X,\ldots, X^*T_n X)\leq  \|X\|^2 w_e(T_1,\ldots, T_n)$ for any 
    operator $X:\cK\to \cH$;
 \item[(vi)]
 $\frac{1}{2} \| (T_1,\ldots, T_n)\|_e\leq
 w_e(T_1,\ldots, T_n)\leq   \| (T_1,\ldots, T_n)\|_e;$
 \item[(vii)]
 $r_e(T_1,\ldots, T_n)\leq w_e(T_1,\ldots, T_n);$
 
 \item [(viii)]
 $w_e(I_\cE\otimes T_1,\ldots, I_\cE\otimes T_n)=w_e(T_1,\ldots, T_n)$ for 
 any separable
 Hilbert space $\cE$;
 \item[(ix)]
 $w_e$ is a continuous map in the norm topology.
 \end{enumerate}
  \end{theorem}
  \begin{proof}
  The first five properties can be easily deduced  using the definition of $w_e$.
 Now,  notice that
  \begin{equation*}
  \begin{split}
  w_e(T_1,\ldots, T_n)&=
  \sup_{\|h\|=1}\left(\sum_{i=1}^n |\left< T_i h,h\right>|^2\right)^{1/2}\\
  &=\sup_{\|h\|=1}\sup_{(\lambda_1,\ldots, \lambda_n)\in \BB_n} 
  \left|\sum_{i=1}^n \lambda_i\left< T_i h,h\right>\right|\\
  &=\sup_{(\lambda_1,\ldots, \lambda_n)\in \BB_n} \sup_{\|h\|=1}
  \left|\left< \sum_{i=1}^n \lambda_i T_i h,h\right>\right|\\
  &=\sup_{(\lambda_1,\ldots, \lambda_n)\in \BB_n}
  w(\lambda_1  T_1+\cdots + \lambda_n  T_n).
  \end{split}
  \end{equation*}
  Consequently, we obtain
  \begin{equation}\label{we-sup}
  w_e(T_1,\ldots, T_n)=\sup\limits_{(\lambda_1,\ldots, \lambda_n)\in \BB_n} 
 w(\lambda_1  T_1+\cdots + \lambda_n  T_n).
  \end{equation}
  It is well known that 
  $$\frac {1} {2} \|X\|\leq w(X)\leq \|X\|\quad \text{ and } \quad r(X)\leq \|X\|
  $$
  for any $X\in B(\cH)$. Applying these inequalities to the operator 
  $X:=\lambda_1 T_1+\cdots + \lambda_n T_n$ for
   $(\lambda_1,\ldots, \lambda_n)\in \BB_n$, and using relation \eqref{we-sup}, we 
   deduce (vi) and (vii).
   
  To prove (viii),  we use relation \eqref{we-sup} and the fact that the
   classical numerical radius
   satisfies the equation $w(I_\cG\otimes X)=w(X)$.
   Indeed, we have
   \begin{equation*}
   \begin{split}
   w_e(I_\cE\otimes T_1,\ldots, I_\cE\otimes T_n)&=
   \sup_{(\lambda_1,\ldots, \lambda_n)\in \BB_n}
   w\left(I_\cE\otimes \sum_{i=1}^n \lambda_iT_i \right)\\
   &=\sup_{(\lambda_1,\ldots, \lambda_n)\in \BB_n}
   w\left(\sum_{i=1}^n \lambda_iT_i\right)\\
  &=
   w_e(T_1,\ldots, T_n).
   \end{split}
   \end{equation*}
   According to (vi) and Theorem \ref{norm-ine}, we  obtain
   $$
  w_e(T_1,\ldots, T_n)\leq \left\| \sum_{i=1}^n T_i T_i^*\right\|^{1/2}.
  $$
  Hence, we deduce that $w_e$ is continuous in the norm topology.
   The proof is complete.
  \end{proof}
  
  \begin{corollary}\label{we-sup2}
  If $T_1,\ldots, T_n\in B(\cH)^{(n)}$, then 
 $$
 w_e(T_1,\ldots, T_n)=\sup\limits_{(\lambda_1,\ldots, \lambda_n)\in \BB_n} 
 w(\lambda_1  T_1+\cdots + \lambda_n  T_n).
 $$
  \end{corollary}

  The next result shows that the euclidean operator radius is equivalent 
  to the operator norm and the joint numerical radius on $B(\cH)^{(n)}$.
  \begin{proposition}\label{2-ineq}
  If $T_1,\ldots, T_n\in B(\cH)^{(n)}$, then 
  \begin{equation}\label{-wew}
  w_e(T_1,\ldots, T_n)\leq w(T_1,\ldots, T_n)
\text{ and }
  \end{equation}
  \begin{equation}\label{-wew2}
  \frac {1} {2\sqrt{n}} \|[T_1,\ldots, T_n]\|\leq  w_e(T_1,\ldots, T_n)
  \leq \|[T_1,\ldots, T_n]\|.
  \end{equation}
  Moreover, the inequalities are sharp.
  \end{proposition}
  \begin{proof}
  Given $\lambda:=(\lambda_1,\ldots \lambda_n)\in \BB_n$, define the vector
 $z_\lambda\in F^2(H_n)$ by setting
 \begin{equation}\label{zl}
 z_\lambda:=\frac{1}{ \sqrt{1- \|\lambda\|^2}}
 \left( \sum_{\alpha\in \FF_n^+}
 \overline{\lambda}_\alpha e_\alpha \right).
 \end{equation}
  We remark that $S_i^*z_\lambda =\overline{\lambda}_i
  z_\lambda$ for any $i=1,\ldots, n$, and $\|z_\lambda\|=1$.
 For any vector  $h\in \cH$, $\|h\|=1$,  and $\lambda\in \BB_n$,  we have
 \begin{equation*}
 \begin{split}
 \left< \sum_{i=1}^n ( S_i\otimes T_i^*)
 (z_\lambda\otimes h), z_\lambda \otimes h \right>&
 = \sum_{i=1}^n \left< z_\lambda \otimes h, \overline{\lambda}_i z_\lambda 
 \otimes T_i h\right>\\
 &=\sum_{i=1}^n \lambda_i \left< h, T_i h\right>.
 \end{split}
 \end{equation*} 
 Since $\|z_\lambda\otimes h\|=1$ and 
 $w(T_1,\ldots, T_n)=w(S_1\otimes T^*_1+\cdots +S_n\otimes T^*_n)$, we infer that
 $$
 \sup_{\lambda\in B_n}\left|\sum_{i=1}^n \lambda_i \left< h, T_i h\right>\right|\leq 
 w(T_1,\ldots, T_n).
 $$
   Hence, we deduce that
  $$
  \left(\sum_{i=1}^n |\left< T_i h, h\right>|^2\right)^{1/2}\leq 
  w(T_1,\ldots, T_n).
  $$
  Taking now the supremum over all $h\in \cH$ with $\|h\|=1$, 
   we obtain  inequality \eqref{-wew}.
   Note that we have equality  in  \eqref{-wew} if $T_i=S_i$, $i=1,\ldots, n$. 
    Indeed
   $w_e(S_1,\ldots, S_n)=w(S_1,\ldots, S_n)=1$ (see Theorem \ref{example}).
   On the other hand,
    the inequality
\eqref{-wew} is strict
   if $T_i=S_i^*$, $i=1,\ldots, n$.  More precisely, we have
   $$\sqrt{n}=r(S_1^*,\ldots, S_n^*)\leq 
   w(S_1^*,\ldots, S_n^*)>w_e(S_1^*,\ldots, S_n^*)=w_e(S_1,\ldots, S_n)=1.
   $$
   
   Combining inequality (vi) of Theorem \ref{propri2} with inequality \eqref{radical},
we obtain \eqref{-wew2}. Notice also that 
$$w_e(S_1,\ldots, S_n)=\|[S_1,\ldots, S_n]\|=1.
 $$
    On the other hand, here is  an example when we have equality 
     in the first inequality of
   \eqref{-wew2}.
   Take $T_i:=\frac{1}{\sqrt{n}}S^{(2)}$, $i=1,\ldots, n$,  where
    $S^{(2)}:=\left[\begin{matrix}0&0\\
   1&0
   \end{matrix}\right].
   $
   It is easy to see that $\|S^{(2)}\|=1$ and $w(S^{(2)})=\frac{1} {2}$.
   We also have
   $$
   \frac{1} {2\sqrt{n}}\|[T_1,\ldots, T_n]\|=\frac{1} {2}\|S^{(2)}\|=\frac{1} {2}.
   $$
   and
   $$
   w_e(T_1,\ldots, T_n)=\frac{1}{\sqrt{n}} w_e(S^{(2)}, \ldots, S^{(2)})=
   w(S^{(2)})=\frac{1} {2}.
   $$
   This completes the proof of the theorem.
  \end{proof}

Combining inequalities \eqref{-wew}, \eqref{-wew2}, and using Theorem \ref{propri} part (vi), 
we deduce that the euclidean operator radius is equivalent to the joint numerical radius.
We obtain
$$\frac{1} {2\sqrt{n}}w(T_1,\ldots, T_n)\leq w_e(T_1,\ldots, T_n)\leq w(T_1,\ldots, T_n).$$
  
  Now, we can prove the following von Neumann type inequality
   for the numerical range.
 
 \begin{theorem}\label{von-2}
 If $T_1,\ldots, T_n\in B(\cH)$,
 then 
 $$
 w(p(T_1,\ldots, T_n))\leq \|p\|_2\left(\frac{1-w(T_1,\ldots, T_n)^{2(m+1)}}
 {1-w(T_1,\ldots, T_n)^2} \right)^{1/2}
 $$
 for any   polynomial $p=\sum_{|\alpha|\leq m} a_\alpha e_\alpha$ in $F^2(H_n)$. 
 \end{theorem}
 \begin{proof}
  According to 
Proposition \ref{2-ineq}, we have
\begin{equation} \label{sww}
\sup_{\|h\|=1} \left(\sum_{|\alpha|=k}
  \left|\left< T_\alpha h,h\right>\right|^2\right)^{1/2}=
  w_e(T_\alpha:\ |\alpha|=k)
  \leq w(T_\alpha:\ |\alpha|=k)
\end{equation}
for $1\leq k\leq m$.
    On the other hand, the joint numerical radius is a norm 
    (see  Theorem \ref{propri}). Consequently, using \eqref{sww},  the multivariable power inequality
    \eqref{power-ineq}, and 
 applying Cauchy's inequality twice,   we obtain
 \begin{equation*}
 \begin{split}
 w\left(p(T_1,\ldots, T_n) \right) &\leq \sum_{k=0}^{m}
 w\left( \sum_{|\alpha|=k} a_\alpha T_\alpha\right)\\
 &=\sum_{k=0}^{m}\sup_{\|h\|=1} \left|\sum_{|\alpha|=k}
  a_\alpha \left< T_\alpha h,h\right>\right|\\
 &\leq \sum_{k=0}^{m}\left[ \left(\sum_{|\alpha|=k} |a_\alpha|^2\right)^{1/2}
 \sup_{\|h\|=1} \left(\sum_{|\alpha|=k}
  \left|\left< T_\alpha h,h\right>\right|^2\right)^{1/2}\right]\\
  &\leq  
  \sum_{k=0}^{m}\left[\left(\sum_{|\alpha|=k} |a_\alpha|^2\right)^{1/2}
   w(T_\alpha:\ |\alpha|=k)\right]\\
  &\leq 
  \sum_{k=0}^{m}\left[\left(\sum_{|\alpha|=k} |a_\alpha|^2\right)^{1/2}
   w(T_1,\ldots, T_n)^k\right]\\
  &\leq 
  \left(\sum_{|\alpha|\leq m} |a_\alpha|^2\right)^{1/2}
 \left(\sum_{k=0}^{m} w(T_1,\ldots, T_n)^{2k} 
 \right)^{1/2}.
 \end{split}
 \end{equation*}
  The proof is complete.
\end{proof}

  Theorem \ref{von-2} and the  $\cA_n$-functional calculus for $n$-tuples of operators with joint numerical radius 
  $w(T_1,\ldots, T_n)<1$ (see Section 1) can be used  to deduce the following result.
 
 \begin{corollary}\label{BKS-2}
 If $T_1,\ldots, T_n\in B(\cH)$, $w(T_1,\ldots, T_n)<1$, and $f(S_1,\ldots, S_n)\in \cA_n$, then
 $$
 w(f(T_1,\ldots, T_n))\leq \frac{1} {\sqrt{1-w(T_1,\ldots, T_n)^2}} \|f\|_2,
 $$
 where $\|f\|_2:=\|f(S_1,\ldots, S_n)(1)\|$.
 Moreover, if $f(0)=0$, then
 $$
 w(f(T_1,\ldots, T_n))\leq \frac{w(T_1,\ldots, T_n)}
  {\sqrt{1-w(T_1,\ldots, T_n)^2}} \|f\|_2.
 $$
 \end{corollary}
 Notice that these inequalities are quite different from the multivariable Berger-Kato-Stampfli
 type inequality obtained in Section 1 (see Corollary \ref{BKS}).

\smallskip

  \section{Joint numerical  range and  spectrum}
  \label{JNRS}

The joint (spatial) numerical range of  an $n$-tuple of operators
 $(T_1,\ldots, T_n)\in B(\cH)^{(n)}$   is the subset of $\CC^n$ defined by
 $$
 W(T_1,\ldots, T_n):=\{(\langle T_1h, h\rangle,\ldots, \langle T_nh, h\rangle):
 \   h\in \cH, ~\|h\|=1\}.
 $$
There is a large literature relating to the classical numerical 
 range ($n=1$). For general proprieties, we refer to \cite{Ha},
 \cite{BD1}, \cite{BD2},  and \cite{GR}. 
 In this section we study
 the joint numerical range  of an $n$-tuple of operators $(T_1,\ldots, T_n)$
 in connection with unitary invariants such as  the {\it right spectrum} $\sigma_r(T_1,\ldots, T_n)$,
  the joint numerical radius $w(T_1,\ldots, T_n)$, the euclidean operator 
  radius $w_e(T_1,\ldots, T_n)$, and the {\it joint spectral radius} $r(T_1,\ldots, T_n)$.
On the other hand, we  obtain 
an analogue of Toeplitz-Hausdorff theorem (\cite{T}, \cite{Hau}) on the convexity 
 of the spatial numerical range of an operator on a Hilbert space, 
 for the joint numerical range
 of operators in the noncommutative analytic Toeplitz algebra $F_n^\infty$.

 We recall   that the joint  right spectrum 
  $\sigma_r(T_1,\ldots, T_n)$ of an   $n$-tuple  
 $(T_1,\ldots, T_n)$ of operators
   in $B(\cH)$ is the set of all $n$-tuples
    $(\lambda_1,\ldots, \lambda_n)$  of complex numbers such that the
     right ideal of $B(\cH)$  generated by the operators
     $\lambda_1I-T_1,\ldots, \lambda_nI-T_n$ does
      not contain the identity operator. Notice that 
the joint left spectrum
  $\sigma_l(T_1^*,\ldots, T_n^*)$ is the complex conjugate
   of $\sigma_r(T_1,\ldots, T_n)$.
  We should mention the connection between the joint   right spectrum 
   of an   $n$-tuple  
 $(T_1,\ldots, T_n)$   and the solutions of  linear equations of the form
 \begin{equation}
 \label{lin-eq}
 T_1x_1+T_2x_2+\cdots +T_nx_n=y,
 \end{equation}
 where $y\in \cH$ is a given vector.  Notice that 
 if $(0,\ldots, 0)\notin 
 \sigma_r(T_1,\ldots, T_n)$ then the equation \eqref{lin-eq}
  has solutions in $\cH$.

 For each $r>0$, define 
 $$
 (\CC^n)_r:=\{(\lambda_1,\ldots, \lambda_n)\in \CC^n:
 \ |\lambda_1|^2+\cdots +|\lambda_n|^2<r^2\}
 $$
 and let $\BB_n:=(\CC^n)_1$ be the open unit ball of $\CC^n$.

 \begin{theorem}\label{conv}
 If $T_1,\ldots, T_n\in B(\cH)$, then
 \begin{enumerate}
 \item[(i)]
  $ \sigma_r(T_1,\ldots, T_n)\subseteq \overline{W(T_1,\ldots, T_n)} \subseteq 
 \overline{(\CC^n)}_{w_e(T_1,\ldots, T_n)}\subseteq 
 \overline{(\CC^n)}_{w(T_1,\ldots, T_n)}
 \subseteq 
 \overline{(\CC^n)}_{\|[T_1,\ldots, T_n]\|};$
 \item[(ii)]
 $\sigma_r(T_1,\ldots, T_n)\subseteq 
 \overline{(\CC^n)}_{r(T_1,\ldots, T_n)}\subseteq 
 \overline{(\CC^n)}_{w(T_1,\ldots, T_n)};$
 \item[(iii)] $ \sigma_r(T_1,\ldots, T_n)$ is a compact subset of $\CC^n$.
 \end{enumerate}
 \end{theorem}
 
 \begin{proof}
 First we prove that 
 $(\lambda_1, \ldots, \lambda_n)\notin \sigma_r(T_1,\ldots, T_n)$ if and only if 
 there exists $\delta>0$ such that 
 \begin{equation}\label{dine}
 \sum_{i=1}^n\|(\overline{\lambda}_i I-T_i^*)h\|^2\geq \delta
  \|h\|^2\quad \text{ for any  } h\in \cH.
 \end{equation}
 Assume that $(\lambda_1, \ldots, \lambda_n)\notin \sigma_r(T_1,\ldots, T_n)$. Then there 
 exist  some operators $X_i\in B(\cH)$, $i=1,\ldots, n$, such that 
 $\sum_{i=1}^n X_i(\overline{\lambda}_i I-T_i^*)=I$.
 Therefore, at least one of the operators $X_i$ is different from zero and 
 $$
 \|h\|\leq \|[X_1,\ldots, X_n]\| 
 \left( \sum_{i=1}^n\|(\overline{\lambda}_i I-T_i^*)h\|^2\right)^{1/2},
 $$
 which proves \eqref{dine}. Conversely, assume that  inequality 
 \eqref{dine} holds for some $\delta>0$.  Then  $AA^*\geq CC^*$, where $A:=[\lambda_1I-T_1,\ldots,
 \lambda_nI-T_n]$ and $C:=\sqrt{\delta} I$. Applying Douglas 
 factorization theorem \cite{Do}, we find a contraction $X:\overline{A^*\cH}\to \cH$ such that 
 $XA^*=C^*$. Extend $X$ to an operator $\tilde{ X}:\oplus_{i=1}^n \cH\to \cH$
 by setting $\tilde{X}k=0$ if $k\in (A^*\cH)^\perp$. Since $\tilde{X}=[\tilde{X}_1, 
 \ldots, \tilde{X}_n]$ for some operators $\tilde{X}_i\in B(\cH)$,
  the equation $\tilde{X}A^*=C^*$ implies $\sum\limits_{i=1}^n (\lambda_i I-T_i)\tilde{X}_i^*=
  \sqrt{\delta} I$,
  which shows that 
   $(\lambda_1, \ldots, \lambda_n)\notin \sigma_r(T_1,\ldots, T_n)$.
   
   Now let us show that 
   $\sigma_r(T_1,\ldots, T_n)\subseteq \overline{W(T_1,\ldots, T_n)}$.
   Notice that if $(\lambda_1, \ldots, \lambda_n)\in \sigma_r(T_1,\ldots, T_n)$, 
   then \eqref{dine}
   implies that there exists a sequence  of vectors $\{h_m\}_{m=1}^\infty$  with 
   $\|h_m\|=1$ and such that 
   $$
   \lim_{m\to \infty} 
   \sum_{i=1}^n\|(\overline{\lambda}_i I-T_i^*)h_m\|^2 =0.
    $$
    Hence, we deduce that
    $\lim\limits_{m\to \infty}\left< (\overline{\lambda}_i I-T_i^*)h_m, h_m\right>=0$
    and therefore 
    $\lim\limits_{m\to \infty}\left< T_ih_m, h_m\right>=\lambda_i$ for any $i=1,\ldots, n$.
    This shows that $(\lambda_1, \ldots, \lambda_n)\in \overline 
    {W(T_1,\ldots, T_n)}$.
    
   The second inclusion  in (i) is is due to the definition of $w_e$, and the
    third (resp. fourth) inclusion  is due to Proposition \ref{2-ineq}
    (resp.  
    Theorem \ref{propri} part (vii)).
 Now let us prove  the inclusions of part (ii).
 Let $(\lambda_1,\ldots, \lambda_n)\in \sigma_r(T_1,\ldots, T_n)$.
 Since the left ideal of $B(\cH)$ generated by the operators 
 $T_1^*-\overline{\lambda}_1 I, \ldots,
 T_n^*-\overline{\lambda}_n I$ does not contain the identity, there is a pure
  state 
 $\mu$ on $B(\cH)$ such that $\mu(X(T_i^*-\overline{\lambda}_i I))=0$ for 
 any $X\in B(\cH)$
 and $i=1,\ldots, n$. In particular, we have $\mu(T_i)= \lambda_i=\overline{
 \mu(T_i^*)}$ and
 $$
 \mu(T_\alpha T_\alpha^*)=\overline{\lambda}_\alpha
  \mu(T_\alpha)=|\lambda_\alpha|^2,\qquad \alpha\in \FF_n^+.
 $$
 Hence, we infer that
 \begin{equation*}
 \begin{split}
 \left( \sum_{i=1}^n |\lambda_i|^2 \right)^{1/2}&=
 \left( \sum_{|\alpha|=m} |\lambda_\alpha|^2\right)^{1/2m}
 =\left|\mu \left(\sum_{|\alpha|=m} T_\alpha T_\alpha^*\right)\right|^{1/2m}\\
& \leq \left\| \sum_{|\alpha|=m} T_\alpha T_\alpha^*\right\|^{1/2m}
 \leq r(T_1,\ldots, T_n).
 \end{split}
 \end{equation*}
 The last inclusion in (ii) is due to Theorem \ref{propri} part (viii).
 To prove (iii),     it is enough
 to show that $\sigma_r(T_1,\ldots, T_n)$ is a closed 
 subset of $\CC^n$. 
 Let $\left\{(\lambda_1^{(m)}, \ldots, 
 \lambda_n^{(m)})\right\}_{m=1}^\infty$ be 
 a sequence of vectors in 
 $\sigma_r(T_1,\ldots, T_n)$ such that $\lim\limits_{m\to\infty}
  \lambda_i^{(m)}=\lambda_i$, $\lambda_i\in \CC$, for each $i=1,\ldots, n$.
  Assume that $(\lambda_1,\ldots, \lambda_n)\notin \sigma_r(T_1,\ldots, T_n)$ and
  choose $m$ large enough such that
  $$
  \gamma:=\delta-\left(\sum_{i=1}^n|\lambda_i^{(m)}-
  \lambda_i|^2\right)^{1/2}>0.
  $$
 Using inequality \eqref{dine}, we obtain
 \begin{equation*}
 \begin{split}
 \left(\sum_{i=1}^n\left\|(\overline{\lambda}_i^{(m)} I-T_i^*)h\right\|^2\right)^{1/2}
 &\geq 
  \left(\sum_{i=1}^n\|(\overline{\lambda}_i I-T_i^*)h\|^2\right)^{1/2}-
   \left(\sum_{i=1}^n | \lambda_i 
   -\lambda_i^{(m)} |^2\right)^{1/2}\|h\|\\
   &\geq \gamma\|h\|
 \end{split}
 \end{equation*}
 for any $h\in \cH$. This shows that 
 $(\lambda_1^{(m)}, \ldots, \lambda_n^{(m)})\notin \sigma_r(T_1,\ldots, T_n)$, which 
 is a contradiction.
 Therefore, we must have 
 $(\lambda_1,\ldots, \lambda_n)\in \sigma_r(T_1,\ldots, T_n)$.  
 The proof is complete.
 \end{proof}

From the proof of Theorem \ref{conv}, we can deduce the following characterization for the right spectrum of an $n$-tuple of operators.
 
 \begin{corollary}\label{approx}
 If $T_1,\ldots, T_n\in B(\cH)$, then the following statements
  are equivalent:
  \begin{enumerate}
  \item[(i)]
  $(\lambda_1,\ldots, \lambda_n)\notin \sigma_r(T_1,\ldots, T_n)$;
  \item[(ii)] there exists $\delta>0$ such that
  $\sum\limits_{i=1}^n (\lambda_iI-T_i)(\overline{\lambda}_iI-T_i^*)\geq \delta I;$
  \item[(iii)]
  $(\lambda_1I-T_1)B(\cH)+\cdots + (\lambda_nI-T_n)B(\cH)=B(\cH)$.
  \end{enumerate}
 \end{corollary}

Another consequence of Theorem \ref{conv} is the following result.

 \begin{corollary}\label{prop}
 If $(\lambda_1,\ldots, \lambda_n) \in  W(T_1,\ldots, T_n)$ and 
 $$
 \sum\limits_{i=1}^n |\lambda_i|^2=\left\|\sum_{i=1}^n T_iT_i^*\right\|,
 $$ then
  $(\lambda_1,\ldots, \lambda_n)\in \sigma_r(T_1,\ldots, T_n)$.
 \end{corollary}
 \begin{proof}
 Since  $(\lambda_1,\ldots, \lambda_n) \in  W(T_1,\ldots, T_n)$,
 there is $h\in \cH$, $\|h\|=1$, such that $\lambda_i=\left< T_ih,h\right>$ for
 $i=1,\ldots, n$. We have
 \begin{equation*}
 \begin{split}
 \left\|\sum_{i=1}^n T_iT_i^*\right\|&=
  \sum\limits_{i=1}^n |\lambda_i|^2=\sum_{i=1}^n |\left< T_ih,h\right>|^2\\
  &\leq \sum_{i=1}^n\|T_ih\|^2 \|h\|^2
  \leq \left<\sum_{i=1}^n T_iT_i^*h, h\right>\\
  &\leq \left\|\sum_{i=1}^n T_iT_i^*\right\|.
  \end{split}
 \end{equation*}
 Consequently, we  must have
 $$
 |\left< T_ih,h\right>|=\|T_ih\| \|h\|, \qquad i=1,\ldots, n
 $$
 which implies  $T_i^*h=\mu_ih$ for some $\mu_i\in \CC$, $i=1,\ldots, n$.
 On the other hand, we have
 $$
 \overline{\lambda}_i=\left< T_i^*h,h\right>=\left<\mu_ih,h\right>=\mu_i, \quad
 i=1,\ldots, n.
 $$
 Therefore, $T_i^*h=\overline{\lambda}_ih$, $i=1,\ldots, n$.
  Now, Corollary \ref{approx} implies 
  $(\lambda_1,\ldots, \lambda_n)\in \sigma_r(T_1,\ldots, T_n)$.
 The proof is complete.
 \end{proof}

 Notice that if  the Hilbert space $\cH$ is finite dimensional  and 
 $w_e (T_1,\ldots, T_n)=\|\sum\limits_{i=1}T_iT_i^*\|^{1/2}$, then Corollary
 \ref{prop} implies
  $\sigma_r(T_1,\ldots, T_n)\neq\emptyset$.

We remark that
 the spectral inclusion of Theorem \ref{conv} part (i) enables us to locate the 
 right spectrum of the sum of two $n$-tuples of operators.
 More precisely, we have
 \begin{equation*}
 \begin{split}
 \sigma_r(T_1+T_1',\ldots, T_n+T_n')&\subseteq 
 \overline{W(T_1+T_1',\ldots, T_n+T_n')}\\
 &\subseteq \overline{W(T_1,\ldots, T_n)+W(T_1',\ldots, T_n')}\\
 &\subseteq \overline{(\CC^n)}_{w_e(T_1,\ldots, T_n)+w_e(T_1',\ldots, T_n')}.
 \end{split}
 \end{equation*}

 Assume now that $T_1,\ldots, T_n\in B(\cH)$ are mutually commuting operators.
 Let $\cB$ be a closed subalgebra of $B(\cH)$ containing $T_1,\ldots, T_n$, and the identity.
 Denote by $\sigma(T_1,\ldots, T_n)$ the Harte spectrum, i.e., 
 $(\lambda_1,\ldots, \lambda_n)\in \sigma(T_1,\ldots, T_n)$ if and only if 
 $$(\lambda_1I-T_1)X_1+\cdots +(\lambda_n I-T_n)X_n\neq I
 $$
 for all $X_1,\ldots, X_n\in \cB$. According to \cite{M}, the joint numerical radius
 has the following property;
 $$
 r(T_1,\ldots, T_n)=\max\{\|(\lambda_1,\ldots, \lambda_n)\|_2:\ 
 (\lambda_1,\ldots, \lambda_n)\in \sigma(T_1,\ldots, T_n)\}.
 $$
 Hence, and using Theorem \ref{propri} part (vii), we deduce that
 \begin{equation}\label{harte}
 \sigma(T_1,\ldots, T_n)\subseteq 
 \overline{(\CC^n)}_{r(T_1,\ldots, T_n)} \subseteq 
 \overline{(\CC^n)}_{w(T_1,\ldots, T_n)}.
 \end{equation}

   In what follows,  we  calculate
  the joint right spectrum and  joint numerical range for certain classes of
  $n$-tuples of operators.

 \begin{theorem}\label{example}
 Let $(V_1,\ldots, V_n)$ be an $n$-tuple of  operators  
 on a Hilbert space $\cK$.  
 \begin{enumerate}
 \item[(i)] If $V_1,\ldots, V_n$ are isometries with orthogonal ranges, 
 then
 $V_1V_1^*+\cdots +V_nV_n^*\neq I$  if and only if
 \begin{equation*}
 \sigma_r(V_1,\ldots, V_n)=\overline{W(V_1,\ldots, V_n)}=\overline{\BB}_n.
 \end{equation*}

 \item[(ii)] If  $V_1V_1^*+\cdots +V_nV_n^*=I$,   then
 $\sigma_r(V_1,\ldots, V_n)\subseteq \partial \BB_n$    and   
  $w(V_1,\ldots, V_n)=1.$
  \end{enumerate}
 In particular, 
 if $S_1,\ldots, S_n$ are the left creation operators
  on the Fock space $F^2(H_n)$,
  then 
 $$
 \sigma_r(S_1,\ldots, S_n)=
 \overline{W(S_1,\ldots, S_n)}=\overline{(\CC^n)}_{w_e(S_1,\ldots, S_n)}=
 \overline{(\CC^n)}_{w(S_1,\ldots, S_n)}=
 \overline{\BB}_n.
 $$
 \end{theorem}
 
  \begin{proof}
  According to the  Wold decomposition for isometries with orthogonal ranges
   \cite{Po-isometric}, the Hilbert space $\cK$ admits an 
   orthogonal decomposition $\cK=\cK_s\oplus \cK_c$ such that $\cK_s, \cK_c$
   are  are reducing subspaces for each isometry
   $V_1,\ldots, V_n$, the row isometry $[V_1|\cK_s, \ldots, V_n|\cK_s]$ is unitarily 
   equivalent to $[S_1\otimes I_\cG,\ldots, S_n\otimes I_\cG]$ for some Hilbert space
   $\cG$, and $[V_1|\cK_c, \ldots, V_n|\cK_c]$ is a Cuntz row isometry, i.e.,
   $\sum_{i=1}^n (V_i|\cK_c)(V_i|\cK_c)^*=I_{\cK_c}$.
 Notice that 
   $V_1V_1^*+\cdots +V_nV_n^*\neq I$ if and only if $\cG\neq \{0\}$.
   
   According to the proof of Theorem \ref{conv}, given
    $\lambda:=(\lambda_1,\ldots, \lambda_n)\in \BB_n$, 
    we have $S_i^* z_\lambda=\overline{\lambda}_i z_\lambda$, \ 
    $i=1,\ldots, n$, where $z_\lambda$ is defined by 
    \eqref{zl}.   Therefore, if $h\in\cG$, $\|h\|=1$, then
    $$
    \left<(S_i\otimes I_\cG)(z_\lambda\otimes h), z_\lambda\otimes h\right>=
    \left<z_\lambda\otimes h, \overline{\lambda}_iz_\lambda\otimes h\right>=
    \lambda_i
    $$
    for any $i=1,\ldots, n$.
    This shows that
    $$
    \BB_n\subseteq W(S_1\otimes I_\cG,\ldots, S_n\otimes I_\cG)\subseteq
    W(V_1,\ldots, V_n).
    $$
    Since $\|[V_1,\ldots, V_n]\|=1$, Theorem \ref{conv} implies
    $\overline{W(V_1,\ldots, V_n)}= \overline\BB_n$.
    To complete the proof of part (i), it is enough to prove that 
    $\BB_n\subseteq \sigma_r(V_1,\ldots, V_n)$. To this end,
    let $(\lambda_1,\ldots, \lambda_n)\in \BB_n$ and $h\in \cG$ with $\|h\|=1$.
    Using the above-mentioned Wold type decomposition, we have
    $$
    \sum_{i=1}^n\|(\overline{\lambda}_i I_\cK-V_i^*)(z_\lambda\otimes h)\|^2=0.
    $$
    Due to Corollary \ref{approx}, we deduce that 
    $(\lambda_1,\ldots, \lambda_n)\in  \sigma_r(V_1,\ldots, V_n)$.
    
    To prove part (ii), fix  $\lambda:=(\lambda_1,\ldots, \lambda_n)\in \BB_n$.
    Since 
    $$\left\| \sum_{i=1}^n \lambda_i V_i^*\right\|\leq 
    \|\lambda\|_2\|[V_1,\ldots, V_n]\|=\|\lambda\|_2,
    $$
    we have
    $$
     \sum_{i=1}^n (\overline{\lambda}_i V_i^*+\lambda_i V_i)\leq 2\|\lambda\|_2.
     $$
     Using this inequality, we obtain
     \begin{equation*}\begin{split}
     \sum_{i=1}^n (\lambda_iI_\cK-V_i)(\overline{\lambda}_i I_\cK-V_i^*)
     &= \sum_{i=1}^n|\lambda_i|^2 I_\cK-\sum_{i=1}^n (\overline{\lambda}_i
      V_i^*+\lambda_i V_i)+\sum_{i=1}^n V_iV_i^*\\
      &\geq \left(\|\lambda\|_2^2-2\|\lambda\|_2+1\right)
       I_\cK=(1-\|\lambda\|_2)^2 I_\cK.
     \end{split}
     \end{equation*}
     Since $\|\lambda\|_2\neq 1$, Corollary \ref{approx} implies
      $(\lambda_1,\ldots, \lambda_n)\notin  \sigma_r(V_1,\ldots, V_n)$.
      Using part (i) of Theorem \ref{conv}, we deduce that 
      $\sigma_r(V_1,\ldots, V_n)\subseteq \partial \BB_n$. 
      
       To prove that
      $w(V_1,\ldots, V_n)=1$, notice that the operator
      $Z:=\sum_{i=1}^n S_i\otimes V_i^*$ is a nonunitary isometry.
      Applying part (ii) of this theorem (in the particular case of a single
      isometry) to $Z$, we deduce  that
      $w(V_1,\ldots, V_n)=w(Z)=1$.
       
      The particular case when $S_1,\ldots, S_n$ are the left creation operators
      follows from part (i) of this theorem and  Theorem
      \ref{conv} (part (i) and (ii)).
       The proof is complete.
  \end{proof}

We established a strong connection between the algebra $F_n^\infty$
  and the function theory on the open unit ball  
  $$\BB_n:=\{(\lambda_1,\ldots, \lambda_n)\in \CC^n:\  
  |\lambda_1|^2+\cdots +|\lambda_n|^2<1\},
  $$
   through the noncommutative von Neumann inequality \cite{Po-von} 
   (see also  \cite{Po-funct}, \cite{Po-disc}, \cite{Po-poisson}, and \cite{Po-tensor}).
    In particular, we proved that there is a completely contractive
     homomorphism $\Phi:F_n^\infty\to H^\infty(\BB_n)$
     defined by
 $$
       [\Phi(f(S_1,\dots, S_n))](\lambda_1,\dots, \lambda_n)=f(\lambda_1,\dots,
       \lambda_n)$$ for any
 $ f(S_1,\dots, S_n)\in F_n^\infty$ and
  $ (\lambda_1,\dots, \lambda_n)\in \BB_n.
$
 A characterization of the analytic functions in the range of the map
   $\Phi$ was obtained  in \cite{ArPo2}, and independently in \cite{DP}.   
   Moreover, it was proved that the quotient $F_n^\infty/{\ker\Phi}$ is an
  operator algebra which can be identified with
   $W_n^\infty:=P_{F_s^2(H_n)} F_n^\infty|_{F_s^2(H_n)}$,
  the compression  of $ F_n^\infty$ to the symmetric Fock space
   ${F_s^2(H_n)} \subset F^2(H_n)$.

 It is well-known that the joint numerical range  of an $n$-tuple of operators 
 in not convex in 
 general if $n\geq 2$. The only known exceptions are in
  the commutative case, when  
 $T_1,\ldots, T_n$ are either bounded analytic operators  or double commuting
 operators
 (see \cite{De}).
 In what follows, we show that there are important
  classes of noncommuting operators for which the joint numerical range 
  or its closure are convex.

If  $f_i  \in F_n^\infty$, \ $i=1,\ldots, k$, we denote by 
$\sigma_r(f_1,\ldots, f_k)$ the right joint spectrum  
  with respect to the
   noncommutative analytic Toeplitz algebra $F_n^\infty$.

  \begin{theorem}\label{spectru}
Let $(f_1,\ldots, f_k)$ be a $k$-tuple of operators in 
  the noncommutative analytic Toeplitz algebra $F_n^\infty$. 
     Then
 the following properties hold:
  \begin{enumerate}
  \item[(i)]
  $ W(P_{\cP_{m}}f_1|\cP_{ m},\ldots,
   P_{\cP_{m}}f_k|\cP_{m})$ is a convex compact subset of $\CC^k$,
   where $\cP_{m}$ is the set of all polynomials in $F^2(H_n)$ of degree
   $\leq m$;
  \item[(ii)]
  $\overline{W(f_1,\ldots, f_k)}$ is a convex compact subset of $\CC^k$;
  \item[(iii)] $(\lambda_1,\ldots, \lambda_n)\notin \sigma_r(f_1,\ldots, f_k)$ if and only
  if there is $\delta>0$ such that
  $$(\lambda_1 I-f_1)(\overline{\lambda}_1 I-f_1^*)+\cdots +
  (\lambda_k I-f_k)(\overline{\lambda}_k I-f_k^*)\geq \delta^2 I;$$
  \item[(iv)]
  $\sigma_r(f_1,\ldots, f_k)$ is a  compact subset of $\CC^k$ and
  $$
  \{(f_1(\lambda ), \ldots, 
  f_k(\lambda )):\  \lambda \in\BB_n\}^{-}\subseteq \sigma_r(f_1,\ldots, f_k)\subseteq
   \overline{W(f_1,\ldots, f_k)}.
   $$
   \end{enumerate}
  \end{theorem}
  \begin{proof}
  We recall that the flipping operator $U\in B(F^2(H_n))$
  is defined by setting $U(1)=1$ and $U(e_\alpha)=e_{\tilde \alpha}$, where
  $\tilde\alpha:= g_{i_m}\cdots g_{i_2}g_{i_1}$ is the reverse
  of 
  $\alpha:=g_{i_1}g_{i_2}\cdots g_{i_m}$.  If 
  $f=\sum_{\alpha\in \FF_n^+} a_\alpha e_\alpha$ is a vector in $F^2(H_n)$, 
   we denote
   $\tilde f:=\sum_{\alpha\in \FF_n^+} a_\alpha e_{\tilde \alpha}.
   $ 
  Now, let $p, q\in \cP_m $ such that $\|p\|_2=\|q\|_2=1$, and let $t\in (0,1)$.
  For each $j=1,\ldots, k$, we have
  \begin{equation*}
  \begin{split}
  t&\left<f_j(S_1,\ldots, S_n)p, p\right>+
  (1-t) \left<f_j(S_1,\ldots, S_n)q, q\right>=
  t\left<\tilde{p}\otimes \tilde{f_j}, \tilde{p}\right>+(1-t)
  \left<\tilde{q}\otimes \tilde{f_j}, \tilde{q}\right>\\
  &=t\left<\tilde{f_j}, \tilde{p}(S_1,\ldots, S_n)^* 
  \tilde{p}(S_1,\ldots, S_n)(1)\right>+(1-t)
  \left<\tilde{f_j}, \tilde{q}(S_1,\ldots, S_n)^* 
  \tilde{q}(S_1,\ldots, S_n)(1)\right>\\
  &=
  \left<\tilde{f_j}, [t\tilde{p}(S_1,\ldots, S_n)^* 
  \tilde{p}(S_1,\ldots, S_n)+(1-t)
  \tilde{q}(S_1,\ldots, S_n)^* 
  \tilde{q}(S_1,\ldots, S_n)](1)\right>.
  \end{split}
  \end{equation*}
  Define the operator $X\in B(F^2(H_n)$ by setting 
  $$
  X:=t\tilde{p}(S_1,\ldots, S_n)^* 
  \tilde{p}(S_1,\ldots, S_n)+(1-t)
  \tilde{q}(S_1,\ldots, S_n)^* 
  \tilde{q}(S_1,\ldots, S_n)
  $$
  and notice that $U^* XU$ is a positive  multi-Toeplitz operator.
    Applying Theorem 1.6 from \cite{Po-analytic} to $U^* XU$,
   we find a polynomial $s\in \cP_m$ such that
   $$
   U^* XU=\tilde{s}(R_1,\ldots, R_n)^* 
  \tilde{s}(R_1,\ldots, R_n),
  $$
  where $R_1,\ldots, R_n$ are the right creation operators
   on the full Fock space
  $F^2(H_n)$.
  Notice also that
  \begin{equation*}
  \begin{split}
  \|s\|_2^2&=
  \left< \tilde{s}(S_1,\ldots, S_n)^* 
  \tilde{s}(S_1,\ldots, S_n)(1),1\right>\\
  &=\left<
  [t\tilde{p}(S_1,\ldots, S_n)^* 
  \tilde{p}(S_1,\ldots, S_n)+(1-t)
  \tilde{q}(S_1,\ldots, S_n)^* 
  \tilde{q}(S_1,\ldots, S_n)](1),1 \right>\\
  &=t\|p\|^2_2+  (1-t) \|q\|^2_2=1.
  \end{split}
  \end{equation*}
  The above calculations reveal that
  \begin{equation*}\begin{split}
  t\left<f_j(S_1,\ldots, S_n)p, p\right>&+
  (1-t) \left<f_j(S_1,\ldots, S_n)q, q\right>\\
  &=\left<\tilde{f_j},
  U^* \tilde{s}(R_1,\ldots, R_n)^* 
  \tilde{s}(R_1,\ldots, R_n)U (1)\right>\\
  &=\left< \tilde{f_j},\tilde{s}(S_1,\ldots, S_n)^* 
  \tilde{s}(S_1,\ldots, S_n))(1)\right>\\
  &=\left<\tilde{s}\otimes\tilde{f_j}, \tilde{s}\right>=
  \left< f_j\otimes s, s\right>\\
 &= \left<f_j(S_1,\ldots, S_n)s,s\right>,
  \end{split}
  \end{equation*}
  for any $j=1,\ldots, k$.
 This shows that  the joint numerical range of the $n$-tuple  \linebreak 
 $(P_{\cP_{m}}f_1|\cP_{ m},\ldots,
   P_{\cP_{m}}f_k|\cP_{m})$ is convex.
   
  Now, we prove (ii). Let $\phi, \psi\in F^2(H_n)$ be such that 
  $\|\psi\|_2=\|\psi\|_2=1$, and define the vectors
  $$
  \phi_m:=\frac{1} {\|P_{\cP_m}\phi\|} P_{\cP_m}\phi \quad \text{ and } \quad 
  \psi_m:=\frac{1} {\|P_{\cP_m}\psi\|} P_{\cP_m}\psi.
  $$
  Notice that $\phi_m,\psi_m\in \cP_m$ and $\|\phi_m\|_2=\|\psi_m\|_2=1$ for any
  $m=1,2,\ldots $. Moreover, it is clear that
  $\|\phi_m-\phi\|_2\to 0$ and $\|\psi_m-\psi\|_2\to 0$ as $m\to\infty$.
  According to (i), there exist polynomials
   $s_m\in \cP_m$ with $\|s_m\|_2=1$ such that
 $$
 t\left<f_j(S_1,\ldots, S_n)\phi_m, \phi_m\right>+
  (1-t) \left<f_j(S_1,\ldots, S_n)\psi_m, \psi_m\right> 
  =\left< f_j(S_1,\ldots, S_n) s_m,s_m\right>
  $$
  for any $j=1,\ldots, k$ and $m=1,2\ldots$.
  Taking $m\to\infty$, we deduce that
  $$
 t\left<f_j(S_1,\ldots, S_n)\phi, \phi\right>+
  (1-t) \left<f_j(S_1,\ldots, S_n)\psi, \psi\right> 
  =\lim_{m\to\infty}\left< f_j(S_1,\ldots, S_n) s_m,s_m\right>
  $$
  for any $j=1,\ldots, k$.
  Therefore 
  \begin{equation}\label{conve}
  \text{\rm conv}\,W(f_1,\ldots, f_k)\subseteq \overline{W(f_1,\ldots, f_k)}.
   \end{equation}
   Now, let 
 $ {\lambda},  {\mu}\in \overline{W(f_1,\ldots, f_k)}$ 
 and let
 $\left\{ {\lambda}^{(m)}\right\}_{m=1}^\infty,
  \left\{ {\mu}^{(m)}\right\}_{m=1}^\infty$ be sequences of vectors
  in $W(f_1,\ldots, f_k)$ such that 
  $ {\lambda}^{(m)}\to {\lambda}$ and 
  $ {\mu}^{(m)}\to {\mu}$ as $m\to\infty$.
  Using relation \eqref{conve}, we deduce that
  $$
  t {\lambda}^{(m)}+(1-t) {\mu}^{(m)}\in 
  \overline{W(f_1,\ldots, f_k)}.
  $$
  Taking $m\to\infty$,  we complete the proof of part   (ii). The property (iii)
  follows from Theorem 3.3 of \cite{Po-analytic}.
  
  We prove now part (iv). According to Theorem \ref{conv}, 
  $\sigma_r(f_1,\ldots, f_k)$ is a compact subset of $\CC^k$.
  Let $(\lambda_1,\ldots, \lambda_n)\in \BB_n$, 
  $f(S_1,\ldots, S_n)\in F_n^\infty$,
  and $p\in \cP$. Using the properties of $z_\lambda$  (see 
  the proof of Theorem \ref{conv}), we have
  \begin{equation*}\begin{split}
  \left<f(S_1,\ldots, S_n)^* z_\lambda,p\right>&=
  \left< z_\lambda, f\otimes p\right>=\overline{f(\lambda_1,\ldots,\lambda_n)}
  \overline{p(\lambda_1,\ldots,\lambda_n)}\\
  &=\left<\overline{f(\lambda_1,\ldots,\lambda_n)} z_\lambda, p\right>
  \end{split}
  \end{equation*}
  for any polynomial $p \in \cP$ . Since $\cP$ is dense 
  in $F^2(H_n)$, we deduce that
  $$
  f(S_1,\ldots, S_n)^* z_\lambda=\overline{f(\lambda_1,\ldots,\lambda_n)}
   z_\lambda.
  $$
 Denote $\mu_j:= f_j(\lambda_1,\ldots,\lambda_n)$, \ $j=1,\ldots, k$, and notice that
 $$
 \sum_{j=1}^k \|(\overline{\mu}_jI-f(S_1,\ldots, S_n)^*) z_\lambda\|^2=0.
 $$
 Since $\|z_\lambda\|=1$, we can use again Corollary
 \ref{approx} and deduce that
  $(\mu_1,\ldots, \mu_k)\in \sigma_r(f_1,\ldots, f_k)$.
  This completes the proof of part (iv).
 \end{proof}

  The next result is a commutative version of Theorem 3.1 from \cite{Po-analytic}.
   The proof follows the same lines. We include it for completeness.

  \begin{theorem}\label{facto}
  If $A\in W_n^\infty\bar\otimes B(\cH, \cH')$ 
  and $B\in W_n^\infty\bar\otimes B(\cH'', \cH')$, then 
  there exists a contraction  $C\in W_n^\infty\bar\otimes B(\cH, \cH'')$ 
  such that $A=BC$  if and only if 
  $AA^*\leq BB^*$.
  \end{theorem}
  \begin{proof} 
  One implication is clear, so assume that $AA^*\leq BB^*$.
  Therefore, there is a contraction
  $X:\cM:=\overline{B^*(F_s^2(H_n)\otimes \cH')}\to F_s^2(H_n)\otimes\cH$
   satisfying $XB^*=A^*$.
  Since $B\in W_n^\infty\bar\otimes B(\cH'', \cH')$, we have 
  $$
  B (B_i\otimes I_{\cH''})=(B_i\otimes I_{\cH'})B, \quad i=1,\ldots, n,
  $$ 
  where $B_i:=P_{F_s^2(H_n)} S_i | F_s^2(H_n)$, $i=1,\ldots, n$.
    Hence, the subspace $\cM\subseteq F_s^2(H_n)$ is invariant  under each
  operator
  $B_i^*\otimes I_{{\cH''}}$, $i=1,\ldots, n$.
  Now, we define the operators
   $T_i:= (B_i^*\otimes I_{\cH''})|\cM$, $i=1,\ldots, n$, 
   acting from $\cM$ to $\cM$.
  Since 
  $A\in W_n^\infty\bar\otimes B(\cH, \cH')$,
   $B\in W_n^\infty\bar\otimes B(\cH'', \cH')$, and $XB^*=A^*$, we have 
  \begin{equation*}
  \begin{split}
  X(B_i^*\otimes I_{\cH''}) B^*k&=XB^* (B_i^*\otimes I_{\cH'})k=
  A^*(B_i^*\otimes I_{\cH'})k\\
  &=
  (B_i^*\otimes I_{\cH})A^*k=
  (B_i^*\otimes I_{\cH}) X B^*k 
  \end{split}
  \end{equation*}
  for any $k\in F_s^2(H_n)\otimes \cH'$.
  Hence,
  \begin{equation}
  \label{intert}
  X(B_i^*\otimes I_{\cH''})=(B_i^*\otimes I_{\cH}) X,\quad i=1,\ldots, n.
  \end{equation}
  Since $F_s^2(H_n)$ is an invariant subspace under each operator
   $S_1^*,\ldots, S_n^*$,  and 
  $B_i^*=S_i^*|F_s^2(H_n)$ for $i=1,\ldots, n$, relation \eqref{intert} implies
  $$
  [P_\cM(S_i\otimes I_{\cH''})|\cM] Y^*=Y^*(S_i\otimes I_\cH), \quad i=1,\ldots, n,
  $$
  where $Y:\cM\to F^2(H_n)\otimes \cH$ and $Yh=Xh$ for any  $h\in \cM$.
 Since $[ S_1\otimes I_{\cH''},\ldots, S_n\otimes I_{\cH''}]$
 is an isometric dilation of the row contraction
 $[P_\cM(S_1\otimes I_{\cH''})|\cM,\ldots, P_\cM(S_n\otimes I_{\cH''})|\cM]$,
  we use the noncommutative commutant lifting theorem \cite{Po-isometric} 
  (see \cite{SzF-book} for the classical case)
  to find
  a contraction $\tilde{C}\in B(F^2(H_n)\otimes 
  \cH, F^2(H_n)\otimes \cH'')$
   such that
  $$
  (S_i\otimes I_{\cH''})\tilde{C}=\tilde{C}(S_i\otimes I_{\cH}), 
  \quad i=1,\ldots, n,
  $$
  and $\tilde{C}^*k=Yk=Xk$ for any $k\in \cM$.
 Hence, we have $X=P_{F_s^2(H_n)\otimes \cH} \tilde{C}^*|\cM $ and, 
 taking into account that $B\cM^\perp =0$, we obtain
 \begin{equation*}
 \begin{split}
 A&=BX^*=BP_\cM \tilde{C}|F_s^2(H_n)\otimes \cH\\
 &=BP_\cM
 P_{F_s^2(H_n)\otimes \cH} \tilde{C}^*|F_s^2(H_n)\otimes \cH=
 BC,
 \end{split}
 \end{equation*}
 where $C:=P_{F_s^2(H_n)\otimes \cH} \tilde{C}^*|F_s^2(H_n)\otimes \cH$. 
  According to \cite{Po-analytic},
  since $\tilde{C}$ is a multi-analytic operator, we have
   $\tilde{C}\in F_n^\infty\overline{\otimes} B(\cH,\cH'')$ and therefore
   $C\in W_n^\infty\overline{\otimes} B(\cH,\cH'')$.
   This completes the proof of the theorem.
 \end{proof}

  Applying Theorem \ref{facto} to  the particular case
   when $\cH=\cH'$ and  $A=\delta I_{F_s^2(H_n)\otimes \cH'}$,
   $\delta>0$,  one can
   easily obtain the following consequence.

  \begin{corollary}\label{coron}
  Let  $B\in W_n^\infty\bar\otimes B(\cH'', \cH')$. The 
  following statements are equivalent:
  \begin{enumerate}
  \item[(i)] There is $D\in B(F_s^2(H_n)\otimes \cH', F_s^2(H_n)\otimes \cH'')$ such that 
  $BD=I;$
  \item[(ii)] There is $\delta>0$ such that $\|B^*k\|\geq \delta\|k\|$, for any
  $k\in F_s^2(H_n)\otimes \cH';$
  \item[(iii)] There is $C\in W_n^\infty\bar\otimes B(\cH', \cH'')$
  such that
  $BC=I$.
  \end{enumerate}
  \end{corollary}
  
  Now, we can obtain as a consequence the following corona type result
   for the algebra $W_n^\infty\bar\otimes B(\cH)$.
  
  \begin{corollary}\label{corona2}
  If  $f_1,\ldots, f_k\in W_n^\infty\bar\otimes B(\cH)$, then there exist 
  $g_1,\ldots, g_k\in W_n^\infty\bar\otimes B(\cH)$ such that
  $$
  f_1g_1+\cdots +f_kg_k=I
  $$
  if and only if there exists $\delta>0$ such that
  $$
  f_1f_1^*+\cdots +f_kf_k^*\geq \delta^2 I.
  $$
   \end{corollary}
  \begin{proof}
  Take $B:=[f_1,\ldots, f_k]\in W_n^\infty\bar\otimes B(\cH^{(k)}, \cH)$,  
  $C:=\left[\begin{matrix}g_1\\\vdots\\
  g_k\end{matrix}\right]\in W_n^\infty\bar\otimes B(\cH,\cH^{(k)})$, and apply 
  Corollary \ref{coron}.
  \end{proof}

The proof of the following result is straightforward, so we omit it.
  
  \begin{lemma}\label{prop2}  
  Let $(T_1,\ldots, T_n)$, $T_i\in B(\cH)$,   and  
  $(V_1,\ldots, V_n)$,
  $V_i\in B(\cK)$, be   $n$-tuples of 
   operators  such that
   $\cH \subseteq \cK$.
   \begin{enumerate}
   \item[(i)] If ~$T_i=P_\cH V_i|\cH$, \ $i=1,\ldots, n$, then
   $W(T_1,\ldots, T_n)\subseteq W(V_1,\ldots, V_n)$ 
  
$$w_e(T_1,\ldots, T_n)\leq
   w_e(V_1,\ldots, V_n), \  \text{ and } \ w(T_1,\ldots, T_n)\leq
   w(V_1,\ldots, V_n).
   $$
   \item[(ii)] If $T_i^*=V_i^*|\cH$, \ $i=1,\ldots, n$,
   then 
   $$\sigma_r(T_1,\ldots, T_n)\subseteq \sigma_r(V_1,\ldots, V_n)\ \text{  and
}\ 
   r(T_1,\ldots, T_n)\leq r(V_1,\ldots, V_n).
   $$
   \end{enumerate}
  \end{lemma}

  We recall that the   Harte spectrum of  a $k$-tuple $(g_1,\ldots, g_k)$, 
  $g_i\in W_n^\infty$,  relative to the commutative algebra algebra
  $W_n^\infty$, is defined by
  $$
  \sigma(g_1,\ldots, g_k):=\left\{(\lambda_1,\ldots, \lambda_k)\in \CC^n:\ 
  (\lambda_1 I-g_1)W_n^\infty+\cdots +
  (\lambda_k I-g_k)W_n^\infty\neq W_n^\infty\right\}.
  $$
  
  Here are some of the properties of the Harte spectrum of $(g_1,\ldots, g_k)$.
    
 \begin{theorem}\label{corona}
 If  $g_1,\ldots, g_k\in W_n^\infty$, then
   \begin{enumerate}
   \item[(i)]
  $(\lambda_1,\ldots, \lambda_k)\notin\sigma(g_1,\ldots, g_k)$ if and
  only if there is $\delta>0$ such
  that
  $$
  (\lambda_1 I-g_1)(\overline{\lambda}_1I-g_1^*)+\cdots + (\lambda_k I-g_k) 
  (\overline{\lambda}_k I-g_k^*)\geq \delta^2 I.
  $$
 \item[(ii)]
   The  Harte spectrum
  $\sigma(g_1,\ldots, g_k)$ is a  compact subset of $\CC^k$ and 
  $$\{(g_1(\lambda), \ldots, 
  g_k(\lambda)):\  \lambda\in\BB_n\}^{-}\subseteq \sigma(g_1,\ldots, g_k).
   $$
   \item[(iii)] $\sigma(g_1,\ldots, g_k)\subseteq 
   \left(\CC^n\right)_{r(g_1,\ldots, g_k)}
   \subseteq 
   \left(\CC^n \right)_{w_e(g_1,\ldots, g_k)}.$
  \item[(iv)] The following inclusions hold:
  \begin{equation} \label{sig-s1}
  \sigma(g_1,\ldots, g_k)
  \subseteq\bigcap \sigma_r(f_1,\ldots, f_k)
  \subseteq\bigcap\overline{[W(f_1,\ldots, f_k)]}
  \end{equation}
  and 
  \begin{equation}\label{sig-s2}
  \sigma(g_1,\ldots, g_k)\subseteq \overline{W(g_1,\ldots, g_k)}
  \subseteq \bigcap
  \overline{W(f_1,\ldots, f_k)} \subseteq \overline{(\CC^k)}_q,
  \end{equation}
  where the intersections are taken over all $f_i\in F_n^\infty$ with 
  $P_{F_s^2(H_n)} f_i|F^2_s(H_n)=g_i$, \ $i=1,\ldots, k$,
    and 
  $$
  q:=\inf \{w_e(f_1,\ldots, f_k):\  f_i\in F_n^\infty, \ P_{F_s^2(H_n)}
   f_i|F^2_s(H_n)=g_i,\  i=1,\ldots, k\}.
  $$
   \end{enumerate}
   \end{theorem}
  \begin{proof}
  Part (i) follows from Corollary \ref{corona2} if we take
   $f_i=\lambda_iI-g_i$, $i=1,\ldots, n$.
   To prove  part (ii),  we recall that any element of $W_n^\infty$ has the 
   form $g(B_1,\ldots, B_n)$, where $g(S_1,\ldots, S_n)\in F_n^\infty$, 
   $B_i=P_{F_s^2(H_n)} S_i|F_s^2(H_n)$, $i=1,\ldots, n$, and $g(B_1,\ldots, B_n)$ is defined by 
   the $F_n^\infty$-functional calculus for row contractions (see \cite{Po-funct}).
   Using the proof of part (iv) of Theorem \ref{spectru}, we have
   $$
   g(S_1,\ldots, S_n)^* z_\lambda =\overline{g(\lambda_1,
    \ldots, \lambda_n)} z_\lambda,
   $$
   where  the vector $z_\lambda$ is defined by \eqref{zl}.
   Since $z_\lambda\in F_s^2(H_n)$ and the symmetric Fock space is  an invariant 
   subspace under each operator $S_1^*, \ldots, S_n^*$, and  $F_n^\infty$ is the 
    WOT-closure of the polynomials in $S_1,\ldots, S_n$ and the identity,
    we deduce that
    \begin{equation}\label{*zl}
    g(B_1,\ldots, B_n)^* z_\lambda =\overline{g(\lambda_1, 
    \ldots, \lambda_n)} z_\lambda.
    \end{equation}
   Let $g_1,\ldots, g_n\in W_n^\infty$ and denote $\mu_j:=
   g_j(\lambda_1, \ldots, \lambda_n)$.
   Using relation \eqref{*zl}, we get
   $$
   \sum_{j=1}^k \|(\overline{\mu}_j I-g_j(B_1,\ldots, B_n)^*) z_\lambda\|^2=0.
   $$
   Since $\|z_\lambda\|=1$, part (i) of this theorem shows that 
   $(\mu_1,\ldots, \mu_k)\in \sigma(g_1,\ldots, g_k)$,
    which completes the proof of (ii).
   
    Part  (iii) follows from  Theorem \ref{conv} and the fact that, in
    the commutative case (see \cite{M}),
    $$
    r(g_1,\ldots, g_k)=\sup \{\|\lambda\|_2:\ \lambda\in \sigma(g_1,\ldots, g_k)\}.
    $$
To prove (iv), note that if $\lambda\notin \bigcap \sigma_r(f_1,\ldots, f_k)$,
   then there exist  operators $f_1,\ldots, f_k\in F_n^\infty$ such that
   $\lambda\notin \sigma_r(f_1,\ldots, f_k)$. According
   to Theorem \ref{spectru} part (iii),  there is $\delta>0$ such that 
   $$
  (\lambda_1 I-f_1)(\overline{\lambda}_1I-f_1^*)+\cdots + (\lambda_k I-f_k) 
  (\overline{\lambda}_k I-f_k^*)\geq \delta^2 I.
  $$
   Since $F_s^2(H_n)$ is invariant under $S_1^*, \ldots, S_n^*$ and taking the
    compression  to the symmetric Fock space, the latter inequality implies
  \begin{equation}
\label{de2}
  (\lambda_1 I-g_1)(\overline{\lambda}_1I-g_1^*)+\cdots + (\lambda_k I-g_k) 
  (\overline{\lambda}_k I-g_k^*)\geq \delta^2 I.
  \end{equation}  
   Using now part (i) of the theorem we deduce that $\lambda\notin 
   \sigma(g_1,\ldots, g_k)$, which proves the first inclusion
   in \eqref{sig-s1}. The second inclusion 
   is implied by Theorem \ref{conv} part (i).
   Now, notice that 
   \begin{equation}
   \label{sisi}
   \sigma(g_1,\ldots, g_k)\subseteq \sigma_r(g_1,\ldots, g_k).
   \end{equation}
   Indeed, if $(\lambda_1,\ldots, \lambda_n)\notin \sigma_r(g_1,\ldots, g_k)$, then,
   according to Corollary \ref{approx}, there is $\delta>0$ such that
    \eqref{de2} holds.
   Using part (i) of the theorem, we obtain \eqref{sisi}. According to Theorem
   \ref{conv} part (i), we get the first inclusion in \eqref{sig-s2}.
     The second inclusion
   follows 
   from
    Lemma \ref{prop2}.
    Once again, if we apply Theorem \ref{conv} part (i) to our setting,
     we deduce the
    last inclusion in \eqref{sig-s2}.
The proof is complete.
 \end{proof}

 Notice that if  
 $B_1,\ldots, B_n$  are the creation operators  acting on the symmetric Fock space, 
 then Theorem \ref{corona} implies
$$
 \sigma(B_1,\ldots, B_n)=
 \overline{W(B_1,\ldots, B_n)}=\overline{(\CC^n)}_{w_e(B_1,\ldots, B_n)}=
 \overline{(\CC^n)}_{w(B_1,\ldots, B_n)}=
 \overline{\BB}_n.
 $$

  \smallskip
    
 \section{$\rho$-operator radius}
 \label{norms}

 The results of this section can be seen as the unification of
     the theory of isometric dilations for row contractions \cite{Sz1},
      \cite{SzF-book},
      \cite{Fr}, \cite{Bu}, \cite{Po-models}, \cite{Po-isometric}, 
      \cite{Po-charact}
      (which corresponds to the case $\rho=1$) 
     and Berger type dilations  of Section 1 for $n$-tuples $(T_1,\ldots, T_n)$ with the joint numerical
     radius $w(T_1,\ldots, T_n)\leq 1$ (which corresponds to the case $\rho=2$). 
 We obtain
 several intrinsic characterizations for $n$-tuples of operators of class
     $\cC_\rho$,
  present   basic properties of the joint
     $\rho$-operator radius  $\omega_\rho$, and extend
       to our multivariable setting (noncommutative and commutative)
     several classical results obtained  by Sz.-Nagy and Foia\c s, Halmos,
       Berger and Stampfli, Holbrook, Paulsen, and  others   
     (\cite{B}, \cite{BS}, \cite{BS1},  \cite{Ha}, \cite{Ha1}, 
     \cite{Ho1}, \cite{Ho2}, \cite{K}, \cite{Pa-book}, \cite{Pe}, \cite{SzF},
       and \cite{W}).
       
  An  $n$-tuple of operators $(T_1,\ldots, T_n)$, \ $T_i\in B(\cH)$,
   belongs to the 
   class $\cC_\rho$, $\rho>0$, if there exist a Hilbert 
  space $\cK\supseteq \cH$
 and isometries  $V_i\in B(\cK)$, $i=1,\ldots, n$, with orthogonal
  ranges,  such that
  \begin{equation}\label{ro}
  T_\alpha =\rho P_\cH V_\alpha |\cH\quad \text{ for any } 
\alpha\in \FF_n^+\backslash\{g_0\},\end{equation}where
  $P_\cH$ is the orthogonal projection of $\cK$ onto $\cH$.
 We call $(V_1,\ldots, V_n)$ a $\rho$-isometric dilation of the   
 $n$-tuple  $(T_1,\ldots, T_n)$.
 Notice that relation \eqref{ro} implies that  any $n$-tuple $(T_1,\ldots, T_n)$ 
 of class $C_\rho$ is power bounded   and $r(T_1,\ldots, T_n)\leq 1$.
 We recall  that an $n$-tuple  of operators  $(T_1,\ldots, T_n)$ 
 is called  power bounded if there exists $M>0$ such that
 $$\left\|\sum_{|\alpha|=k} T_\alpha T_\alpha^*\right\|^{1/2}\leq M
 $$
 for any $k=1,2,\ldots$.

The first result of this section  provides characterizations for 
$n$-tuples of operators of class $\cC_\rho$. We recall that $R_1,\ldots, R_n$ denote
the right creation operators acting on the full Fock space.

  \begin{theorem}\label{rho}
  If  $T_1,\ldots, T_n\in B(\cH)$, then
    the following statements are equivalent:
  \begin{enumerate}
 \item[(i)] 
  $(T_1\ldots, T_n)\in \cC_\rho$;
  \item[(ii)] The inequality
  \begin{equation}\label{ro-2}
  (\rho-2)\left\|\left(I-\sum_{i=1}^n z_i R_i\otimes T_i^*\right)k\right\|^2 +2~
  \text{\rm Re}\left< \left(I-\sum_{i=1}^n z_iR_i\otimes T_i^*\right)k,
  k\right>\geq 0
  \end{equation}
  holds for any $k\in F^2(H_n)\otimes \cH$ and $z_i\in {\DD}$, 
  $i=1,\ldots, n$.
   \item[(iii)]The inequality
   \begin{equation}\label{1-ro}
   \left(1-\frac{1}{\rho}\right)\sum_{i=1}^n(\overline{z}_i R_i^*\otimes T_i+
   z_i R_i\otimes T_i^*)+
   \left(\frac{2}{\rho} -1\right) \sum_{i=1}^n |z_i|^2(I\otimes T_iT_i^*)
   \leq I
   \end{equation}
   holds for any $z_i\in {\DD}$, 
  $i=1,\ldots, n$.
  \item[(iv)]  
   The spectral radius $r(T_1\ldots, T_n)\leq 1$ and 
   the multi-Toeplitz operator
  $$
   \sum_{k=1}^\infty \sum_{|\alpha|=k}\overline{z}_\alpha 
  R_\alpha^*\otimes T_{\tilde \alpha} +\rho I+
  \sum_{k=1}^\infty \sum_{|\alpha|=k} z_\alpha 
  R_\alpha\otimes T_{\tilde \alpha}^*
  $$
  is positive for any ~$z_i\in {\DD}$, 
  $i=1,\ldots, n$, where the convergence is in the operator norm.
   \item[(v)] The spectral radius $r(T_1\ldots, T_n)\leq 1$ and 
  the multi-Toeplitz operator
  $M_r(T_1,\ldots, T_n)$ defined by
  \begin{equation}\label{mr}
  M_r(T_1,\ldots, T_n):= \sum_{k=1}^\infty \sum_{|\alpha|=k}r^k 
  R_\alpha^*\otimes T_{\tilde \alpha} +\rho I+
  \sum_{k=1}^\infty \sum_{|\alpha|=k}r^k 
  R_\alpha\otimes T_{\tilde \alpha}^*
  \end{equation}
  is positive for any $0<r<1$, where the convergence is in the operator norm.
  \end{enumerate}
  \end{theorem}
  \begin{proof}
 Assume that  $(T_1,\ldots, T_n)\in \cC_\rho$.  Then there
 exist a Hilbert space $\cK\supseteq \cH$
 and isometries  $V_i\in B(\cK)$, $i=1,\ldots, n$, with orthogonal
  ranges,  such that relation \eqref{ro} holds. According to Lemma \ref{spec-ra}, if 
  $z_i\in \DD$, $i=1,\ldots, n$, then the series 
 $\sum\limits_{k=1}\sum\limits_{|\alpha|=k} 
   z_\alpha R_\alpha\otimes V_{\tilde \alpha}^* $ is convergent 
   in the operator norm and 
   $$
   I+2\sum_{k=1}\sum_{|\alpha|=k} 
   z_\alpha R_\alpha\otimes V_{\tilde \alpha}^*=
  \left(I+\sum_{i=1}^n z_iR_i\otimes 
   V_i^*\right) \left(I-\sum_{i=1}^n z_iR_i\otimes 
   V_i^*\right)^{-1}.
   $$
   Hence, using  relation \eqref{ro} and Lemma \ref{spec-ra}, we deduce that
   \begin{equation}\label{roo}
   \begin{split}
   \left( 1-\frac{2}{\rho}\right) I&+ \frac{2}{\rho}
   \left(I-\sum_{i=1}^n z_iR_i\otimes 
   T_i^*\right)^{-1}\\
   &=P_{F^2(H_n)\otimes \cH}
  \left(I+\sum_{i=1}^n z_iR_i\otimes 
   V_i^*\right) \left(I-\sum_{i=1}^n z_iR_i\otimes 
   V_i^*\right)^{-1}|F^2(H_n)\otimes \cH,
   \end{split}
   \end{equation}
   where $P_{F^2(H_n)\otimes \cH}$ is the orthogonal projection
   of $F^2(H_n)\otimes \cK$ onto $F^2(H_n)\otimes \cH$.
  Setting $y:=\left(I-\sum_{i=1}^n z_iR_i\otimes 
   V_i^*\right)x$, where $x\in F^2(H_n)\otimes \cK$,  straightforward calculations
    show that
   \begin{equation*}
   \begin{split}
   \text{\rm Re} \left<\left(I+\sum_{i=1}^n z_iR_i\otimes 
   V_i^*\right)\right. &\left.\left(I-\sum_{i=1}^n z_iR_i\otimes 
   V_i^*\right)^{-1} y, y \right> \\
  & =\text{\rm Re}\left<
   \left(I+\sum_{i=1}^n z_iR_i\otimes 
   V_i^*\right)x, \left(I-\sum_{i=1}^n z_iR_i\otimes 
   V_i^*\right) x\right>\\
   &=\|x\|^2-\left\| \sum_{i=1}^n (z_i R_i\otimes V_i^*)x\right\|^2\\
   &=
   \left<\left[I\otimes \left(I-\sum_{i=1}^n 
   |z_i|^2V_iV_i^*\right)\right]x,x\right>\geq 0.
   \end{split}
   \end{equation*}
  Hence and using relation \eqref{roo}, we  deduce that
  \begin{equation}
  \label{Re}
  \text{\rm Re}\left\{ \left<\left[\left(1-\frac{2} {\rho}\right) I+\frac{2}
   {\rho}
  \left(I-\sum_{i=1}^n z_iR_i\otimes T_i^*\right)^{-1}\right]y,y\right>\right\}
  \geq 0
  \end{equation}
  for any $y\in F^2(H_n)\otimes \cH$.
  If $k\in F^2(H_n)\otimes \cH$, then 
  $y:= \left(I-\sum_{i=1}^n z_iR_i\otimes T_i^*\right)k$ is in 
  $F^2(H_n)\otimes \cH$ and inequality \eqref{Re} implies \eqref{ro-2}
   for any $z_i\in \DD$.
  Now,  note  that relation \eqref{ro-2}  is equivalent to the inequality
  \begin{equation}\label{2/ro}
  \left(\frac{2}{\rho} -1\right) \sum_{i=1}^n |z_i|^2 \|(R_i\otimes T_i^*)k\|^2+
  \left( 2-\frac{2}{\rho} \right)~
   \text{\rm Re}\left< \left(\sum_{i=1}^n z_i R_i\otimes T_i^*\right)k,k\right>
   \leq \|k\|^2
   \end{equation}
     for any $k\in F^2(H_n)\otimes \cH$ and $z_i\in \DD$, 
  $i=1,\ldots, n$.  Moreover, the latter inequality is equivalent to
   \eqref{1-ro}, which proves that  (ii) is equivalent to (iii).
  
  Now assume that condition (iii) holds.
    We show that the spectrum
   of the
  operator
  $\sum_{i=1}^n z_iR_i\otimes T_i^*$ is contained in $\DD$. Suppose that 
  there exists $\lambda_0\in \DD\backslash\{0\}$ such that $\lambda_0^{-1}$ is 
  in the boundary  of the spectrum
   of $\sum_{i=1}^n z_iR_i\otimes T_i^*$. Since  $\lambda_0^{-1}$ is in the 
   approximative spectrum, there exists a sequence
   $\{k_m\}\subset F^2(H_n)\otimes \cH$ such that $\|k_m\|=1$ for any $m=1,2,\ldots$, and 
   $$
   \lim_{m\to\infty}\left(I-\lambda_0\sum_{i=1}^n z_iR_i\otimes 
   T_i^*\right)k_m= 0.
   $$
   Hence, we have:
   \begin{enumerate}
   \item[(1)]
    $\lim\limits_{m\to \infty}\left<\lambda_0\left(\sum_{i=1}^n z_iR_i\otimes 
   T_i^*\right)k_m, k_m\right>= 1$;
   \item[(2)] $\lim\limits_{m\to \infty}\|\lambda_0\left(\sum_{i=1}^n z_iR_i\otimes 
   T_i^*\right)k_m\|=1$.
   \end{enumerate}
   Let $0<\gamma<\min \left\{\frac{1} {|\lambda_0||z_i|}-1:
    \ i=1,\ldots, n\right\}$ 
   and note that $w_i:=(\lambda_0+\gamma \lambda_0)z_i\in \DD$ for any $i=1,\ldots, n$.
   Since (ii)$\Leftrightarrow$(iii), we can use inequality \eqref{ro-2} to obtain 
   $$
   (\rho-2)\left\|\left(I-\sum_{i=1}^n 
    w_i R_i\otimes T_i^*\right)k_m\right\|^2 +2~
  \text{\rm Re}\left< \left(I-\sum_{i=1}^n 
   w_iR_i\otimes T_i^*\right)k_m,
  k_m\right>\geq 0
  $$
  for any $m=1,2,\ldots$. Taking $m\to\infty$ in this
   inequality and using relations $(1)$ and $(2)$, we get
   $(\rho-2)\gamma-2\geq2$. Letting $\gamma\to 0$, we get a contradiction.
   Therefore, the operator
   $I-\sum_{i=1}^n z_iR_i\otimes T_i^*$ is invertible for any $z_i\in \DD$, \ 
   $i=1,\ldots, n$.
   Using Lemma \ref{spec-ra}, we deduce that
   $$
   r(z_1T_1,\ldots, z_nT_n)=r\left(\sum_{i=1}^n z_iR_i\otimes T_i^*\right)<1
   $$
   and  the equality \eqref{inverse} holds.
     Let
    $y=(I-\sum_{i=1}^n z_iR_i\otimes T_i^*)k\in F^2(H_n)\otimes \cH$,
    where $k\in F^2(H_n)\otimes \cH$.
    The inequality \eqref{ro-2} implies 
    \begin{equation}\label{ro-inv}
    (\rho-2)\|y\|^2+\text{\rm Re}\left<y, 
    \left(I-\sum_{i=1}^n z_iR_i\otimes T_i^*\right)^{-1} y\right>\geq 0 
    \end{equation}
   for any $y\in F^2(H_n)\otimes \cH$ and $z_i\in \DD$, \ $i=1,\ldots, n$.
  It is clear that the inequality \eqref{ro-inv} and Lemma \ref{spec-ra}
  imply (iv). The implication (iv)$\implies$(v) is obvious.
  
  Assume now that condition (v) holds.  
  As in the proof of Theorem \ref{vN} (the implication (iii)$\implies$(i)),
   one can  show that
  \begin{equation*} 
  \left< \frac {1} {\rho}M_r\left( \sum_{|\beta|\leq q} e_\beta\otimes
   h_\beta \right), 
   \sum_{|\gamma|\leq q} e_\gamma\otimes h_\gamma\right> 
   = \sum_{|\beta|, |\gamma|\leq q}\left< K_{T,\rho,r}(\gamma, \beta)
    h_\beta, h_\gamma\right>,
  \end{equation*} 
  where
  the multi-Toeplitz kernel  $K_{T,\rho,r}:\FF_n^+\times \FF_n^+\to B(\cH)$
  is defined  by  
   $$
   K_{T,\rho,r}(\alpha, \beta):= 
   \begin{cases}
  \frac {1} {\rho} r^{|\beta\backslash \alpha|}  T_{\beta\backslash \alpha}
   &\text{ if } \beta>\alpha\\
   I  &\text{ if } \alpha=\beta\\
   \frac {1} {\rho} r^{|\alpha\backslash \beta|}(T_{\alpha\backslash \beta})^* 
     &\text{ if } \alpha>\beta\\
    0\quad &\text{ otherwise}.
   \end{cases} 
   $$ 
   Now, since $  M_r\geq 0$, we must have
  $[K_{T,\rho,r}(\alpha,\beta)]_{|\alpha|, |\beta|\leq q}\geq 0$ for any $0<r<1$
  and $q=0,1,\ldots$.
  Taking $r\to 1$, we  get
   $[K_{T,\rho}(\alpha,\beta)]_{|\alpha|, |\beta|\leq q}\geq 0$ for any 
    $q=0,1,\ldots$.
    Using the Naimark type dilation theorem of   \cite{Po-posdef},
     we find a row isometry $[V_1,\ldots, V_n]$  such that condition \eqref{ro}
     holds.   This proves the implication (v)$\implies$(i)
      and completes the proof
     of the theorem.
  \end{proof}

We remark that if $n\geq 2$, then $(T_1,\ldots, T_n)\in \cC_\rho$ does not imply that $(T_1^*,\ldots, T_n^*)\in \cC_\rho$. Given an arbitrary $n$-tuple of operators $(T_1,\ldots, T_n)$ and $\rho>0$, it is easy to see that there is $t>0$ such that $r\left(\frac{1}{t}T_1,\ldots, 
  \frac{1}{t}T_n\right)\leq 1$ and
  \begin{equation*}
  \begin{split}
  M_r\left(\frac{1}{t}T_1,\ldots, 
  \frac{1}{t}T_n\right)&\geq
  \rho I -2\sum_{k=1}^\infty\left(\frac{r} {t}\right)^k\left\| 
  \sum_{|\alpha|=k} R_\alpha\otimes T_{\tilde \alpha}^*\right\| I\\
  &
  \geq\rho I-2\sum_{k=1}^\infty\left(\frac{r} {t}\right)^k\left\| 
  \sum_{i=1}^n  T_i T_i^*\right\|^{k/2} I\geq 0,
  \end{split}
  \end{equation*}
  where the operator $M_r$ is defined by  \eqref{mr}.
  Consequently, Theorem \ref{rho} shows that  the $n$-tuple of operators
  $\left(\frac{1}{t}T_1,\ldots, 
  \frac{1}{t}T_n\right)$ is of class $ \cC_\rho$.
 Following the classical case (see \cite{Ho1}, \cite{W}), we
 define the map $\omega_\rho: B(\cH)^{(n)}\to [0,\infty)$ for  $\rho>0$,
  by setting
  \begin{equation}\label{om}
  \omega_\rho(T_1,\ldots, T_n):=\inf \left\{t:\ t>0,
   \left(\frac{1}{t}T_1,\ldots, 
  \frac{1}{t}T_n\right)\in \cC_\rho\right\}.
  \end{equation}
 When $\rho=\infty$, we set $\omega_\infty:=\lim\limits_{\rho\to \infty}
  \omega_\rho(T_1,\ldots, T_n)$.
  We present now some of the basic properties  of the joint $\rho$-operator radius.
  \begin{theorem} \label{more-propri}
  Let $(T_1,\ldots, T_n)$ be an $n$-tuple of operators $T_i\in B(\cH)$, 
   $i=1,\ldots, n$. The $\rho$-operator radius has the following properties:
   \begin{enumerate}
   \item[(i)]
   $
  \omega_\rho(zT_1,\ldots, zT_n)=|z| ~\omega_\rho(T_1,\ldots, T_n)
  $
  for any $z\in \CC;$
  \item[(ii)]
  $
  \|[T_1,\ldots, T_n]\|\leq \rho \omega_\rho(T_1,\ldots, T_n);
  $
  \item[(iii)]
  $\omega_\rho(T_1,\ldots, T_n)=0$ if and only if $T_1=\cdots= T_n=0;$
  \item[(iv)] 
  $(T_1\ldots, T_n)\in \cC_\rho$ if and only if
  $\omega_\rho(T_1,\ldots, T_n)\leq 1;$
  \item[(v)]
  $\omega_\infty(T_1,\ldots, T_n)=r(T_1,\ldots, T_n)\leq \omega_\rho(T_1,\ldots, T_n);$
  \item[(vi)] If $0<\rho\leq \rho'$, then $\cC_\rho\subseteq \cC_{\rho'}$
 and 
 $$\omega_{\rho'}(T_1,\ldots, T_n)\leq\omega_{\rho}(T_1,\ldots, T_n).
 $$
   \end{enumerate}
  \end{theorem}
  \begin{proof}
Part (i) follows imediately from Theorem \ref{rho} and 
  relation \eqref{om}.
  If $(\frac {1}{t}T_1,\ldots, \frac {1}{t}T_n)\in \cC_\rho$, then, 
  according to 
  \eqref{ro}, we have $\|[T_1,\ldots, T_n]\|\leq t\rho$.
  Consequently, using the definition \eqref{om}, we obtain the  inequality (ii).
  Part (iii) follows from (ii), and (iv) is a simple consequence of 
\eqref{om}.
  We prove now (v).
  If $t>0$ and  $(\frac {1}{t}T_1,\ldots, \frac {1}{t}T_n)\in \cC_\rho$,  then
  $r(T_1,\ldots, T_n)\leq t$. Hence, we deduce that
  \begin{equation}
  \label{r-om}
  r(T_1,\ldots, T_n)\leq \omega_\rho (T_1,\ldots, T_n).
  \end{equation}
  Now, we prove that if $r(T_1,\ldots, T_n)< 1$, 
  then there is $\rho>0$ such 
  $(T_1,\ldots, T_n)\in \cC_\rho$.
  Indeed, since the spectral radius $r(\cdot)$ is
   homogeneous, there is $t>1$ such that $r(tT_1,\ldots, tT_n)< 1$.
   Therefore, there is $M>0$ such that
   $$t^{2k} \left\| \sum_{|\alpha|=k} T_\alpha T_\alpha^*\right\|\leq M, 
   \quad \text{ for any } k=1,2\ldots.
   $$
  Hence, we deduce   that, for $0<r<1$, 
  \begin{equation*}
  \begin{split}
  M_r\left(T_1,\ldots, 
  T_n\right)&\geq
  \rho I -2\sum_{k=1}^\infty {r}^k\left\| 
  \sum_{|\alpha|=k} R_\alpha\otimes T_{\tilde \alpha}^*\right\| I\\
  &
  \geq\rho I-2\sum_{k=1}^\infty {r}^k\left\| 
  \sum_{i=1}^n  T_i T_i^*\right\|^{k/2} I\\
  &\geq \rho I-2\sum_{k=1}^\infty\left(\frac {r} {t}\right)^k\sqrt{M}I\geq 0
  \end{split}
  \end{equation*}
  for sufficiently large $\rho>0$. According to Theorem \ref{rho} part (v), we infer that 
   $(T_1,\ldots, T_n)\in \cC_\rho$.
  Now, if $r(T_1,\ldots, T_n)\neq 0$, let $\gamma>1$ and denote 
  $$X_i:=\frac{1} {\gamma r(T_1,\ldots, T_n)} T_i,\qquad i=1,\ldots, n.
  $$
  Since $r(X_1,\ldots, X_n)<1$, there is $\rho>0$ such that 
  $(X_1,\ldots, X_n)\in \cC_\rho$, i.e., $\omega_\rho(X_1,\ldots, X_n)\leq 1$.
  Hence, $\omega_\rho(T_1,\ldots, T_n)\leq \gamma r(T_1,\ldots, T_n)$, for any $\gamma>1$,
   which together with \eqref{r-om} imply (v).
   If $r(T_1,\ldots, T_n)=0$, then $r(mT_1,\ldots, mT_n)=0$ for any
    $m=1,2,\ldots$. Therefore, 
  $(mT_1,\ldots, mT_n)\in\cC_\rho$ for some $\rho>0$, i.e., 
  $\omega_\rho(T_1,\ldots, T_n)\leq \frac {1} {m}$. This implies 
  $$
  \lim_{\rho\to \infty}\omega_\rho(T_1,\ldots, T_n)=0.
  $$
  Part (vi) follows from the equivalence (i)$\Leftrightarrow$(ii) of Theorem \ref{rho}.
  The proof is complete.
  \end{proof}
\begin{corollary} $\bigcup\limits_{\rho>0}\cC_\rho$ is dense in the set of all 
power bounded 
$n$-tuples of operators in $B(\cH)$.
\end{corollary}
\begin{proof}
If $(T_1,\ldots, T_n)$ is power bounded $n$-tuple of operators, then $r(T_1,\ldots, T_n)\leq 1$ and 
therefore
$r(sT_1,\ldots, sT_n)<1$ for  $0\leq s<1$. According to 
the proof of Theorem \ref{more-propri} part (v), there is $\rho> 0$ such that 
$(sT_1,\ldots, sT_n)\in \cC_\rho$. This proves our assertion.
\end{proof}

\begin{corollary}
If $(A_1,\ldots, A_n), (B_1,\ldots, B_n)\in B(\cH)^{(n)}$  and $\rho\geq 1$,
then 
$$
\omega_\rho(A_i B_j:\ i,j=1,\ldots, n)\leq \rho^2\omega_\rho(A_1,\ldots, A_n)
\,\omega_\rho(B_1,\ldots, B_n).
$$
\end{corollary}
\begin{proof}
Since $\rho\geq 1$, we can use 
Theorem \ref{more-propri} part (ii) and (vi),  and deduce that
\begin{equation*}
\begin{split}
\omega_\rho(A_i B_j:\ i,j=1,\ldots, n)&\leq 
\omega_1(A_i B_j:\ i,j=1,\ldots, n)\\
&\leq \|[A_1,\ldots, A_n]\|\|B_1,\ldots, B_n]\|\\
&\leq\rho^2\omega_\rho(A_1,\ldots, A_n)\omega_\rho(B_1,\ldots, B_n).
\end{split}
\end{equation*}
\end{proof}

More properties of the $\rho$-operator radius are considered in the following

\begin{theorem}\label{admis}
Let $(T_1,\ldots, T_n)\in B(\cH)^{(n)}$ be an $n$-tuple of operators  
and let
$\omega_\rho: B(\cH)^{(n)}\to [0,\infty)$  be the $\rho$-operator radius.
 Then
\begin{enumerate}
\item[(i)]
$\omega_\rho (U^*T_1U,\ldots, U^*T_nU)= \omega_\rho(T_1,\ldots, T_n)$ 
for any unitary operator
  $U:\cK\to \cH;$
  \item [(ii)]
  $\omega_\rho(T_1|\cE,\ldots, T_n|\cE)\leq  \omega_\rho(T_1,\ldots, T_n)$ for any 
     subspace $\cE\subseteq \cH$ invariant under each  operator 
     $T_1,\ldots, T_n;$
 \item [(iii)]
 $\omega_\rho(I_\cG\otimes T_1,\ldots, I_\cG\otimes T_n)=
 \omega_\rho(T_1,\ldots, T_n)$ for 
 any separable
 Hilbert space $\cG$;
  \item [(iv)]
  
If $0<\rho\leq 2$, then  $\omega_\rho$ is a norm on $B(\cH)^{(n)}$ and 
 \begin{equation}\label{ine-iso}
 \omega_\rho (V^*T_1V,\ldots, V^*T_nV)\leq \omega_\rho(T_1,\ldots, T_n)
 \end{equation} 
for any  isometry
  $V:\cK\to \cH$.
 \end{enumerate}
\end{theorem}
\begin{proof}
According to Theorem \ref{rho},  if $(T_1,\ldots, T_n)\in C_\rho$,
then  the $n$-tuples $(U^*T_1U,\ldots, U^*T_nU)$, $(T_1|\cE,\ldots, T_n|\cE)$, and
$(I_\cG\otimes T_1,\ldots, I_\cG\otimes T_n)$ are of  class $C_\rho$.
Fix $(T_1,\ldots, T_n)\in B(\cH)$ and $\epsilon >0$. 
According to the definition of $\omega_\rho$, there exists $r>0$ such that
$$
r<\omega_\rho(T_1,\ldots, T_n)+\epsilon
$$
 and 
$(\frac {1} {r} T_1,\ldots, \frac {1} {r}T_n)\in C_\rho$.
 Using the above observation,
we have
$(\frac {1} {r} U^*T_1U,\ldots, \frac {1} {r}U^*T_nU)\in C_\rho$.
This shows that
$$
\omega_\rho (U^*T_1U,\ldots, U^*T_nU)\leq
 r<\omega_\rho(T_1,\ldots, T_n)+\epsilon.
 $$
 Taking $\epsilon \to 0$, we obtain
 \begin{equation}
\label{leq}
 \omega_\rho (U^*T_1U,\ldots, U^*T_nU)\leq
  \omega_\rho(T_1,\ldots, T_n).
 \end{equation}
Consequently, replacing $(T_1,\ldots, T_n)$ with $(U^*T_1U,\ldots, U^*T_nU)$,
 we obtain the reverse inequality. Therefore, the property (i) holds.
 The properties (ii) and (iii) can be proved in a similar manner. 

  To show that $\omega_\rho$ is a norm  if $0<\rho\leq 2$, it is enough
  to prove that $\cC_\rho$ is a convex body in $B(\cH)^{(n)}$.
  Let $(A_1,\ldots, A_n)$,  $(B_1,\ldots, B_n)$ be $n$-tuples of operators in
  $\cC_\rho$ and let $a,b\geq 0$ be such that $a+b=1$.
  Applying Theorem \ref{rho} (see the equivalence (i)$\Leftrightarrow$(ii))
  to each of the  $n$-tuples, and summing up, we get
  \begin{equation*}\begin{split}
  2\text{\rm Re}&\left<\left[I-\sum_{i=1}^n z_i R_i\otimes(aA_i^*+bB_i^*)\right] k,k
  \right>\\
  &\geq 
  (2-\rho) \left[a\left\|\left(I-\sum_{i=1}^n z_i R_i\otimes A_i^*\right)
   k\right\|^2+
   b\left\|\left(I-\sum_{i=1}^n z_i R_i\otimes B_i^*\right)
   k\right\|^2\right]
  \end{split}
  \end{equation*}
  for any $k\in F^2(H_n)\otimes \cH$.
  Since 
  $$a\|x\|^2+b\|y\|^2\geq \|ax+by\|^2, \qquad x,y\in F^2(H_n)\otimes \cH,
  $$
  $a+b=1$, and $2-\rho\geq 0$, the above inequality implies
  \begin{equation*}\begin{split}
  2\text{\rm Re}~&\left<\left[I-\sum_{i=1}^n z_i R_i\otimes(aA_i^*+bB_i^*)\right] k,k
  \right>\\
  &\geq 
  (2-\rho)\left\|I-\sum_{i=1}^n z_i R_i\otimes (aA_i+bB_i)^* 
   k\right\|^2
  \end{split}
  \end{equation*}
  for any $k\in F^2(H_n)\otimes \cH$.
  Using again Theorem \ref{rho} part (ii), we deduce that
  $$a(A_1,\ldots, A_n)+b(B_1,\ldots, B_n)\in \cC_\rho.
  $$
  
To prove the last part of the theorem,
note that, according to Theorem \ref{rho},    $(T_1,\ldots, T_n)\in C_\rho$ if
 and only if the operator
 $$F(T_1,\ldots, T_n):=
I-\left(1-\frac{1}{\rho}\right)\sum_{i=1}^n(\overline{z}_i R_i^*\otimes T_i+
   z_i R_i\otimes T_i^*)-
   \left(\frac{2}{\rho} -1\right) \sum_{i=1}^n |z_i|^2(I\otimes T_iT_i^*)
   $$
    is positive for any $z_i\in \DD$, 
  $i=1,\ldots, n$. 
 Notice that
 \begin{equation*}\begin{split}
 F(V^*T_1V,\ldots, V^*T_nV)=(I\otimes V^*)&\left[F(T_1,\ldots, T_n)\right.\\
 &+
 \left.\left(\frac{2}{\rho} -1\right) \sum_{i=1}^n |z_i|^2(I\otimes 
 T_i(I-VV^*)T_i^*)\right] (I\otimes V).
 \end{split}
 \end{equation*}
 Therefore, if $0<\rho\leq 2$, then $F(V^*T_1V,\ldots, V^*T_nV)\geq 0$
 for any $z_i\in \DD$, 
  $i=1,\ldots, n$. Using again Theorem \ref{rho}, we infer that
  $(V^*T_1V,\ldots, V^*T_nV)\in C_\rho$.
  Notice that the inequality \eqref{ine-iso} can be deduced similarly to  \eqref{leq}.
  The proof is complete.
\end{proof}

According to the noncommutative dilation theory  for row contractions,  
it is clear that 
  $$
  \omega_1(T_1,\ldots, T_n)=\left\|\sum_{i=1}^n T_iT_i^*\right\|^{1/2}
  =\omega_1(R_1\otimes T_1^*+\cdots + R_n\otimes T_n^*).
  $$

 In what follows we obtain  a generalization of this fact
 and more characterizations
for the $n$-tuples of operators of class $C_\rho$.

\begin{theorem}\label{ro-contr}
Let $T_1,\ldots, T_n\in B(\cH)$ and let $\cS\subset C^*(S_1,\ldots, S_n)$
  be the operator system defined by
  \begin{equation}\label{Sys1}
  \cS:=\{p(S_1,\ldots, S_n)+q(S_1,\ldots, S_n)^*:\ p,q\in \cP_n\}.
  \end{equation}
 Then the following statements are equivalent:
 \begin{enumerate}
 \item[(i)] $(T_1,\ldots, T_n)\in \cC_\rho;$  
  \item[(ii)]
$R_1\otimes T_1^*+\cdots+ R_n\otimes T_n^*\in \cC_\rho;$
\item[(iii)]
   The map $\Psi:\cS\to B(\cH)$ defined by
 \begin{equation*}
 \Psi(p(S_1,\ldots, S_n)+q(S_1,\ldots, S_n)^*):=
 p(T_1,\ldots, T_n)+q(T_1,\ldots,T_n)^*
 + (\rho-1)(p(0)+\overline{q(0)})I
 \end{equation*}
 is completely positive;
 \item[(iv)] The spectral radius $r(T_1,\ldots, T_n)\leq 1$ and
 $$
 \rho I+(1-\rho)r\sum_{i=1}^n (R_i^*\otimes T_i+ R_i\otimes T_i^*)+
 (\rho-2)r^2\left( I\otimes \sum_{i=1}^n T_iT_i^*\right)\geq 0
 $$
 for any $0<r<1$.
\end{enumerate}
Moreover, for any $T_1,\ldots, T_n\in B(\cH)$, 
$$
\omega_\rho(T_1,\ldots, T_n)=\omega_\rho
(R_1\otimes T_1^*+\cdots+ R_n\otimes T_n^*).
$$
\end{theorem}

\begin{proof}  
Assume that  $(T_1,\ldots, T_n)\in \cC_\rho$.   
In the particular case when 
$z_1=\cdots =z_n=z\in \DD$ and 
$X:= R_1\otimes T_1^*+\cdots+ R_n\otimes T_n^*$, the inequality \eqref{ro-2}
from Theorem \ref{rho}  implies
$$
(\rho-2)\|(I-zX)k\|^2+
2\,\text{\rm Re}\left< (I-zX)k,k\right>\geq 0
$$
for any $k\in F^2(H_n)\otimes \cH$ and $z\in \DD$.
According to \cite{SzF-book}, we deduce that $X\in C_\rho$.

Conversely, assume that $X\in C_\rho$.
Using again \cite{SzF-book}, we have $r(X)\leq 1$ and
$$
\sum_{k=1}^\infty r^k e^{-ik\theta} {X^*}^k +\rho I +
\sum_{k=1}^\infty r^k e^{ik\theta} {X}^k\geq 0
$$
for any $0<r<1$ and $\theta\in [0,2\pi]$.
In particular, when  $\theta=0$, we get the inequality
\begin{equation}
\label{CC}
\sum_{k=1}^\infty \sum_{|\alpha|=k} r^k R_\alpha^*\otimes T_{\tilde \alpha}+
\rho I + \sum_{k=1}^\infty \sum_{|\alpha|=k} r^k R_\alpha\otimes 
T_{\tilde \alpha}^*\geq 0
\end{equation}
 for any $0<r<1$. 
 On the other hand, Corollary \ref{pro} shows that
 $r(T_1,\ldots, T_n)=r(X)\leq 1$.
According to   Theorem \ref{rho}, the inequality \eqref{CC} implies
 $(T_1,\ldots, T_n)\in \cC_\rho$. Therefore,  
(i)$\Leftrightarrow$(ii).

Now, we prove that (i)$\Leftrightarrow$(iii). Assume that  
$(T_1,\ldots, T_n)\in \cC_\rho$. Using Theorem \ref{rho}, we have 
$r(T_1,\ldots, T_n)\leq 1$
 and the operator
$$\sum_{k=1}^\infty \sum_{|\alpha|=k}r^k 
  R_\alpha^*\otimes T_{\tilde \alpha} +\rho I+
  \sum_{k=1}^\infty \sum_{|\alpha|=k}r^k 
  R_\alpha\otimes T_{\tilde \alpha}^*
  $$
  is positive for any $0<r<1$, where the convergence is in the operator norm.
 This shows that, for any $q=0,1,\ldots,$ the operator
 $$
 P_{\cP_q\otimes \cH}\left(
 \sum_{1\leq |\alpha|\leq q}r^{|\alpha|} 
  R_\alpha^*\otimes T_{\tilde \alpha} +\rho I+
    \sum_{1\leq |\alpha|\leq q}r^{|\alpha} 
  R_\alpha\otimes T_{\tilde \alpha}^*\right)|\cP_q\otimes \cH
  $$
  is positive. Taking $r\to 1$, we infer that 
  $$
  M_{\rho, q}:=P_{\cP_q\otimes \cH}\left(
 \sum_{1\leq |\alpha|\leq q}  
  R_\alpha^*\otimes T_{\tilde \alpha} +\rho I+
    \sum_{1\leq |\alpha|\leq q}  
  R_\alpha\otimes T_{\tilde \alpha}^*\right)|\cP_q\otimes \cH
  $$
  is positive. Now, the proof is similar to that of the implication (i)$\implies$(ii) of Theorem
  \ref{C2}, where the multi-Toeplitz kernel
  $K_{T,2}$ is replaced by the kernel
$$
   K_{T,\rho}(\alpha, \beta):= 
   \begin{cases}
   T_{\beta\backslash \alpha} &\text{ if } \beta>\alpha\\
   \rho I  &\text{ if } \alpha=\beta\\
    T_{\alpha\backslash \beta}^*  &\text{ if } \alpha>\beta\\
    0\quad &\text{ otherwise}.
   \end{cases} 
   $$  
  We leave the details to the reader.
  Conversely, assume condition (iii) holds. As in the proof of Theorem
  \ref{C2} (implication (ii)$\implies$(i)) one can show that the kernel
  $K_{T,\rho}$ is positive. Using the Naimark type dilation theorem of \cite{Po-posdef},
  we deduce that $[T_1,\ldots, T_n]\in C_\rho$.
  
  Now, it remains to prove that (iv)$\Leftrightarrow$(i).
  Note that the operator $M_r$ of Theorem \ref{rho}  satisfies the equalities
  \begin{equation*}
  \begin{split}
  M_r&=
  \left(I-\sum_{i=1}^n r R_i^*\otimes T_i\right)^{-1} +(\rho-2)I +
  \left(I-\sum_{i=1}^n r R_i\otimes T_i^*\right)^{-1}\\
  &=\left(I-\sum_{i=1}^n r R_i^*\otimes T_i\right)^{-1}\Gamma_{\rho, r} 
  \left(I-\sum_{i=1}^n r R_i\otimes T_i^*\right)^{-1},
  \end{split}
  \end{equation*}
  where
  $$\Gamma_{\rho, r}:=
 \rho I+(1-\rho)r\sum_{i=1}^n (R_i^*\otimes T_i+ S_i\otimes T_i^*)+
 (\rho-2)r^2\left( I\otimes \sum_{i=1}^n T_iT_i^*\right). 
 $$
 This shows that the condition (v) of Theorem \ref{rho} is equivalent to 
 $r(T_1,\ldots, T_n)\leq 1$ and $\Gamma_{\rho, r}\geq 0$ for any
 $0<r<1$.
Therefore, (iv)$\Leftrightarrow$(i).
 Using the equivalence (i)$\Leftrightarrow$(ii),
  we deduce the last part of the theorem.
 The proof is complete.
 \end{proof}

\begin{corollary}\label{num-ra}
Let $(T_1,\ldots, T_n)$ be an $n$-tuple of operators and  
let $w(T_1,\ldots, T_n)$ be the joint numerical
    radius. 
 The following statements are equivalent:
  \begin{enumerate}
  \item[(i)] $T\in \cC_2;$
  \item[(ii)] $w(T_1,\ldots, T_n)\leq 1;$
  \item[(iii)] $r(T_1,\ldots, T_n)\leq 1$ and
  $\sum\limits_{i=1}^n (R_i^*\otimes T_i+R_i\otimes T_i^*)\leq 2I;$
  \item[(iv)] \text{\rm Re}\,$\left<\left(I-\sum\limits_{i=1}^n z_iR_i\otimes T_i^*\right)
   k,  k\right>\geq 0$ for any $k\in F^2(H_n)\otimes \cH$ and
    $z_i\in {\DD}$.
  \item[(v)] $\sum_{i=1}^n |\left<(R_i\otimes T_i^*)k, k\right>|\leq \|k\|^2$
  for any $k\in F^2(H_n)\otimes \cH$.
  \end{enumerate}
 Moreover,  
 $
   \omega_2(T_1,\ldots, T_n)=
    w(T_1,\ldots, T_n).
   $ 
 \end{corollary}
  
  \begin{proof}
  Theorem \ref{berger} shows that (i)$\Leftrightarrow$(ii).
  Consequently, using   Theorem
\ref{more-propri} part (iv)
  \ref{ro-contr} (in the particular case $\rho=2$) and the homogeneity
  of $\omega_2$ and $w$, we deduce that 
  $$
   \omega_2(T_1,\ldots, T_n)=
    w(T_1,\ldots, T_n).
   $$ 
  According to Theorem \ref{ro-contr}, $(T_1,\ldots, T_n)\in \cC_2$ if 
  and only if $r(T_1,\ldots, T_n)\leq 1$ and 
  $$t\sum_{i=1}^n (R_i^*\otimes T_i + R_i\otimes T_i^*)\leq 2I
  $$
  for any $0<t<1$. This proves that (i)$\Leftrightarrow$(iii).
  The equivalence (i)$\Leftrightarrow$(iv) is a particular case of Theorem
  \ref{rho} ( see (i)$\Leftrightarrow$(ii)).
  The equivalence (iv)$\Leftrightarrow$(v) is  obvious.
  The proof is complete.
  \end{proof}

 As we have seen in this  paper, 
an important role in  multivariable operator theory for  
 $n$-tuples of operators
 $T:=(T_1,\ldots, T_n)$ is played by the multi-analytic
  operator
 $$
 R_1\otimes T_1^*+\cdots +R_n\otimes T_n^*.
 $$
 We can prove that  this operator is  a complete  
 unitary invariant.
 More precisely, if
  $T':=(T_1',\ldots, T_n')$, \ $T_i'\in B(\cH')$, is another
  $n$-tuple of operators,  
   one can prove that $T$ is unitarily equivalent  
  to $T'$ if and only if the
   multi-analytic operators
  $R_1\otimes T_1^*+\cdots +R_n\otimes T_n^*$ and 
  $R_1\otimes T_1'^*+\cdots +R_n\otimes T_n'^*$ are unitarily equivalent, i.e.,
  there exists a unitary  multi-analytic operator
  $W$ from $F^2(H_n)\otimes \cH$ to $F^2(H_n)\otimes \cH'$ such that
  $$
  W(R_1\otimes T_1^*+\cdots +R_n\otimes T_n^*)= 
  (R_1\otimes T_1'^*+\cdots +R_n\otimes T_n'^*)W.
  $$
  One implication is clear, so assume  
   that  the above equality holds.
   Since $W$ is a unitary multi-analytic operator,
   we must have $W=I\otimes U^*$ for some unitary 
  $U$ operator from $\cH$ to $\cH'$ (see \cite{Po-charact}).
   Hence, and using the fact that 
   $R_1,\ldots, R_n$ are  the right creation operators on $F^2(H_n)$, we deduce 
   that $UT_i=  T_i' U$ for any $i=1,\ldots, n$.

The proof of the following result is similar to that of Corollary \ref{simi}.
  We shall omit  it.
\begin{theorem}\label{BKS-ro}
Let $(T_1,\ldots, T_n)$ be an $n$-tuple of operators
 $T_i\in B(\cH)$ with $\omega_\rho(T_1,\ldots, T_n)\leq 1$.
 Then there is a completely bounded map
 $\Phi:\cA_n+\cA_n^*\to  B(\cH)$ such that
 $$
 \Phi(p(S_1,\ldots, S_n)+q(S_1,\ldots, S_n)^*)=p(T_1,\ldots, T_n)+q(T_1,\ldots, T_n)
$$
for any $p(S_1,\ldots, S_n), q(S_1,\ldots, S_n)\in \cP$, 
and $\|\Phi\|_{cp}\leq |1-\rho|+\rho$.
Consequently, the $n$-tuple $(T_1,\ldots, T_n)$ is simultaneously similar to 
a row contraction.
\end{theorem}

\begin{corollary}
If $R_1\otimes T_1^*+\cdots +R_n \otimes T_n^*$ has the spectrum in the open
 unit disc,
 then $(T_1,\ldots, T_n)$ is simultaneously similar to 
a row contraction.
\end{corollary}
\begin{proof}
Since 
$$r(T_1,\ldots, T_n)=r(R_1\otimes T_1^*+\cdots +R_n \otimes T_n^*)<1,$$
 the condition (v) of Theorem \ref{rho} holds for some large $\rho>0$.
 Now, the result follows from Theorem \ref{BKS-ro}.
\end{proof}

Now, as in the case of the joint numerical range ($\rho=2$), 
Theorem \ref{BKS-ro} can be used to show that if $(T_1,\ldots, T_n)\in \cC_\rho$, then
$$f(T_1,\ldots, T_n)=\lim_{r\to 1} f(rT_1,\ldots, rT_n),
$$
where $f(S_1,\ldots, S_n)\in \cA_n$ and the limit exists in the operator norm.
Moreover,
$$
f(T_1,\ldots, T_n)+(\rho-1) f(0)I=\rho P_\cH f(V_1,\ldots, V_n)|\cH,
$$
where $V_1,\ldots V_n$  are the isometries from \eqref{ro}. 
Our multivariable  Berger-Kato-Stampfli type theorem for the operator radius $\omega_\rho$ is 
the following. The proof is similar to the proofs of Theorem \ref{BKS1} and 
Theorem \ref{Real-ine}, but uses Theorem \ref{ro-contr}.

\begin{theorem}

If $\omega_\rho(T_1,\ldots, T_n)\leq 1$ and $F_j \in M_m\otimes \cA_n$ with  $F_j(0)=0$, $j=1,\ldots, k$, then
$$\omega_\rho(F_1(T_1,\ldots, T_n),
\ldots, F_k(T_1,\ldots, T_n))\leq \|[F_1,\ldots, F_k]\|.$$
Moreover, if $\omega_\rho(T_1,\ldots, T_n)\leq 1$ and $f\in \cA_n$ with 
$\text{\rm Re}\, f\geq 0$, then
$$
\text{\rm Re}\,W( f(T_1,\ldots, T_n))\geq (1-\rho)\text{\rm Re}\, f(0).
$$

\end{theorem}

We can now provide the following power inequality for the $\rho$-operator radius.

\begin{corollary}
If $(T_1,\ldots, T_n)\in B(\cH)^{(n)}$ and $\rho\in (0,\infty)$,
then
$$
\omega_\rho(T_\alpha:\ |\alpha|=k)\leq \omega_\rho(T_1,\ldots, T_n)^k.
$$
Moreover, $\omega_\infty(T_\alpha:\ |\alpha|=k)=
\omega_\infty(T_1,\ldots, T_n)^k$.
\end{corollary}
\begin{proof}
Since ~$\omega_\rho$ is homogeneous, we can assume that
$\omega_\rho(T_1,\ldots, T_n)=1$. Hence, $(T_1,\ldots, T_n)\in \cC_\rho$ and 
therefore there is a Hilbert space $\cK\supseteq \cH$ and
$V_1,\ldots, V_n\in B(\cK)$ isometries with orthogonal ranges such that 
$T_\beta=\rho P_\cH V_\beta|\cH$ for any $\beta\in \FF_n^+$.
Notice that the operators $V_\alpha$, $|\alpha|=k$,  are isometries with
 orthogonal ranges.  Now, the above equation  implies  that
  $(T_\alpha:\ |\alpha|=k)\in \cC_\rho$,
 which is equivalent to
 $\omega_\rho(T_\alpha:\ |\alpha|=k)\leq 1$.
 According to Theorem \ref{more-propri} part (v), we have
 \begin{equation*}
 \begin{split}
 \omega_\infty(T_\alpha:\ |\alpha|=k)&=r(T_\alpha:\ |\alpha|=k)\\
 &=r(T_1,\ldots, T_n)^k=\omega_\infty(T_1,\ldots, T_n)^k.
 \end{split}
 \end{equation*}
The proof is complete.
\end{proof}

In the commutative case, we can obtain the following 
characterization  for $n$-tuples of operators of class $\cC_\rho$.

  \begin{theorem}\label{ro-C2-commut}
  Let $T_1,\ldots, T_n\in B(\cH)$ be commuting operators and let 
  $\cS_c\subset C^*(B_1,\ldots, B_n)$
  be the operator system defined by
  \begin{equation}\label{ro-Sys-c}
  \cS_c:=\{p(B_1,\ldots, B_n)+q(B_1,\ldots, B_n)^*:\ p,q\in \cP\},
  \end{equation}
  where $B_1,\ldots, B_n$ are the creation operators on the symmetric Fock space.
  Then the following statements are equivalent:
  \begin{enumerate}
 \item[(i)] 
  $(T_1,\ldots, T_n)\in \cC_\rho$;
 \item[(ii)] The map $\Psi_c:\cS_c\to B(\cH)$ defined by
 \begin{equation*} 
 \Psi_c(p(B_1,\ldots, B_n)+q(B_1,\ldots, B_n)^*):=
 p(T_1,\ldots, T_n)+q(T_1,\ldots,T_n)^*
 +(\rho-1) (p(0)+\overline{q(0)})I
 \end{equation*}
 is completely positive.
 \item[(iii)] There is a Hilbert space $\cG\supseteq \cH$ and a $*$-representation
    $\pi:C^*(B_1,\ldots, B_n)\to B(\cG)$ such that 
    $$T_\alpha =\rho P_\cH \pi(B_\alpha)| \cH\quad \text{ for any }  
    \alpha\in \FF_n^+\backslash \{g_0\}.
    $$
  \end{enumerate}
  \end{theorem}
  \begin{proof}
 Assume that $(T_1,\ldots, T_n)\in \cC_\rho$. According the the proof of
  Theorem \ref{rho}, 
  the multi-Toeplitz  kernel 
   $ K_{T,\rho}:\FF_n^+\times \FF_n^+\to B(\cH)$
 defined
   by  
   $$
    K_{T,\rho}(\alpha, \beta):= 
   \begin{cases}
   \frac{1} {\rho}T_{\beta\backslash \alpha} &\text{ if } \beta>\alpha\\
   I  &\text{ if } \alpha=\beta\\
    \frac{1} {\rho}T_{\alpha\backslash \beta}^*  &\text{ if } \alpha>\beta\\
    0\quad &\text{ otherwise}.
   \end{cases} 
   $$
  is positive definite.  Since $T_1,\ldots, T_n$ are commuting operators, 
  according to Theorem 3.3 of \cite{Po-moment}, there exists a completely
   positive linear map
  $\Phi:C^*(B_1,\ldots, B_n)\to B(\cH)$ such that $\Phi(I)=I$ and
  $$
  \Phi(B_\sigma)=\frac {1} {\rho} T_\sigma, \quad
   \sigma\in \FF_n^+\backslash\{g_0\}.
  $$
  Consequently, using  Stinespring's theorem,  there exists a Hilbert
   space $\cG\supseteq \cH$ 
  and a $*$-representation $\pi:C^*(B_1,\ldots, B_n)\to B(\cH)$ such  
  that (iii) holds.
  
  Assume now that condition  (iii) holds and let $W_i:=\pi(B_i)$, $i=1,\ldots, n$. Since 
   $[W_1,\ldots, W_n]$
  is a row contraction, there exists a Hilbert space $\cK\supseteq \cG$
   and $V_1,\ldots, V_n\in B(\cK)$
  isometries with orthogonal ranges such that
  $W_\alpha^*=V_\alpha^*|\cK$, $\alpha\in \FF_n^+$. Since $\cH\subseteq \cG$,
   (iii) implies
  $T_\alpha=\rho P_\cH V_\alpha|\cH$, $\alpha\in \FF_n^+$. Using Theorem \ref{rho}, we deduce that
   $(T_1,\ldots, T_n)\in \cC_\rho$,
  therefore (iii)$\implies$(i). 
  Now we show that 
(iii)$\implies$(ii).
  It is clear that (iii)  implies
  \begin{equation*}
  \begin{split}
  p(T_1,\ldots, T_n)+ q(T_1,\ldots, T_n)^*&+(\rho-1)(p(0)+
  \overline{q(0)})I\\
  &=\rho P_\cH\pi[p(B_1, \ldots, B_n)+
  q(B_1, \ldots, B_n)^*]|\cH.
  \end{split}
  \end{equation*}
 Hence, we deduce  (ii).
 To prove that  (ii)$\implies$(i), note that there is a completely positive linear map
 $\gamma:C^*(S_1,\ldots, S_n)\to C^*(B_1,\ldots, B_n)$ such that 
 $\gamma(S_\alpha S_\beta^*)=B_\alpha B_\beta^*$
 (see \cite{Po-poisson}).
 Since the map $\Psi_c\circ \gamma$ is completely positive and 
 $(\Psi_c\circ \gamma)|\cS=\Psi$,  where $\cS$ and $\Psi$ are defined in Theorem 
 \ref{ro-contr},
 the latter theorem  implies $(T_1,\ldots, T_n)\in \cC_\rho$.
 The proof is complete.
  \end{proof}

The next result contains more properties for the $\rho$-operator radius.

\begin{proposition}\label{om-ineq}
\begin{enumerate}
\item[(i)]
Let $0<\rho\leq \rho'$ and $(T_1,\ldots, T_n)\in B(\cH)^{(n)}$. Then
\begin{equation}\label{ro'}
\omega_{\rho'}(T_1,\ldots, T_n)\leq \omega_{\rho}(T_1,\ldots, T_n)\leq 
\left(\frac {2\rho'} {\rho}-1\right) \omega_{\rho'}(T_1,\ldots, T_n).
\end{equation}
In particular, the map $\rho\mapsto \omega_{\rho}(T_1,\ldots, T_n)$ is 
continuous on $(0,\infty)$.
\item[(ii)]
The operator radius $\omega_\rho:B(\cH)^{(n)}\to [0,\infty)$ is continuous in the norm topology of $B(\cH)^{(n)}$.
\item[(iii)] For any $\rho>0$,
$$\omega_\rho(T_1,\ldots, T_n)\geq 
  \max\left\{ 1,\frac {2} {\rho}-1\right\} r(T_1,\ldots, T_n). 
$$
\end{enumerate}
\end{proposition}
\begin{proof} It is well-known  (see \cite{Ho1}) that for any $X\in B(\cH)$,
$$
\omega_{\rho'}(X)\leq \omega_{\rho}(X)\leq 
\left(\frac {2\rho'} {\rho}-1\right) \omega_{\rho'}(X).
$$
Applying this result to the operator
 $X:=R_1\otimes T_1^*+\cdots +R_n\otimes T_n^*$
and taking into account that
$$\omega_\rho(T_1,\ldots, T_n)=
\omega_\rho (R_1\otimes T_1^*+\cdots +R_n\otimes T_n^*)
$$
(see Theorem \ref{ro-contr}),  we deduce \eqref{ro'}.
The continuity of the map $\rho\mapsto \omega_{\rho}(T_1,\ldots, T_n)$ is now obvious.
 According to  \eqref{ro'}, we have 
\begin{equation*}\|[T_1,\ldots, T_n]\|\leq \omega_\rho(T_1,\ldots, T_n)\leq
\left(\frac {2} {\rho}-1\right)\|[T_1,\ldots, T_n]\| \end{equation*} if $0<\rho\leq 1$, and \begin{equation*} \frac{1}{2\rho-1} \|[T_1,\ldots, T_n]\|\leq \omega_\rho(T_1,\ldots, T_n)\leq \|[T_1,\ldots, T_n]\|\end{equation*}if $\rho>1$. Hence, the continuity of $\omega_\rho$ follows.
To prove (iii), we recall from \cite{Ho1}  the inequality  $w_\rho(X)\geq \max\left\{ 1,\frac {2} {\rho}-1\right\} r(X)$  for  $X\in B(\cH) $. Applying this result to the operator $X $ and 
using Theorem \ref{ro-contr}, we deduce the inequality (iii).
The proof is complete. \end{proof}

In what follows we calculate the $\rho$-operator radius for a class of row contractions.

\begin{example}   Let $(T_1,\ldots, T_n)\in B(\cH)^{(n)}$ be an $n$-tuple of
 operators  with   $\|[T_1,\ldots, T_n]\|=1$ and
$T_\alpha=0$  for any $\alpha\in \FF_n^+$ with  $|\alpha|=2$. 
 Then
$$\omega_{\rho}(T_1,\ldots, T_n)=\frac{1}{\rho}, \quad \text{ for }  \rho\in (0,\infty),
$$
and $\omega_\infty(T_1,\ldots, T_n)=0$.
Indeed,  since $[T_1,\ldots, T_n]$ is a row contraction, there is a Hilbert
 space 
$\cK\supseteq \cH$ and $V_1,\ldots, V_n\in B(\cK)$ isometries with orthogonal
 ranges
such that $T_\beta=P_\cH V_\beta|\cH$ for any $\beta\in \FF_n^+$. 
Set $X_i:=\rho T_i$, \ $i=1,\ldots, n$, and note that
 the condition
$T_\alpha=0$   if   $|\alpha|=2$
 implies

$$X_\beta=\rho T_\beta=\rho P_\cH V_\beta|\cH$$for any $\beta\in \FF_n^+\backslash \{g_0\}$.
Therefore, $(X_1,\ldots, X_n)\in \cC_\rho$, i.e., 
$\omega_\rho(X_1,\ldots, X_n)\leq 1$, which implies 
$\omega_{\rho}(T_1,\ldots, T_n)\leq \frac {1} {\rho}$.
The reverse inequality follows from Theorem $\ref{more-propri}$ part (ii).  Now, it is clear that $$\omega_\infty(T_1,\ldots, T_n)=\lim_{\rho\to\infty}\omega_\rho(T_1,\ldots, T_n)=0.$$ 

\end{example}

\smallskip

  \section*{Part II. Joint operator radii, inequalities, and applications}
  \label{JOR}

  \section{von Neumann inequalities}
  \label{von}
 
 In this section, we prove  
     von Neumann type inequalities \cite{vN} 
     for arbitrary admissible or strongly admissible operator radii
     $\omega:B(\cH)^{(k)}\to [0,\infty)$.  
      We show  that,   given a row contraction
      $[T_1,\ldots, T_n]$,  an
        inequality of the form
      $$\omega(f_1(T_1,\ldots, T_n),\ldots, f_k(T_1,\ldots, T_n))\leq
      \omega(f_1(S_1,\ldots, S_n),\ldots, f_k(S_1,\ldots, S_n))
      $$
      holds if $f_1(S_1,\ldots, S_n),\ldots, f_k(S_1,\ldots, S_n)$ belong
       to   operator algebras 
      or operator systems generated by the left creation operators
       $S_1,\ldots, S_n$ and the identity such as
      the noncommutative disc algebra $\cA_n$,  the Toeplitz $C^*$-algebra
       $C^*(S_1,\ldots, S_n)$,  the noncommutative analytic Toeplitz algebra
        $F_n^\infty$, 
        the noncommutative Douglas type algebra ${\cD}_n$, and the operator system
         \ $F_n^\infty (F_n^\infty)^*$. 
        The operator  $f_j(T_1,\ldots, T_n)$ 
       is defined by an appropriate functional calculus for row contractions.
       Actually, we obtain matrix-valued generalizations of the above inequality.

Let $\omega:=\{\omega_{\cH,n}\}$ ($\cH$ is any separable Hilbert space and 
$n=1,2,\ldots$) be a family of mappings 
 $\omega_{\cH, n}: B(\cH)^{(n)}\to [0,\infty)$.
 If $(T_1,\ldots, T_n)\in B(\cH)^{(n)}$ we denote
 $\omega(T_1,\ldots, T_n):=\omega_{\cH,n}(T_1,\ldots, T_n)$.
 We say that $\omega$ is a {\it joint operator radius} if,
  for any 
 $(T_1,\ldots, T_n)\in B(\cH)^{(n)}$,  
 \begin{enumerate}
 \item[(i)] $\omega(T_1,\ldots, T_n)=\omega(T_{\sigma(1)},\ldots, T_{\sigma(n)})$
 for any permutation $\sigma:\{1,\ldots, n\}\to \{1,\ldots, n\}$;
 \item[(ii)]
 $\omega(T_1,\ldots, T_n,0)=\omega(T_1,\ldots, T_n)$ .
 \end{enumerate}
A joint operator radius $\omega$ is called {\it admissible} if, for 
any 
 $(T_1,\ldots, T_n)\in B(\cH)^{(n)}$,  it satisfies 
the following 
properties:
\begin{enumerate}
\item[(iii)] $\omega(U^*T_1U,\ldots, U^*T_n U)=\omega(T_1,\ldots, T_n)$
for any 
  unitary operator $U:\cK\to \cH$;
 \item[(iv)]
$\omega(I_\cG\otimes T_1,\ldots,I_\cG\otimes  T_n)=\omega(T_1,\ldots, T_n)$
for  any separable Hilbert space $\cG$;
\item[(v)] $\omega(T_1|\cM,\ldots, T_n|\cM)\leq \omega(T_1,\ldots, T_n)$
 for any invariant
  subspace $\cM\subseteq \cH$   under each operator $T_1,\ldots T_n$.  
\end{enumerate}

A joint operator radius $\omega$ is called {\it strongly admissible} if it satisfies 
the  
properties (iii), (iv), and 
\begin{enumerate}
\item[(v$'$)]
$\omega(P_\cM T_1|\cM,\ldots, P_\cM T_n|\cM)\leq \omega(T_1,\ldots, T_n)$
for any $(T_1,\ldots, T_n)\in B(\cH)^{(n)}$  and any
 closed subspace $\cM\subseteq \cH$.
\end{enumerate}

We remark that one can easily prove that conditions (iii) and (v$'$) hold 
 if and only if 
  $$\omega(V^*T_1V,\ldots, V^*T_n V)\leq \omega(T_1,\ldots, T_n)$$
for any 
 $(T_1,\ldots, T_n)\in B(\cH)^{(n)}$ and any  isometry $V:\cK\to \cH$.
Obviously, any strongly admissible radius is admissible.

We say that a 
 joint operator radius $\omega:B(\cK)^{(n)}\to [0,\infty)$ 
  is norm (resp. strongly, $*$-strongly, weakly) {\it quasi-continuous} if,  
 for any $b\geq 0$, the ball
 $$B_{\omega, b}:=\left\{ (T_1,\ldots, T_n)\in B(\cK)^{(n)}:\ 
 \omega (T_1,\ldots, T_n)\leq b\right\}
 $$
 is closed  in the product topology of $B(\cK)^{(n)}$, where $B(\cK)$  
 is endowed with 
 the norm (resp. strong, $*$-strong,  weak) operator topology.

\begin{proposition}\label{quasi} 

\begin{enumerate}

\item[(i)] The joint operator radii $\|\cdot\|$, $\|\cdot\|_e$, $w(\cdot)$, $w_e(\cdot)$, $r(\cdot)$, $r_e(\cdot)$, and $\omega_\rho(\cdot)$,  \ $0<\rho\leq 2$,  are strongly admissible;

\item[(ii)] The joint operator radii $\|\cdot\|$, $\|\cdot\|_e$, $w(\cdot)$, and $w_e(\cdot)$ are   WOT quasi-continuous;

\item[(iii)] The joint operator radius $\omega_\rho(\cdot)$, $(\rho>0)$, is admissible and $*$-SOT quasi-continuous.  

\end{enumerate}

\end{proposition}

\begin{proof}

Part (i) of the theorem concerning $w(\cdot)$, $w_e(\cdot)$, and 
$\omega_\rho(\cdot)$,  \ $0<\rho\leq 2$, follows from Theorem \ref{propri},
 Theorem \ref{propri2}, and  Theorem \ref{admis}, respectively.
To prove that $\|\cdot\|$, $\|\cdot\|_e$, $r(\cdot)$, and  $r_e(\cdot)$
are strongly admissible is a simple exercise, so we leave it to the reader.
Now we prove (ii).
For each $i=1,\ldots, n$, let $\{T_i^{(\alpha)}\}_{\alpha\in J}$ be a net
 of operators in $B(\cH)$ which is WOT convergent to $T_i\in B(\cH)$.
 If $\|[T_1^{(\alpha)},\ldots, T_n^{(\alpha)}]\|\leq 1$ for any $\alpha\in J$, then
 \begin{equation*} 
 \left|\left< \sum_{i=1}^n T_i^{(\alpha)}h_i, k\right>\right| 
 \leq \|(h_1,\ldots, h_n)\| \|k\|
 \end{equation*}
 for any $h_1,\ldots, h_n, k\in \cH$.
This implies that $[T_1,\ldots, T_n]$ is a contraction,
 which proves that $\|\cdot\|$ is WOT quasi-continuous.
  Similarly one can prove that $\|\cdot\|_e$ is WOT quasi-continuous.
  
  Now, assume that $w(T_1^{(\alpha)},\ldots, T_n^{(\alpha)})\leq 1$
  for any $\alpha\in J$. According to Theorem \ref{propri} part (vi), we deduce that
  $\|T_i^{(\alpha)}\|\leq 2$. Since the net $\{T_i^{(\alpha)}\}_{\alpha\in J}$
  is WOT convergent to $T_i$, it is also $w^*$-convergent to $T_i$. This is equivalent to the WOT
  convergence of the net $\{I\otimes T_i^{(\alpha)}\}_{\alpha\in J}$ to $I\otimes T_i$.
  Therefore, the net $\{R_i\otimes {T_i^{(\alpha)}}^*\}_{\alpha\in J}$ is WOT
  convergent to $R_i\otimes T_i^*$ for $i=1,\ldots, n$.
  Using  Corollary \ref{num-ra}, we deduce that
  $w(T_1,\ldots, T_n)\leq 1$, which proves that the joint numerical range 
  is WOT
  quasi-continuous. Similarly, one can prove that the euclidean numerical
  radius   is WOT
  quasi-continuous.
  
  To prove part (iii),  notice first that  $\omega_\rho(\cdot)$  is admissible due to  Theorem
  \ref{admis}. Now, assume that the net 
  $\{{T_i^{(\alpha)}}^*\}_{\alpha\in J}$
  is SOT convergent to $T_i^*$ for each $i=1,\ldots, n$,  and
   $\omega_\rho(T_1^{(\alpha)}, \ldots, T_n^{(\alpha)})\leq 1$.
   According to Theorem \ref{more-propri}, we have $\|T_i^{(\alpha)}\|\leq \rho$.
   Since the map $A\mapsto I\otimes A$ is SOT-continuous on bounded subsets 
   of $B(\cH)$, we can use Theorem \ref{rho} (the equivalence (i)$\Leftrightarrow$(ii)) to deduce that 
   $(T_1,\ldots, T_n)\in \cC_\rho$. Consequently,   Theorem \ref{more-propri}
   implies $\omega_\rho(T_1,\ldots, T_n)\leq 1$, which completes the proof.\end{proof}

We need to recall from \cite{Po-poisson} a few facts about noncommutative
Poisson transforms associated with row contractions $T:=[T_1,\ldots, T_n]$,
\ $T_i\in B(\cH)$. For  each $0<r\leq 1$, define the defect operator
$\Delta_r:=(I-r^2T_1T_1^*-\cdots -r^2 T_nT_n^*)^{1/2}$.
The Poisson  kernel associated with $T$ is the family of operators
$K_{T,r} :\cH\to F^2(H_n)\otimes \cH$, \ $0<r\leq 1$, defined by
$$
K_{T,r}h:= \sum_{k=0}^\infty \sum_{|\alpha|=k} e_\alpha\otimes r^{|\alpha|} 
\Delta_r T_\alpha^*h,\quad h\in \cH.
$$
When $r=1$, we denote $\Delta:=\Delta_1$ and $K_T:=K_{T,1}$.
The operators $K_{T,r}$ are isometries if $0<r<1$ and 
$$
K_T^*K_T=I-
\text{\rm SOT-}\lim_{k\to\infty} \sum_{|\alpha|=k} T_\alpha T_\alpha^*.
$$
This shows that $K_T$ is an isometry if and only if $T$ is a $C_0$-row
 contraction (\cite{Po-isometric}),
i.e., 
$$
\text{\rm SOT-}\lim_{k\to\infty} \sum_{|\alpha|=k} T_\alpha T_\alpha^*=0.
$$
A key property of the Poisson kernel
is that
\begin{equation}\label{eq-ker}
K_{T,r}(r^{|\alpha|} T_\alpha^*)=(S_\alpha^*\otimes I)K_{T,r}
\end{equation}
for any $0<r<1$ and $\alpha\in \FF_n^+$. If $T$ is $C_0$-row
 contraction, then equality \eqref{eq-ker} holds also for $r=1$.
 In \cite{Po-poisson}, we introduced  the Poisson transform associated with 
 $T:=[T_1,\ldots, T_n]$ as the unital completely contractive  linear map
 $\Phi_T:C^*(S_1,\ldots, S_n)\to B(\cH)$ defined by
 \begin{equation}\label{Po-tran}
 \Phi_T(f):=\lim_{r\to 1} K_{T,r}^* (f\otimes I)K_{T,r},
 \end{equation}
 where the limit exists in the norm topology of $B(\cH)$. Moreover, we have
 $$
 \Phi_T(S_\alpha S_\beta^*)=T_\alpha T_\beta^*, \quad \alpha,\beta\in \FF_n^+.
 $$
 When $T$ is a completely noncoisometric (c.n.c.) row-contraction, i.e.,
 there is no $h\in \cH$, $h\neq 0$, such that 
 $$
 \sum_{|\alpha|=k}\|T_\alpha^* h\|^2=\|h\|^2, 
 \quad \text{\rm for any } \ k=1,2,\ldots,
 $$
an  $F_n^\infty$-functional calculus was developed  in \cite{Po-funct}.
 More precisely, if $f=\sum\limits_{\alpha\in \FF_n^+} a_\alpha S_\alpha$ is 
 in $F_n^\infty$, then
 \begin{equation}
 \label{funct-cal}
 \Phi_T(f)=f(T_1,\ldots, T_n):=
 \text{\rm SOT-}\lim_{r\to 1}\sum_{k=0}^\infty
  \sum_{|\alpha|=k} r^{|\alpha|} a_\alpha T_\alpha 
\end{equation}
exists and $\Phi_T:F_n^\infty\to B(\cH)$ is a WOT-continuous 
completely contractive homomorphism.  
More about   noncommutative Poisson transforms  on $C^*$-algebras generated
by isometries can be found in 
\cite{Po-poisson}, \cite{ArPo2}, \cite{Po-tensor}, \cite{Po-curvature}, and \cite{Po-similarity}.

I what follows we will see that
  the Poisson transform $\Phi_T$ can be extended to the noncommutative
 analytic Toeplitz algebra $F_n^\infty$ by the formula \eqref{Po-tran}, where 
 $f\in F_n^\infty$ and the limit is taken in the SO-toplogy.  Moreover, 
  $\Phi_T| F_n^\infty$   coincides
 with the $F_n^\infty$-functional calculus.  Let $\cF_n$ be the operator system defined by
$$
\cF_n:=\text{\rm span}\{fg^*:\ f, g\in F_n^\infty\}.
$$
Assume that $T:=[T_1,\ldots, T_n]$, $T_i\in B(\cH)$,  is a $C_0$-row
 contraction and let $\Phi_T:\cF_n\to B(\cH)$ be the linear map defined by
 \begin{equation}
 \label{poiss-1}
 \Phi_T(fg^*):=f(T_1,\ldots, T_n)g(T_1,\ldots, T_n)^*, \quad  f,g\in F_n^\infty,
 \end{equation}
 where  $f(T_1,\ldots, T_n)$  and $g(T_1,\ldots, T_n)$ are
  defined by  formula \eqref{funct-cal}.

 Throughout this section, $\omega:B(\cK)^{(k)}\to [0,\infty)$ is a 
 joint operator radius, where $\cK$ is a separable Hilbert space.
 For each $p\in \NN:=\{1,2,\ldots\}$, let $M_p$ denote the $p\times p$ complex
  matrices and set $\Phi_{T,p}:=  I_{M_p}\otimes \Phi_T$.
  
  We present now  von Neumann type inequalities for arbitrary
admissible or strongly admissible joint operator radii.

 \begin{theorem}\label{c0}
  Let $T:=[T_1,\ldots, T_n]$, $T_i\in B(\cH)$,  be  a $C_0$-row
 contraction and let  $\omega$ be a joint operator radius.
 \begin{enumerate}
 \item[(i)]
 If ~$\Phi_T$~ is the $F_n^\infty$-functional calculus for 
  $C_0$-row
 contractions and  $\omega$ is admissible, then
 \begin{equation}\label{ine-om1}
 \omega\left([\Phi_{T,p}(F_1)]^*, \ldots, [\Phi_{T,p}(F_k)]^*\right)\leq 
 \omega\left(F_1^*, \ldots, F_k^*\right)
 \end{equation}
 for any $F_1,\ldots, F_k\in M_p\otimes F_n^\infty$, and  $p\in \NN$.
 \item[(ii)]
 The map ~$\Phi_T$~  defined by \eqref{poiss-1} is a unital  completely contractive linear map 
which extends the $F_n^\infty$-functional 
calculus for $C_0$-row contractions, and we have
\begin{equation}\label{ext-c0}
\Phi_T(fg^*)=K_T^* (fg^*\otimes I) K_T,\quad f,g\in F_n^\infty,
\end{equation}
where $K_T$ is the Poisson kernel associated with $T$.
 Moreover, if $\omega$ is  strongly admissible, then
 \begin{equation}\label{ine-om2}
  \omega\left(\Phi_{T,p}(X_1), \ldots, \Phi_{T,p}(X_k)\right)\leq 
 \omega\left(X_1, \ldots, X_k\right)
 \end{equation}
 for any $X_1,\ldots, X_k\in M_p\otimes \cF_n $, and $p\in \NN$.
 \end{enumerate}
 \end{theorem}
 
 \begin{proof}
 Since $T$ is $C_0$-row contraction, equality \eqref{eq-ker} (case $r=1$)
  implies
 \begin{equation}\label{ker-pol}
 K_T^*(p(S_1,\ldots, S_n)\otimes I_\cH) =p(T_1,\ldots, T_n)K_T^*
 \end{equation}
 for any polynomial $p(S_1,\ldots, S_n)\in \cP$. According to \cite{Po-funct}, 
 if $f(S_1,\ldots, S_n):= \sum\limits_{k=0}^\infty 
 \sum\limits_{|\alpha|=k} a_\alpha S_\alpha$
 is
 in $F_n^\infty$, then, for any $0<r<1$, 
 $f_r(S_1,\ldots, S_n):= \sum\limits_{k=0}^\infty \sum\limits_{|\alpha|=k}
  r^{|\alpha|}a_\alpha S_\alpha$  is in the noncommutative disc algebra $\cA_n$
  and $\|f_r(S_1,\ldots, S_n)\|\leq \|f(S_1,\ldots, S_n)\|$.
  Since 
  $$
 \lim_{m\to \infty} \sum_{k=0}^m \sum_{|\alpha|=k}
  r^{|\alpha|} a_\alpha S_\alpha = f_r(S_1,\ldots, S_n)
  $$
  in the norm topology, the noncommutative von Neumann inequality 
  (see Corollary \ref{VN}) and relation 
    \eqref{ker-pol} imply
  \begin{equation*} 
 K_T^*(f_r(S_1,\ldots, S_n)\otimes I_\cH) =f_r(T_1,\ldots, T_n)K_T^*
 \end{equation*}
 for any $f(S_1,\ldots, S_n)\in F_n^\infty$.
 Using the $F_n^\infty$-functional calculus (see \eqref{funct-cal})
  and taking $r\to 1$ in 
  the above equality,
  we obtain
 \begin{equation}\label{ker-pol3}
 K_T^*(f(S_1,\ldots, S_n)\otimes I_\cH) =f(T_1,\ldots, T_n)K_T^*
 \end{equation}
 for any $f(S_1,\ldots, S_n)\in F_n^\infty$, where $f(T_1,\ldots, T_n)$ 
 is defined
 by formula \eqref{funct-cal}.
 Hence, we deduce that
 \begin{equation*}\begin{split}
 \left[\begin{matrix}
 K_T& 0&\cdots &0\\
 0& K_T&\cdots &0\\
 \vdots&\vdots&\ddots&\vdots\\
 0&0&\cdots & K_T
 \end{matrix}\right]
 &\left[f_{ij}(T_1,\ldots, T_n)\right]_{p\times p}^*\\
 &=
 \left(\left[f_{ij}(S_1,\ldots, S_n)\right]_{p\times p}^*\otimes I_\cH\right)
 \left[\begin{matrix}
 K_T& 0&\cdots &0\\
 0& K_T&\cdots &0\\
 \vdots&\vdots&\ddots&\vdots\\
 0&0&\cdots & K_T
 \end{matrix}\right]
 \end{split}
 \end{equation*}
 for any operator matrix 
 $\left[f_{ij}(S_1,\ldots, S_n)\right]_{p\times p}\in M_p\otimes F_n^\infty$.
 Let $\tilde{K}_T:\cH^{(p)}\to [F^2(H_n)\otimes \cH]^{(p)}$
 denote the $p$-fold ampliation of the Poisson kernel $K_T$.
 Since $\tilde{K}_T$ is an isometry,  the above equation  shows that
 the operator matrix
 $\left[f_{ij}(T_1,\ldots, T_n)\right]_{p\times p}^*$ is unitarily 
  equivalent to the restriction of 
  $\left[f_{ij}(S_1,\ldots, S_n)\right]_{p\times p}^*\otimes I$ to
   its invariant subspace  $\tilde{K}_T (\cH^{(p)})$.
  Taking into account that 
  $$\Phi_{T,p}\left(\left[f_{ij}(S_1,\ldots, S_n)\right]_{p\times p}\right)
  =  \left[f_{ij}(T_1,\ldots, T_n)\right]_{p\times p}
  $$
  and
     $\omega$ is an admissible operator radius,  we have
   \begin{equation*}
   \begin{split}
   \omega\left([\Phi_{T,p}(F_1)]^*,\ldots, [\Phi_{T,p}(F_k)]^*\right)&=
   \omega\left((F_1^*\otimes I_\cH)|\tilde{K}_T (\cH^{(p)}),\ldots, 
    (F_k^*\otimes I_\cH)|\tilde{K}_T (\cH^{(p)})\right)\\
    &\leq\omega(F_1^*\otimes I_\cH, \ldots, F_k^*\otimes I_\cH)\\
    &=\omega(F_1^*, \ldots, F_k^*)
   \end{split}
   \end{equation*}
   for any
   $F_1,\ldots, F_k\in M_p\otimes F_n^\infty$, $p\in \NN$.
  This completes the proof of part (i) of the theorem.
 
 To prove part (ii), note that  relation \eqref{ker-pol3} implies
 \begin{equation}
 \label{fg}
 K_T^*[f(S_1,\ldots, S_n)g(S_1,\ldots, S_n)^*\otimes I_\cH]K_T= 
 f(T_1,\ldots, T_n)g(T_1,\ldots, T_n)^*
 \end{equation}
 for any $f(S_1,\ldots, S_n)$, $g(S_1,\ldots, S_n)\in F_n^\infty$, 
 which proves
 \eqref{ext-c0}.
 Since $K_T$ is an isometry, it is clear  that $\Phi_T$
 is a unital completely contractive linear map.
 Moreover, using relation \eqref{fg}, we get
 $$\tilde{K}_T^* (X_j\otimes I_\cH)\tilde{K}_T=\Phi_{T,p} (X_j)
 $$
 for any  $X_j \in M_p\otimes \cF_n $, $j=1,\ldots, k$,  and $p\in \NN$. 
 Hence, and taking into account that $\omega$ is strongly admissible,
 we  deduce that
 \begin{equation*}
   \begin{split}
   \omega\left(\Phi_{T,p}(X_1),\ldots, \Phi_{T,p}(X_k)\right)&=
   \omega\left(\tilde{K}_T^*(X_1\otimes I_\cH)\tilde{K}_T,\ldots, 
    \tilde{K}_T^*(X_k\otimes I_\cH)\tilde{K}_T \right)\\
    &\leq\omega(X_1\otimes I_\cH, \ldots, X_k\otimes I_\cH)\\
    &=\omega(X_1, \ldots, X_k).
   \end{split}
   \end{equation*}
  This completes the proof of the theorem.
 \end{proof}

We remark that 
 the map $\Phi_T$ defined by \eqref{poiss-1} can be extended  to a completely contractive linear map from $B(F^2(H_n))$
 to $B(\cH)$ by setting 
$$
 \Phi_T(X)=K_T^* (X\otimes I) K_T,\quad  X\in B(F^2(H_n)).
 $$

\smallskip

 If $[T_1,\ldots, T_n]$ is a row contraction and 
 $f(S_1,\ldots, S_n):=\sum\limits_{k=0}^\infty \sum\limits_{|\alpha|=k} a_\alpha S_\alpha$ 
 is in $F_n^\infty$, we denote
 $$
 f_r(T_1,\ldots, T_n):=\sum\limits_{k=0}^\infty \sum\limits_{|\alpha|=k} 
 a_\alpha
  r^{|\alpha|}T_\alpha, \quad
  0<r<1,
  $$
   where the convergence is in 
  the  operator norm topology. 
 
 \begin{lemma}\label{le-fo}
 If $[T_1,\ldots, T_n]$, $T_i\in B(\cH)$,  is a row contraction, then
 \begin{equation}\label{KTr}
 K_{T,r}^* (f(S_1,\ldots, S_n)\otimes I_\cH)=f_r(T_1,\ldots, T_n) K_{T,r}^*
 \end{equation}
 for any $f(S_1,\ldots, S_n)\in F_n^\infty$ and $0<r<1$.
 Moreover, if $T$ is c.n.c., then
 \begin{equation}\label{sot-lim}
 f(T_1,\ldots, T_n)=\text{\rm SOT-}\lim_{r\to 1} K_{T,r}^* 
 (f(S_1,\ldots, S_n)\otimes I_\cH)K_{T,r}
 \end{equation}
 for any $f(S_1,\ldots, S_n)\in F_n^\infty$.
 \end{lemma}
 \begin{proof}
 The  relation \eqref{eq-ker} implies
 \begin{equation}\label{eq-ker2}
 K_{T,r}^*[p(S_1,\ldots, S_n)\otimes I_\cH]=p(rT_1,\ldots, rT_n) K_{T,r}^*
 \end{equation}
 for any polynomial $p(S_1,\ldots, S_n)$ in $\cP$.
 If $0<t<1$, then $f_t(S_1,\ldots, S_n)$ is in noncommutative disc algebra
  $\cA_n$ and 
 $$
 \lim_{m\to \infty}\sum_{|\alpha|\leq m}t^{|\alpha|} a_\alpha S_\alpha
  =f_t(S_1,\ldots, S_n),
 $$
 where the convergence is in the operator norm.
 Since $[rT_1,\ldots, rT_n]$ is a $C_0$-row contraction, the noncommutative
  von Neuman  inequality  implies
  $$
 \lim_{m\to \infty} \sum_{|\alpha|\leq m}
 t^{|\alpha|} r^{|\alpha|} a_\alpha T_\alpha
   =f_t(rT_1,\ldots, rT_n)
 $$
 in the operator norm.
 Using now relation \eqref{eq-ker2} where 
 $p(S_1,\ldots, S_n):=
 \sum\limits_{|\alpha|\leq m}t^{|\alpha|} a_\alpha S_\alpha$ 
 and taking limit as 
 $m\to \infty$, we get
 \begin{equation}
 \label{eq-ker3}
 K_{T,r}^* (f_t(S_1,\ldots, S_n)\otimes I_\cH)=f_t(rT_1,\ldots, rT_n) 
 K_{T,r}^*.
 \end{equation}
  Let us prove  that 
 \begin{equation}\label{lim-t}
 \lim_{t\to 1} f_t(rT_1,\ldots, rT_n)=f_r(T_1,\ldots, T_n),
 \end{equation}
 where the convergence is in the operator norm. 
  First, notice that, if  $\epsilon>0$,  there is   $m_0\in \NN$
 such that $\left(\sum\limits_{k=m_0}^\infty r^k\right) \|f\|_2<\frac {\epsilon} {2}$,
 where $\|f\|_2:=\|f(S_1,\ldots, S_n)(1)\|$.
 Since $[T_1,\ldots, T_n]$ is a row contraction, the latter inequality implies 
 \begin{equation*}
 \begin{split}
 \sum_{k=m_0}^\infty r^k \left\|\sum_{|\alpha|=k} a_\alpha T_\alpha\right\|
 &\leq 
 \sum_{k=m_0}^\infty r^k \left(\sum_{|\alpha|=k} |a_\alpha|^2\right)^{1/2}  \\
 &\leq 
 \left(\sum_{k=m_0}^\infty r^k\right) \|f\|_2<\frac {\epsilon} {4}.
 \end{split}
 \end{equation*}
 Consequently, we have
 \begin{equation*}
 \begin{split}
 \left\|\sum_{k=0}^\infty \sum_{|\alpha|=k}
  t^{|\alpha|} r^{|\alpha|} a_\alpha T_\alpha -
  \sum_{k=0}^\infty \sum_{|\alpha|=k}
    r^{|\alpha|} a_\alpha T_\alpha\right\|
    &\leq 
    \left\| \sum_{k=1}^{m_0-1} \sum_{|\alpha|
    =k}r^k(t^k-1)a_\alpha T_\alpha\right\|+ \frac{\epsilon}{2}\\
    &\leq \sum_{k=1}^{m_0-1} r^k(t^k-1)\|f\|_2+\frac{\epsilon}{2}.
 \end{split}
 \end{equation*}
 Now, it clear that there exists $0<\delta<1$ such that 
 $\sum\limits_{k=1}^{m_0-1} r^k(t^k-1)\|f\|_2<\frac{\epsilon}{2}$ for any $t\in (\delta, 1)$.
  This proves relation \eqref{lim-t}.
Since SOT-$\lim\limits_{t\to 1} f_t(S_1,\ldots, S_n)=f(S_1,\ldots, S_n)$, 
 we can 
 use relation \eqref{lim-t} and then take the limit in 
 \eqref{eq-ker3}, as $t\to 1$, to obtain
 relation \eqref{KTr}. 
 The last part of the lemma follows from  relations \eqref{KTr} 
  and  
 \eqref{funct-cal}.
 The proof is complete.
 \end{proof}

 Under the additional condition of norm quasi-continuity of the joint 
 operator radius,
 we can prove the following result.

\begin{theorem}\label{von-ra}
Let $T:=[T_1,\ldots, T_n]$, $T_i\in B(\cH)$,  be  a  row
 contraction,
$\Phi_T$ be the Poisson transform associated with $T$,
  and let  $\omega$ be a joint operator radius.
 \begin{enumerate}
 \item[(i)]
 If $\omega$ is norm quasi-continuous and  admissible, then
 \begin{equation}\label{ine-om3}
  \omega\left([\Phi_{T,p}(F_1)]^*, \ldots, [\Phi_{T,p}(F_k)]^*\right)\leq 
 \omega\left(F_1^*, \ldots, F_k^*\right)
 \end{equation}
 for any $F_1,\ldots, F_k\in M_p\otimes \cA_n $,  and $p\in \NN$.
 \item[(ii)]
 If $\omega$ is  norm quasi-continuous and strongly admissible, then
 \begin{equation}\label{ine-om4}
  \omega\left(\Phi_{T,p}(X_1), \ldots, \Phi_{T,p}(X_k)\right)\leq 
 \omega\left(X_1, \ldots, X_k\right)
 \end{equation}
 for any $X_1,\ldots, X_k\in M_p\otimes C^*(S_1,\ldots, S_n) $, and  $p\in \NN$.
 \end{enumerate}
 \end{theorem}
 \begin{proof}
 According to Lemma \ref{le-fo}, we have
 $$(F(S_1,\ldots, S_n)^*\otimes I)\tilde{K}_{T,r}= \tilde{K}_{T,r} 
 F(rT_1,\ldots, rT_n)^*
 $$
 for any $F=F(S_1,\ldots, S_n):=[f_{ij}(S_1,\ldots, S_n)]_{p\times p}\in M_p\otimes \cA_n$,
 where $\tilde{K}_{T,r}$ is the ampliation of the Poisson kernel $K_{T,r}$.
 Since  $K_{T,r}$ is an isometry and $\omega$ is admissible, we have
 \begin{equation*}
 \begin{split}
 \omega\left(F_1(rT_1,\ldots, rT_n)^*\right.,\ldots, &\left.F_k(rT_1,\ldots, rT_n)^*\right)\\
 &=
 \omega\left((F_1^*\otimes I)|\tilde{K}_{T,r}(\cH^{(p)}), \ldots,  
 (F_k^*\otimes I)|\tilde{K}_{T,r}(\cH^{(p)})\right)\\
 &\leq\omega(F_1^*\otimes I, \ldots, F_k^*\otimes I)\\
 &=\omega(F_1^*, \ldots, F_k^*)
 \end{split}
 \end{equation*}
 for any $F_1,\ldots, F_k\in M_p\otimes \cA_n $,  and $p\in \NN$.
 According  to relation \eqref{Po-tran}, 
 $F_j(rT_1,\ldots, rT_n)\to \Phi_{T,p}(F_j)$
 in the operator norm as $r\to 1$.
 Consequently, since $\omega$ is norm quasi-continuous 
  the above inequality implies
  \eqref{ine-om3}.  This completes the proof of  
 part (i) of the theorem.
 
 Let us  prove part (ii).
 Since  $\omega$ is strongly admissible and $K_{T,r}$ is an isometry, 
 we have
 \begin{equation*}
 \begin{split}
 \omega\bigl(\tilde{K}_{T,r}^* (X_1\otimes I)\tilde{K}_{T,r},\ldots,  
 \tilde{K}_{T,r}^* (X_k\otimes I)\tilde{K}_{T,r})&\leq
 \omega(X_1\otimes I, \ldots, X_k\otimes I\bigr)\\
 &=\omega(X_1, \ldots, X_k)
 \end{split}
 \end{equation*}
 for any  $X_1,\ldots, X_k\in M_p\otimes C^*(S_1,\ldots, S_n)$.
  Now, using relation \eqref{Po-tran}, the norm continuity of $\omega$, 
 and taking 
 $r\to 1$, we deduce inequality \eqref{ine-om4}. The proof is complete.
 \end{proof}

  Consider the linear span
$$
\cD_n:=\text{\rm span}\left\{ f(S_1,\ldots, S_n)S_\alpha^*:
\ f(S_1,\ldots, S_n)\in F_n^\infty, \alpha\in \FF_n^+
\right\}.
$$
We remark that the norm closed non-selfadjoint 
algebra generated by the noncommutative analytic Toeplitz algebra $F_n^\infty$ and the 
$C^*$-algebra $C^*(S_1,\ldots, S_n)$ coincides with  the norm closure of
 $ \cD_n$. Indeed,
since $C^*(S_1,\ldots, S_n)$ is the norm closure of
 polynomials of the form  ~$\sum a_{\alpha\beta} S_\alpha S_\beta^*$, it is clear that
 $f(S_1,\ldots, S_n) Y\in \overline{\cD}_n$ for any
  $f(S_1,\ldots, S_n)\in F_n^\infty$
 and $Y\in C^*(S_1,\ldots, S_n)$.  To show that
  $Y f(S_1,\ldots, S_n) \in \overline{\cD}_n$, it is enough to prove that 
  $$S_j^* f(S_1,\ldots, S_n)\in \cP^*+ F_n^\infty.
  $$
   Since
  \begin{equation*}
  f(S_1,\ldots, S_n)=\sum_{i=1}^n S_i\varphi_i(S_1,\ldots, S_n)+ a_0 I
  \end{equation*}
  for any $j=1, \ldots, n$, 
for some elements $\varphi_i(S_1,\ldots, S_n)\in F_n^\infty$ and $a_0\in \CC$, and 
taking into account that $S_iS_j^*=\delta_{ij} I$, $i,j=1,\ldots, n$,  we have
$$
S_j^* f(S_1,\ldots, S_n)=a_0 S_j^*+ \varphi_j(S_1,\ldots, S_n),
$$
which proves our assertion.
We can view ${\cD}_n$ as a noncommutative Douglas
 type subalgebra (see \cite{Do-book}) of $B(F^2(H_n))$.
Let $T:=[T_1,\ldots, T_n]$, $T_i\in B(\cH)$,  be a c.n.c. row contraction
and
 define the linear map $\Phi_T:\cD_n+\cD_n^*\to B(\cH)$
by setting
\begin{equation}\label{def-doug}
\Phi_T\bigl(f(S_1,\ldots, S_n) S_\alpha^*+ S_\beta g(S_1,\ldots, S_n)^*\bigr)
:=f(T_1,\ldots, T_n) T_\alpha^*+ T_\beta
g(T_1,\ldots, T_n)^*
\end{equation}
for any $f(S_1,\ldots, S_n), g(S_1,\ldots, S_n)\in F_n^\infty$, 
and $\alpha,\beta\in \FF_n^+$.

 \smallskip
 
 \begin{theorem}\label{von-ra-cnc}
  Let $T$ be a c.n.c. row contraction and let
    $\omega$ be a joint operator radius. 
  The map $\Phi_T$ defined by \eqref{def-doug}   can be extended 
   to a  unital completely contractive linear map
from the norm closure of $\cD_n+\cD_n^*$  to $B(\cH)$ by setting
\begin{equation}\label{def-ext}
\Phi_T(X):=\text{\rm WOT-}\lim_{r\to 1}K_{T,r}^* [X\otimes I)]
 K_{T,r}, \quad X\in \overline{\cD_n+\cD_n^*},
 \end{equation}
where   $\{K_{T,r}\}_{0<r<1}$  is the Poisson
kernel associated with $T$.
  \begin{enumerate}
  \item[(i)] If $\omega$ is $*$-SOT quasi-continuous and  admissible, then
  \begin{equation}\label{ine-omi}
 \omega\left([\Phi_{T,p}(F_1)]^*, \ldots, [\Phi_{T,p}(F_k)]^* \right)\leq 
 \omega\left(F_1^*, \ldots, F_k^* \right)
 \end{equation}
 for any $F_1,\ldots, F_k\in M_p\otimes F_n^\infty$, and $p\in \NN$.
   \item[(ii)]
   If $\omega$ is SOT quasi-continuous and strongly 
 admissible, then 
 \begin{equation}\label{ine-om5}
 \omega\left(\Phi_{T,p}(Y_1), \ldots, \Phi_{T,p}(Y_k) \right)\leq 
 \omega\left(Y_1, \ldots, Y_k \right)
 \end{equation}
 for any $Y_1,\ldots, Y_k\in M_p\otimes  \overline{\cD}_n $, and $p\in \NN$.
 \item[(iii)]
  If  $\omega$ is WOT quasi-continuous and strongly 
 admissible, then  inequality \eqref{ine-om5}
 holds for any $Y_1,\ldots, Y_k\in M_p\otimes  \overline{\cD_n+\cD_n^*}$, and 
 $p\in \NN$.
 \end{enumerate}
\end{theorem}

 \begin{proof}
 According to Lemma \ref{le-fo}, we have
 \begin{equation}\label{K*K}
 K_{T,r}^* (f(S_1,\ldots, S_n)p(S_1,\ldots, S_n)^* \otimes I_\cH)K_{T,r}=
 f_r(T_1,\ldots, T_n) p_r(T_1,\ldots, T_n)^*  
 \end{equation}
for any  $f(S_1,\ldots, S_n)\in F_n^\infty$ and any polynomial $p\in \cP$.
Since
$$
\|f_r(T_1,\ldots, T_n)\|\leq \|f(S_1,\ldots, S_n)\|\quad \text{ for any } 
\ 0<r<1,
$$
 and 
$p_r(T_1,\ldots, T_n)\to p(T_1,\ldots, T_n)$ in norm as $r\to 1$, 
one can use relations
\eqref{funct-cal} and \eqref{K*K} to deduce that
\begin{equation}\label{SOT1}
\text{\rm SOT-}\lim_{r\to 1} 
K_{T_r} \bigl[f(S_1,\ldots, S_n)p(S_1,\ldots, S_n)^* \otimes I_\cH\bigr]K_{T,r}=
 f(T_1,\ldots, T_n) p(T_1,\ldots, T_n)^*  
\end{equation}
 Hence, we  obtain
 \begin{equation}\label{SOT2}
 \begin{split}
\text{\rm WOT-}\lim_{r\to 1} 
K_{T_r} &\bigl[\bigl(q(S_1,\ldots, S_n)g(S_1,\ldots, S_n)^*+
f(S_1,\ldots, S_n)p(S_1,\ldots, S_n)^*\bigr) \otimes I_\cH\bigr]K_{T,r}\\
&
 =q(T_1,\ldots, T_n)g(T_1,\ldots, T_n)^*+
 f(T_1,\ldots, T_n) p(T_1,\ldots, T_n)^*  
 \end{split}
\end{equation}
for any $f,g\in F_n^\infty$ and $p,q\in \cP$. 
This shows that the linear map $\Phi_T:\cD_n+ \cD_n^*\to B(\cH)$ defined by
\eqref{def-doug} is  completely contractive and  the limit in \eqref{def-ext} exists
for any 
 $X\in \cD_n+ \cD_n^*$.
 
  Now, we show that the limit  in \eqref{def-ext} 
 exists  for any 
 $X\in \overline{\cD_n+ \cD_n^*}$.
 Let $\{X_m\}_{m=1}^\infty$ be a sequence of operators in $\cD_n+ \cD_n^*$
  such that $X_m\to X$ in norm as $m\to \infty$. Since $\Phi|(\cD_n+ \cD_n^*)$
  is a completely contractive linear map, we have
  $$
  \|\Phi_T(X_p)-\Phi(X_m)\|\leq \|X_p-X_m\|, \quad p,m\in\NN,
  $$
  which implies that $\{\Phi_T(X_m)\}_{m=1}^\infty$ is a Cauchy sequence
  in  the operator norm.
  Therefore, there exists $Y\in B(\cH)$ such that $\Phi_T(X_m)\to Y$ in norm as 
  $m\to\infty$.
  Note that, for any  $h,k\in \cH$, we have
  \begin{equation*}
  \begin{split}
  \bigl|\left<K_{T,r}^*X K_{T,r}h,k\right>&-\left< Yh,k\right>\bigr|\\
  &\leq 
  \left|\left<K_{T,r}^*(X-X_m) K_{T,r}h,k\right>\right|+
  \left|\left<K_{T,r}^*X_m K_{T,r}h,k\right>
  -\left< \Phi_T(X_m)h,k\right>\right|\\
  &+\left|\left< \Phi_T(X_m)h,k\right>-\left< Yh,k\right>\right|.
  \end{split}
  \end{equation*}
 Assume that $h\neq 0$, $k\neq 0$, and choose $m_0$ 
 such that 
  $$
  \|X_m-X\|\leq \frac {\epsilon}{3\|h\|\|k\|} \ \text{\ and \ } \ 
  \|\Phi(X_m)-Y\|\leq \frac {\epsilon}{3\|h\|\|k\|}
  $$
  for any  $m\geq m_0$. 
 According to relations \eqref{SOT2} and \eqref{def-doug},
  we can find $0<\delta<1$ such that 
 $$
 \left|\left<K_{T,r}^*X_m K_{T,r}h,k\right>
  -\left< \Phi_T(X_m)h,k\right>\right|\leq \frac{\epsilon}{3}
  $$
  for any $r\in (\delta, 1)$. 
   Now, the above inequalities imply
  $$
 \bigl|\left<K_{T,r}^*X K_{T,r}h,k\right>-\left< Yh,k\right>\bigr|\leq \epsilon
 $$
 for any $r\in (\delta, 1)$. 
 This proves that 
 $$
 \lim_{r\to 1}\left<K_{T,r}^*X K_{T,r}h,k\right>=\left<Yh,k\right>,
 \quad h,k\in \cH.
 $$
 Therefore, the limit in \eqref{def-ext} exists for any
 $X\in \overline{\cD_n+ \cD_n^*}$
 and $\Phi_T$ is a completely contractive linear map.
 Similarly, using equation \eqref{SOT1} instead of \eqref{SOT2}, 
 one can prove that
 \begin{equation}
 \label{sot2}
 \Phi_T(X)=\text{\rm SOT-}\lim_{r\to 1}K_{T,r}^* [X\otimes I)]
 K_{T,r}, \quad \text{ for any }\ X\in \overline{\cD}_n.
 \end{equation}
 
Now,  to prove (i), notice that Lemma \ref{le-fo} implies
  $$
  \tilde{K}_{T,r} [\Phi_{T,p}((F)_r)]^*=(F^*\otimes I_\cH)\tilde{K}_{T,r}
  $$
  for any $F=[f_{ij}]_{p\times p}\in M_p\otimes F_n^\infty$, 
  where $F_r:=[(f_{ij})_r]_{p\times p}$, $0<r<1$. 
 Since $\omega$ is admissible, the latter equality implies
 \begin{equation}
 \label{ome}\begin{split}
 \omega\bigl( [\Phi_{T, p}\left( (F_1)_r\right)]^*,\ldots,
  [\Phi_{T, p}\left((F_k)_r\right)]^*\bigr)
 &\leq \omega(F_1^*\otimes I_\cH, \ldots, F_k^*\otimes I_\cH)\\
 &= \omega(F_1^*, \ldots, F_k^*)
 \end{split}
 \end{equation}
 for any $F_1,\ldots, F_k\in M_p\otimes F_n^\infty$ and $p\in \NN$.
 According to the $F_n^\infty$-functional calculus for c.n.c. row contractions,  
 $[\Phi_{T, p}\left( (F_j)_r\right)]$ is SOT convergent to 
 $[\Phi_{T, p}\left(F_j \right)]$
   as $r\to 1$. Since $\omega$ is $*$-SOT quasi-continuous, inequality
 \eqref{ome} implies \eqref{ine-omi}.

 Let us prove part (ii) of the theorem.
 Since  $\omega$ is strongly admissible and $K_{T,r}$ is an isometry, 
 we have
 \begin{equation*}
 \begin{split}
 \omega\bigl(\tilde{K}_{T,r}^* (Y_1\otimes I)\tilde{K}_{T,r},\ldots,  
 \tilde{K}_{T,r}^* (Y_k\otimes I)\tilde{K}_{T,r}\bigr)&\leq
 \omega(Y_1\otimes I, \ldots, Y_k\otimes I)\\
 &=\omega(Y_1, \ldots, Y_k)
 \end{split}
 \end{equation*}
 for any  $Y_1,\ldots, Y_k\in M_p\otimes\overline{\cD}_n$ and $p\in \NN$.
 Consequently, using relation \eqref{sot2} and the SOT-continuity of $\omega$,
  we deduce 
 inequality  \eqref{ine-om5}. 
 The proof of part (iii) is similar to that of part (ii) but uses 
 relation \eqref{def-ext} instead and the WOT-continuity of $\omega$.
 The proof is complete.
 \end{proof}

 \begin{corollary}\label{cor-sot}
 If $X\in \overline{\cD}_n$, then 
$$
\Phi_T(X)=\text{\rm SOT-}\lim_{r\to 1}K_{T,r}^* [X\otimes I)]
 K_{T,r}.
 $$
 Moreover, the map defined by \eqref{def-ext} extends the
  $F_n^\infty$-functional calculus for c.n.c.
 row contractions and the Poisson transform defined by \eqref{Po-tran}.
\end{corollary}

\smallskip

  \section{Constrained von Neumann inequalities}
  \label{von2}

       If a row contraction  $[T_1,\ldots, T_n]$ satisfies certain   constrains,
        we prove that there is a suitable invariant subspace 
        $\cE\subset F^2(H_n)$ under each operator $S_1^*, \ldots, S_n^*$, such
         that an inequality of the form
        $$\omega(f_1(T_1,\ldots, T_n),\ldots, f_k(T_1,\ldots, T_n))\leq
      \omega(P_\cE f_1(S_1,\ldots, S_n)|\cE,\ldots, P_\cE f_k(S_1,\ldots, S_n)|
      \cE)
      $$
      holds, where $\omega$ is an admissible or strongly admissible 
      joint operator radii.
      
     This type  of constrained von Neumann inequalities as well as matrix-valued
     versions are 
     considered in  this section.
     In particular, we obtain appropriate  multivariable generalizations
     of several von Neumann type inequalities obtained in
     \cite{vN}, \cite{PY}, \cite{Sz2}, \cite{Dr}, \cite{Po-von}, \cite{Po-funct}, 
     \cite{Po-disc}, \cite{Po-poisson},  \cite{Arv1},   \cite{ArPo2},
     and \cite{Po-tensor},
     to joint operator
     radii. 
     For example, 
     if $[T_1,\ldots, T_n]$  is a row contraction with  commuting entries, then we prove that
     $$\omega(f_1(T_1,\ldots, T_n),\ldots, f_k(T_1,\ldots, T_n))\leq
      \omega(f_1(B_1,\ldots, B_n),\ldots, f_k(B_1,\ldots, B_n))
      $$
        if $f_1(B_1,\ldots, B_n),\ldots, f_k(B_1,\ldots, B_n)$  are elements of 
        operator algebras such as
      the commutative disc algebra $\cA_n^c$,  the Toeplitz $C^*$-algebra
       $C^*(B_1,\ldots, B_n)$,   Arveson's algebra
        $W_n^\infty$, 
        the commutative Douglas type algebra $\overline{\cD_n^c}$, and  the operator system
         \ $W_n^\infty (W_n^\infty)^*$. 
        The operator $f_j(T_1,\ldots, T_n)$ 
       is  defined by an appropriate functional calculus for row contractions 
       with commuting entries.
 
Let $\cH$ be a separable complex Hilbert space, 
 $J$ be a WOT-closed  two-sided ideal of $F^\infty_n$, and let 
$T:=[T_1,\ldots, T_n]$, $T_i\in B(\cH)$,  be  a $C_0$-row
 contraction  such that 
 $$
 \varphi(T_1,\ldots, T_n)=0,\quad  \text{ for any } \ \varphi (S_1,\ldots, S_n)\in J.
 $$
Denote by $W_J^\infty$   the WOT-closed algebra generated by the operators
 $B_i:=P_{\cN_J} S_i |\cN_J$, \ $i=1,\ldots, n$, and the identity, where
 $$
 \cN_J: F^2(H_n)\ominus \cM_J\quad \text{ and }\quad 
 \cM_J:=\overline{ J F^2(H_n)}.
 $$
 Notice that
 $\cM_J:=\overline{\{\varphi:\ \varphi(S_1,\ldots, S_n)\in J \}}^{\|\cdot\|_2}$
 where $\varphi:=\varphi(S_1,\ldots, S_n)(1)$, and
 $$
 \cN_J=\bigcap_{\varphi\in J}\ker \varphi(S_1,\ldots, S_n)^*.
 $$
 We proved in \cite{ArPo2} that 
 \begin{equation}\label{w-formulas}
 W_J^\infty=P_{\cN_J} F_n^\infty |\cN_J=\{f(B_1,\ldots, B_n):
 \ f(S_1,\ldots, S_n)\in F_n^\infty\},
 \end{equation}
where $f(B_1,\ldots, B_n)$ is defined by the
 $F_n^\infty$-functional calculus  for $C_0$-row contractions.
Moreover, the map $\Phi_T^J:W_J^\infty\to B(\cH)$
defined by
\begin{equation}\label{J-cal}
\Phi_T^J(f(B_1,\ldots, B_n)):= f(T_1,\ldots, T_n)
\end{equation}
is a WOT-continuous, completely contractive homomorphism (see Section 4 of \cite{ArPo2}),
 which we call the $W_J^\infty$-functional calculus.
Now, consider the operator system
 defined by
$$
\cW_n^J:=\text{\rm span}\{f(B_1,\ldots, B_n)g(B_1,\ldots, B_n)^*:
\ f(B_1,\ldots, B_n), g(B_1,\ldots, B_n)\in W_J^\infty\}.
$$
Let $\Phi_T^J:\cW_n^J\to B(\cH)$ be the linear map defined by
 \begin{equation}
 \label{J-poiss-1}
 \Phi_T^J[f(B_1,\ldots, B_n)g(B_1,\ldots, B_n)^*]:=
 f(T_1,\ldots, T_n)g(T_1,\ldots, T_n)^*,  
 \end{equation}
 where  $f(T_1,\ldots, T_n)$  and $g(T_1,\ldots, T_n)$ are defined by 
 formula
 \eqref{funct-cal}.

  In what follows, we obtain constrained von Neumann inequalities for
   arbitrary admissible
(resp. strongly admissible) joint operator radii. 
 
 \begin{theorem}\label{J-c0}
  Let 
    $\omega$ be a joint operator radius, $J$ be a WOT-closed two-sided ideal
     of $F_n^\infty$, and let 
$T:=[T_1,\ldots, T_n]$, $T_i\in B(\cH)$,  be  a $C_0$-row
 contraction  such that 
 $$
 \varphi(T_1,\ldots, T_n)=0,\quad \varphi(S_1,\ldots, S_n)\in J.
 $$
 Then, 
 the map ~$\Phi_T^J$~  defined by \eqref{J-poiss-1} is a unital  completely contractive 
 linear map 
  and  
\begin{equation}\label{J-ext-c0}
\Phi_T^J\bigl(f(B_1,\ldots, B_n)g(B_1,\ldots, B_n)^*\bigr)
=K_T^* [f(S_1,\ldots, S_n)g(S_1,\ldots, S_n)^*\otimes I] K_T
\end{equation}
for any $f(S_1,\ldots, S_n), g(S_1,\ldots, S_n)\in F_n^\infty$,
where $K_T$ is the Poisson kernel associated with $T$.
 \begin{enumerate}
 \item[(i)]
 If   $\omega$ is admissible, then
 \begin{equation}\label{J-ine-om1}
 \omega\left([\Phi^J_{T,p}(F_1)]^*, \ldots, [\Phi^J_{T,p}(F_k)]^*\right)\leq 
 \omega\left(F_1^*, \ldots, F_k^*\right)
 \end{equation}
 for any $F_1,\ldots, F_k\in M_p\otimes W_J^\infty$, and  $p\in \NN$.
 \item[(ii)]
 If $\omega$ is  strongly admissible, then
 \begin{equation}\label{J-ine-om2}
  \omega\left(\Phi^J_{T,p}(X_1), \ldots, \Phi^J_{T,p}(X_k)\right)\leq 
 \omega\left(X_1, \ldots, X_k\right)
 \end{equation}
 for any $X_1,\ldots, X_k\in M_p\otimes \cW^J_n $, and $p\in \NN$.
 \end{enumerate}
 \end{theorem}
 
 \begin{proof}
 According to Theorem \ref{c0} (see relation \eqref{ker-pol3}), we have
 \begin{equation}
 \label{ker-p}
 \left<(\varphi(S_1,\ldots, S_n)^*\otimes I)K_Th,1\otimes k\right>=
 \left<K_T\varphi(T_1,\ldots, T_n)^*h,1\otimes k\right>
 \end{equation}
 for any $\varphi(S_1,\ldots, S_n)\in F_n^\infty$ and $h,k\in \cH$.
 Note that if $\varphi(S_1,\ldots, S_n)\in J$, then 
 $\varphi(T_1,\ldots, T_n)=0$, and relation \eqref{ker-p} implies
 $\left<K_Th, \varphi\otimes k\right>=0$ for any $h,k\in \cH$.
 Taking into account the definition of $\cM_J$, we deduce that 
 $K_T(\cH)\subseteq \cN_J\otimes \cH$. This implies 
 \begin{equation}
 \label{KT}
 K_T=\left(P_{ \cN_J}\otimes I_\cH\right) K_T.
 \end{equation}
 
 Since $J$ is a left ideal of $F_n^\infty$, $\cN_J$ is an invariant subspace 
 under each operator $S_1^*,\ldots, S_n^*$ and therefore
  $B_\alpha=P_{ \cN_J}S_\alpha|\cN_J$,
 $\alpha\in \FF_n^+$.
 Since $[B_1,\ldots, B_n]$ is a $C_0$-row contraction, we can use 
 the $F_n^\infty$-functional calculus to deduce that
 \begin{equation}
 \label{fnn}
 f(B_1,\ldots, B_n)=P_{ \cN_J}f(S_1,\ldots, S_n)|\cN_J
 \end{equation}
 for any $f(S_1,\ldots, S_n)\in F_n^\infty$.
 Taking into account relations \eqref{ker-pol3} (see Theorem \ref{c0}), 
  \eqref{KT}, and \eqref{fnn},
 we obtain
 \begin{equation*}
 \begin{split}
 K_Tf(T_1,\ldots, T_n)^*&=\left(P_{ \cN_J}\otimes I_\cH\right)
 [f(S_1,\ldots, S_n)^*\otimes I_\cH]\left(P_{ \cN_J}\otimes I_\cH\right)K_T\\
 &=
 \left[\left(P_{ \cN_J}f(S_1,\ldots, S_n)| 
 \cN_J\right)^*\otimes I_\cH\right] K_T\\
&= 
 \left[f(B_1,\ldots, B_n)^*\otimes I_\cH\right]K_T.
 \end{split}
 \end{equation*}
Therefore, we have
\begin{equation}\label{KFFK}
 K_Tf(T_1,\ldots, T_n)^*=\left[f(B_1,\ldots, B_n)^*\otimes I_\cH\right]K_T
 \end{equation}
 for any $f(B_1,\ldots, B_n)\in W_J^\infty$.
 Hence, we obtain
 \begin{equation}
 \label{fg2}
 f(T_1,\ldots, T_n)g(T_1,\ldots, T_n)^*=K_T^*\left[
 f(B_1,\ldots, B_n)g(B_1,\ldots, B_n)^*\otimes I_\cH
 \right]K_T,
 \end{equation}
 which proves relation \eqref{J-ext-c0} and that $\Phi_T^J$ is a unital completely contractive linear map.
 
 To prove part (i) of the theorem, notice that relation
 \eqref{KFFK} implies
 \begin{equation}
 \label{K-matr}
 \tilde{K}_T\left[f_{ij}(T_1,\ldots, T_n)\right]^*_{p\times p}=
 \left(\left[f_{ij}(B_1,\ldots, B_n)\right]^*_{p\times p}\otimes I_\cH\right)
 \tilde{K}_T,
 \end{equation}
 where $\tilde{K}_T$ is the $p$-fold ampliation of $K_T$.
  Since $\tilde{K}_T$ is an isometry, $\omega$ is admissible, and 
  $$
  \Phi_{T,p}^J\left([f_{ij}(B_1,\ldots, B_n)]_{p\times p}\right)=  
  [f_{ij}(T_1,\ldots, T_n)]_{p\times p},
  $$
  one can use relation \eqref{K-matr} to deduce that
  \begin{equation*}
  \begin{split}
  \omega\left(\left[\Phi_{T,p}(F_1)\right]^*,\ldots,
   \left[\Phi_{T,p}(F_k)\right]^*\right)
   &=\omega\left((F_1^*\otimes I_\cH)|\tilde{K}_T(\cH^{(p)}),\ldots,
   (F_k^*\otimes I_\cH)|\tilde{K}_T(\cH^{(p)})\right)\\
   &\leq \omega\left(F_1^*\otimes I_\cH,\ldots, F_k^*\otimes I_\cH\right)\\
   &=\omega(F_1^*,\ldots, F_k^*)
  \end{split}
  \end{equation*} 
  for any $F_1,\ldots, F_k\in M_p\otimes W_J^\infty$, and $p\in \NN$.
 
 To prove the last part of the theorem, notice that relation \eqref{fg2} implies
 $$\Phi_{T,p}(X)=\tilde{K}_T^*(X\otimes I_\cH)\tilde{K}_T
 $$
 for any $X\in M_p\otimes \cW_n^J$ and $p\in \NN$. 
 Since $\omega$ is strongly admissible and $\tilde{K}_T$ is an isometry,
 we  have 
 \begin{equation*} \begin{split}
  \omega\left(\Phi^J_{T,p}(X_1), \ldots, \Phi^J_{T,p}(X_k)\right)
  &=\omega\left(\tilde{K}_T^*(X_1\otimes I_\cH)\tilde{K}_T, \ldots,
  \tilde{K}_T^*(X_k\otimes I_\cH)\tilde{K}_T\right)\\
  &\leq 
 \omega\left(X_1\otimes I_\cH, \ldots, X_k\otimes I_\cH\right)\\
 &=\omega\left(X_1, \ldots, X_k\right)
 \end{split}
 \end{equation*}
 for any $X_1,\ldots, X_k\in M_p\otimes \cW^J_n $, and $p\in \NN$.
 The proof is complete.
 \end{proof}

Let us remark that
 the map $\Phi^J_T$ defined by \eqref{J-poiss-1} can be extended  to a
  completely contractive linear map from $B( \cN_J)$
 to $B(\cH)$ by setting 
$$
 \Phi_T(X)=K_T^* (X\otimes I) K_T,\quad  X\in B(\cN_J).
 $$

We consider now some examples.
 
 \begin{example}
 \label{Sz-ine}
 In the particular case when $n=1$,  $F_1^\infty$ can be identified with 
 the  classical Hardy algebra $H^\infty(\DD)$. It is well-known (see \cite{G})
 that any $w^*$-closed 
 ideal $J$ of $H^\infty$ is generated by an inner function
  $u\in H^\infty(\DD)$, i.e., $J=uH^\infty(\DD)$. 
  In this case, we have 
  $$\cN_J=H^2\ominus uH^2=\ker u(S)^*.
  $$
 Let $T\in B(\cH)$ be a c.n.u. contraction such that $u(T)=0$. Then, according to 
 \cite{SzF-book}, $T$ is of class $C_{\cdot 0}$, which in our notation means
 $T\in C_0$.
 If $\omega$ is  a strongly admissible operator radius, then applying 
 Theorem \ref{J-c0}  to $T$, we obtain
$$
\omega\left(\sum_{j=1}^m f_j(T) g_j(T)^*,\ldots, \sum_{j=1}^m 
\chi_j(T) \psi_j(T)^*\right)\leq
\omega\left( \sum_{j=1}^m f_j(B) g_j(B)^*,\ldots, \sum_{j=1}^m 
\chi_j(B) \psi_j(B)^*\right)
$$
for any functions $f_j, g_j,\ldots,  \chi_j, \psi_j\in H^\infty(\DD)$, $m\in \NN$, 
 where 
$B:=P_{\ker u(S)^*} S| \ker u(S)^*$ and $S$ is the unilateral shift 
 on the Hardy space $H^2$.  
 In the particular case when $k=1$, $m=1$,  and $\omega$ is the operator norm on $B(\cH)$, we get
 $$
 \|f(T) g(T)^*\|\leq \|f(B) g(B)^*\|
 $$
 for any $f,g\in H^\infty$. When $g=1$, we find  the constrained von Neumann
 inequality obtained by Sz.-Nagy \cite{Sz2}.
 
 \end{example}

 \begin{example}\label{nilp}
 Let $J_m$ be the WOT-closed two-sided ideal of $F_n^\infty$
  generated by the monomials $S_\alpha$ for   $\alpha \in \FF_n^+$ with 
  $ |\alpha|=m$.
  It is easy to see that $\cN_{J_m}=\cP_{m-1}$, the set of all polynomials 
  in $F^2(H_n)$
  of degree
  $\leq m-1$. Let $T:=[T_1,\ldots, T_n]$ be a row contraction
  such that 
  $$T_\alpha=0\quad \text{ for any } \quad \alpha \in \FF_n^+ \text{ with } \  |\alpha|=m.
  $$
  Since $T$ is clearly a $C_0$-row contraction, one can apply Theorem $\ref{J-c0}$
  to $T$ and $J_m$. Therefore, if $\omega$ is a  strongly admissible operator radius, then
  \begin{equation*}
  \begin{split}
  \omega\left(\sum_{|\alpha|,|\beta|\leq m-1} a_\alpha T_\alpha T^*_\beta
  \right.,& \left.\ldots,
  \sum_{|\alpha|,|\beta|\leq m-1} b_\alpha T_\alpha T^*_\beta\right)\\
  &\leq 
  \omega\left(\sum_{|\alpha|,|\beta|\leq m-1} a_\alpha 
  S^{(m)}_\alpha {S_\beta^{(m)}} ^*, \ldots,
  \sum_{|\alpha|,|\beta|\leq m-1} b_\alpha S^{(m)}_\alpha {S_\beta^{(m)}}^* \right),
  \end{split}
  \end{equation*}
  where $S_i^{(m)}:= P_{\cP_{m-1}} S_i|\cP_{m-1}$, \ $i=1,\ldots, n$, are the truncated
  left creation operators acting on $\cP_{m-1}$.
  
  Now, let $J_s$ be the WOT-closed two-sided ideal of $F_n^\infty$ generated by the 
  commutators $S_iS_j-S_jS_i$, $i,j=1,\ldots, n$.
  Notice that if $J:=J_m+J_s$, then
  $$
  \cM_J=J_s F^2(H_n)\bigvee J_m F^2(H_n)
  $$
  and  
  $\cN_J=F_s^2(H_n)\cap \cP_{m-1}$,
  where $F_s^2(H_n)\subset F^2(H_n)$ is the symmetric Fock space.
  Now, assume that $T:=[T_1,\ldots, T_n]$ is  a  row contraction such that
  $T_iT_j=T_jT_i$ for any $i,j=1,\ldots, n$, and $T_\alpha=0$ for any 
   $\alpha \in \FF_n^+$, $|\alpha|=m$. Then we can apply  
   Theorem $\ref{J-c0}$
  to $T$ and  the two-sided ideal $J$ and obtain a commutative version 
  of the above inequality.
  More precisely, the operators $S_i^{(m)}$, $i=1,\ldots, n$, are replaced 
  by the mutually commuting
  operators
  $$
  Z_i^{(m)}:= P_{F_s^2(H_n)\cap \cP_{m-1}} S_i |F_s^2(H_n)\cap \cP_{m-1}, 
  \qquad i=1,\ldots, n.
  $$
 \end{example}

  \begin{example}
 Let $W_n^\infty$ be the WOT-closed algebra generated
 by the operators  
 $$
 Z_i:=P_{F^2_s(H_n)} S_i|F^2_s(H_n),\qquad  i=1,\ldots, n,
 $$ and 
 the identity on the symmetric Fock space $F^2_s(H_n)$.
 Fix an operator $\varphi(S_1,\ldots, S_n)\in F_n^\infty$ and let 
 let $T:=[T_1,\ldots, T_n]$ be a $C_0$-row contraction such that
  $T_iT_j=T_jT_i$ for any $i,j=1,\ldots, n$ and
  $$\varphi(T_1,\ldots, T_n)=0.$$
  Let $J_\varphi$ be the WOT-closed two-sided ideal of
   $F_n^\infty$ generated by $\varphi$ and let $J_{\varphi, s}:= J_\varphi + J_s$.
   Straightforward computations reveal that
   \begin{equation*}
   \begin{split}
   \cN_{J_{\varphi, s}}&=F^2_s(H_n)\bigcap \left( \bigcap_{\psi\in J_\varphi}
   \ker \psi(S_1,\ldots, S_n)^*\right)\\
   &=\ker \varphi(Z_1,\ldots, Z_n)^*.
   \end{split}
   \end{equation*}
   Here, we used the fact that the symmetric Fock space is invariant 
   under each operator
   $S_1^*, \ldots, S_n^*$ and $Z_1,\ldots, Z_n$ are mutually commuting.
   Now, we can apply Theorem  \ref{J-c0}  to our setting and obtain the inequality
   \begin{equation*}\begin{split}
   \omega&\left(\sum_{j=1}^m f_j(T_1,\ldots, T_n) g_j(T_1,\ldots, T_n)^*,
   \ldots, \sum_{j=1}^m 
\chi_j(T_1,\ldots, T_n) \psi_j(T_1,\ldots, T_n)^*\right)\\
&\leq
\omega\left( \sum_{j=1}^m f_j(B_1,\ldots, B_n) g_j(B_1,\ldots, B_n)^*,\ldots, 
\sum_{j=1}^m 
\chi_j(B_1,\ldots, B_n) \psi_j(B_1,\ldots, B_n)^*\right)
\end{split}
   \end{equation*}
   for any elements $f_j, g_j,\ldots,  \chi_j, \psi_j\in F_n^\infty$, $m\in \NN$, 
 and
 $$
 B_i:=P_{\ker \varphi(Z_1,\ldots, Z_n)^*} Z_i|\ker \varphi(Z_1,\ldots, Z_n)^*,
 \qquad i=1,\ldots, n.
 $$
 In the particular case when $k=1$, $m=1$,  and $\omega$ is the operator
  norm on $B(\cH)$, we get
 $$
 \|f(T_1,\ldots, T_n) g(T_1,\ldots, T_n)^*\|\leq \|f(B_1,\ldots, B_n)
  g(B_1,\ldots, B_n)^*\|
 $$
 for any $f,g\in F_n^\infty$.
 This inequality is a multivariable commutative version of 
 the constrained von Neumann
 inequality obtained by Sz.-Nagy \cite{Sz2} (case $n=1$ and $g=1$).
 \end{example}

 In what follows, we need the following result.
  
 \begin{lemma}\label{J-le-fo}
  Let $J$ be the WOT-closed two-sided ideal
     of $F_n^\infty$ generated by some homogenous polynomials
      $q_1,\ldots, q_d\in \cP$.  If
$T:=[T_1,\ldots, T_n]$, $T_i\in B(\cH)$,  is  a row
 contraction  such that 
 $$
q_j(T_1,\ldots, T_n)=0,\quad  j=1,\ldots, d,
 $$
   then
 \begin{equation}\label{J-KTr}
 K_{T,r}^* (f(B_1,\ldots, B_n)\otimes I_\cH)=f_r(T_1,\ldots, T_n) K_{T,r}^*
 \end{equation}
 for any $f(S_1,\ldots, S_n)\in  F_n^\infty$ and $0<r<1$.
 Moreover, if $T$ is c.n.c., then
 \begin{equation}\label{J-sot-lim}
 f(T_1,\ldots, T_n)=\text{\rm SOT-}\lim_{r\to 1} K_{T,r}^* 
 (f(B_1,\ldots, B_n)\otimes I_\cH)K_{T,r}
 \end{equation}
 for any $f(S_1,\ldots, S_n)\in F_n^\infty$.
 \end{lemma}
 \begin{proof}
 According to Lemma \ref{le-fo}, we have
 \begin{equation}
 \label{fof}
 [f(S_1,\ldots, S_n)^*\otimes I_\cH]K_{T,r}=K_{T,r} f_r(T_1,\ldots, T_n)^*
 \end{equation}
 for any $f(S_1,\ldots, S_n)\in F_n^\infty$ and $0<r<1$.
 Hence, we deduce that 
 \begin{equation}
 \label{Ktr}
 \left<K_{T,r} k,\varphi\otimes h\right>=
 \left<k, \varphi_r(T_1,\ldots, T_n)\Delta_rh\right>
 \end{equation}
 for any $\varphi(S_1,\ldots, S_n)\in J$ and $h,k\in \cH$.
 Since the $F_n^\infty$- functional calculus is a homomorphism and $q_j$ are homogenous
  polynomials such that $q_j(T_1,\ldots, T_n)=0$, $j=1,\ldots, d$, we have
  \begin{equation}
  \label{omo}
  (fq_j) (rT_1,\ldots, rT_n)=r^{\text{deg}\,(q_j)}
   f(rT_1,\ldots, rT_n)q_j(T_1,\ldots, T_n)=0
  \end{equation}
 for any $f(S_1,\ldots, S_n)\in F_n^\infty$ and $j=1,\ldots, d$, where  deg\,$(q_j)$ denotes 
 the degree
 of the polynomial $q_j$.
 On the other hand,  $J$ is the WOT-closed two-sided ideal of $F_n^\infty$
 generated by  
      $q_1,\ldots, q_d $, and the $F_n^\infty$-functional calculus 
      for $C_0$-row contractions is WOT-continuous. Consequently, since
      $[rT_1,\ldots, rT_n]$ is  a $C_0$-row contraction, relation
      \eqref{omo} implies $\varphi(rT_1,\ldots, rT_n)=0$ for any
       $\varphi(S_1,\ldots, S_n)\in J$ and 
      $0<r<1$.
      Hence,  relation \eqref{Ktr} implies
      $\left<K_{T,r} k,\varphi\otimes h\right>=0$ for any
      $\varphi(S_1,\ldots, S_n)\in J$ and $h,k\in \cH$.
      Since $\{\varphi:\ \varphi(S_1,\ldots, S_n)\in J\}$ is dense in $\cM_J$, the latter equality
      shows that
      \begin{equation}
      \label{inclu}
      K_{T,r}(\cH)\subseteq \cN_J\otimes \cH
      \end{equation}
      for any $0<r<1$.
      Now, relations \eqref{fof}, \eqref{inclu}, and \eqref{fnn} imply
      \begin{equation*}
      \begin{split}
      K_{T,r} f_r(T_1,\ldots, T_n)^*&=\left(P_{\cN_J}\otimes I_\cH\right)
       [f(S_1,\ldots, S_n)^*\otimes I_\cH]\left(P_{\cN_J}\otimes I_\cH\right)
       K_{T,r}\\
       &=[f(B_1,\ldots, B_n)^*\otimes I_\cH]K_{T,r},
      \end{split}
      \end{equation*}
 which proves relation \eqref{J-KTr}.
 Notice that relation \eqref{J-KTr} implies
 \begin{equation}\label{J-fr}
 f_r(T_1,\ldots, T_n)=K_{T,r}^*[f(B_1,\ldots, B_n)\otimes I_\cH] K_{T,r}.
 \end{equation}
 If $[T_1,\ldots, T_n]$ is  a c.n.c. row contraction, then the $F_n^\infty$-functional calculus
 implies 
 $$
 \text{\rm SOT-}\lim_{r\to 1} f_r(T_1,\ldots, T_n)=f(T_1,\ldots, T_n),
 $$
 which together with \eqref{J-fr} imply relation \eqref{J-sot-lim}.
 The proof is complete.
 \end{proof}

 Let $\cA_n^J$  be the norm closed algebra generated by 
$B_i:=P_{\cN_J} S_i|\cN_J$, \ $i=1,\ldots, n$,   and the 
identity $I_{\cN_J}$.
We denote by $\cB$   the linear span of 
$\{B_\alpha B_\beta^*:\ \alpha, \beta\in \FF_n^+\}$ and  by $\overline{\cB}$  its 
norm closure.

Under the additional condition of norm quasi-continuity of the joint operator radius, we can obtain more constrained von Neumann inequalities.

\begin{theorem}\label{J-von-ra}
Let 
    $\omega$ be a joint operator radius, $J$ be the WOT-closed two-sided ideal
     of $F_n^\infty$ generated by some homogenous polynomials
      $q_1,\ldots, q_d\in \cP$.  Let 
$T:=[T_1,\ldots, T_n]$, $T_i\in B(\cH)$,  be  a row
 contraction  such that 
 $$
q_j(T_1,\ldots, T_n)=0,\quad  j=1,\ldots, d.
 $$
 Then  there is a unital completely contractive linear map $\Phi_T^J: \overline{\cB}\to B(\cH)$
 such that $\Phi_T^J(B_\alpha B_\beta^*)=T_\alpha T_\beta^*$, 
 $\alpha, \beta\in \FF_n^+$. Moreover,
 \begin{equation}
 \label{J-wot}
 \Phi_T^J(X)=\lim_{r\to 1} K_{T,r}^* (X\otimes I_\cH) K_{T,r}, 
 \quad X\in \overline{\cB},
 \end{equation}
 where the convergence  is in the operator norm topology,
  and $\{K_{T,r}\}_{0<r<1}$ is the Poisson
  kernel associated with $T$.
 \begin{enumerate}
 \item[(i)]
 If $\omega$ is norm quasi-continuous and  admissible, then
 \begin{equation}\label{J-ine-om3}
  \omega\left([\Phi^J_{T,p}(F_1)]^*, \ldots, [\Phi^J_{T,p}(F_k)]^*\right)\leq 
 \omega\left(F_1^*, \ldots, F_k^*\right)
 \end{equation}
 for any $F_1,\ldots, F_k\in M_p\otimes \cA^J_n $,  and $p\in \NN$.
 \item[(ii)]
 If $\omega$ is  norm quasi-continuous and strongly admissible, then
 \begin{equation}\label{J-ine-om4}
  \omega\left(\Phi^J_{T,p}(X_1), \ldots, \Phi^J_{T,p}(X_k)\right)\leq 
 \omega\left(X_1, \ldots, X_k\right)
 \end{equation}
 for any $X_1,\ldots, X_k\in M_p\otimes \overline{\cB} $ and  $p\in \NN$.
 \end{enumerate}
 \end{theorem}
 \begin{proof}
 According to Lemma \ref{J-le-fo},  we have
 \begin{equation}
 \label{J-KK}
 K_{T,r}^*(B_\alpha B_\beta^*\otimes I_\cH) K_{T,r}=r^{|\alpha|+
 |\beta|} T_\alpha T_\beta^*
 \end{equation}
 for any $\alpha,\beta\in \FF_n^+$ and $0<r<1$. Since $K_{T,r}$ is an isometry,
 relation \eqref{J-KK} implies
 \begin{equation}\label{J-vn}
 \left\|\sum a_{\alpha \beta} T_\alpha T_\beta^*\right\|\leq
  \left\|\sum a_{\alpha \beta} B_\alpha B_\beta^*\right\|
 \end{equation}
 for any   $\sum a_{\alpha \beta} B_\alpha B_\beta^* \in \cB$.
 Moreover, the von Neumann type inequality \eqref{J-vn} implies 
 the existence of a unique unital completely contractive linear map
 $\Phi_T^J: \overline{\cB}\to B(\cH)$
 such that $\Phi_T^J(B_\alpha B_\beta^*)=T_\alpha T_\beta^*$, 
 $\alpha, \beta\in \FF_n^+$. The proof of \eqref{J-wot} is similar to that
  of Theorem 3.8 from \cite{Po-poisson}. We shall omit it.
  
  Now, let us prove part (i) of the theorem.
  Since $\cA_n^J\subset W_J^\infty=P_{\cN_J} F_n^\infty|\cN_J$, Lemma \ref{J-le-fo}
  implies
  \begin{equation}
  \label{KTF}
  \tilde{K}_{T,r} F(rT_1,\ldots, rT_n)^*=[F(B_1,\ldots, B_n)^*\otimes I_\cH]
  \tilde{K}_{T,r}
  \end{equation}
 for any $F(B_1,\ldots, B_n):=[f_{ij}(B_1,\ldots, B_n)]_{p\times p}
 \in M_p\otimes \cA_n^J$, where  $\tilde{K}_{T,r}$ is the ampliation
  of $\tilde{K}_{T,r}$.
  Since  $K_{T,r}$ is an isometry and $\omega$ is admissible, we have
 \begin{equation*}
 \begin{split}
 \omega(&F_1(rT_1,\ldots, rT_n)^*,\ldots, F_k(rT_1,\ldots, rT_n)^*)\\
 &=
 \omega\left((F_1(B_1,\ldots, B_n)^*\otimes I_\cH)|\tilde{K}_{T,r}(\cH^{(p)}),
  \ldots,  
 (F_k(B_1,\ldots, B_n)^*\otimes I_\cH)|\tilde{K}_{T,r}(\cH^{(p)})\right)\\
 &\leq\omega(F_1(B_1,\ldots, B_n)^*\otimes I_\cH, 
 \ldots, F_k(B_1,\ldots, B_n)^*\otimes I_\cH)\\
 &=\omega(F_1(B_1,\ldots, B_n)^*, \ldots, F_k(B_1,\ldots, B_n)^*)
 \end{split}
 \end{equation*}
 for any $F_1(B_1,\ldots, B_n),\ldots, F_k(B_1,\ldots, B_n)\in M_p\otimes \cA^J_n $,  and $p\in \NN$.
 According  to relations \eqref{KTF} and  \eqref{J-wot},  we have
 $$
 \lim_{r\to 1}F_j(rT_1,\ldots, rT_n)=
 \lim_{r\to 1} \tilde{K}_{T,r}^*[F(B_1,\ldots, B_n)^*\otimes I_\cH]
  \tilde{K}_{T,r}=
  \Phi^J_{T,p}(F_j),
  $$
  where the convergence is in operator norm topology.
 Consequently, since $\omega$ is norm quasi-continuous 
  the above inequality implies
  \eqref{J-ine-om3}.  This completes the proof of  
 part (i) of the theorem.
 
 Let us  prove part (ii).
 Since  $\omega$ is strongly admissible and $K_{T,r}$ is an isometry, 
 we have
 \begin{equation*}
 \begin{split}
 \omega(\tilde{K}_{T,r}^* (X_1\otimes I_\cH)\tilde{K}_{T,r},\ldots,  
 \tilde{K}_{T,r}^* (X_k\otimes I_\cH)\tilde{K}_{T,r})&\leq
 \omega(X_1\otimes I_\cH, \ldots, X_k\otimes I_\cH)\\
 &=\omega(X_1, \ldots, X_k)
 \end{split}
 \end{equation*}
 for any  $X_1,\ldots, X_k\in M_p\otimes  \overline{\cB}$.
  Now, using relation \eqref{J-wot}, the norm continuity of $\omega$, 
 and taking 
 $r\to 1$, we deduce inequality \eqref{J-ine-om4}. The proof is complete.
 \end{proof}

Let $T:=[T_1,\ldots, T_n]$, $T_i\in B(\cH)$,  be  a c.n.c. row
 contraction and let
 $J$ be a WOT-closed two-sided ideal
     of $F_n^\infty$ 
 such that 
 $$
 \varphi(T_1,\ldots, T_n)=0,\quad \varphi(S_1,\ldots, S_n)\in J.
 $$
 Let  
$$
\cD^J_n:=\text{\rm span}\left\{ f(B_1,\ldots, B_n)B_\alpha^*:
\ f(S_1,\ldots, S_n)\in F_n^\infty, \alpha\in \FF_n^+
\right\}
$$
  and   define the linear map  $\Phi^J_T:\cD^J_n+(\cD^J_n)^*\to B(\cH)$ 
by setting
\begin{equation}\label{J-def-doug}
\Phi^J_T\bigl(f(B_1,\ldots, B_n) B_\alpha^*+ B_\beta g(B_1,\ldots, B_n)^*\bigr)
:=f(T_1,\ldots, T_n) T_\alpha^*+ T_\beta
g(T_1,\ldots, T_n)^*
\end{equation}
for any $f(S_1,\ldots, S_n), g(S_1,\ldots, S_n)\in F_n^\infty$, 
and $\alpha,\beta\in \FF_n^+$.

\begin{theorem}\label{J-cnc}
Let 
    $\omega$ be a joint operator radius and  $J$ be the WOT-closed two-sided ideal
     of $F^\infty_n$ generated by some homogenous polynomials
      $q_1,\ldots, q_d\in \cP$.  Let 
$T:=[T_1,\ldots, T_n]$, $T_i\in B(\cH)$,  be  a c.n.c. row
 contraction  such that 
 $$
q_j(T_1,\ldots, T_n)=0,\quad  j=1,\ldots, d.
 $$
  Then, the map $\Phi^J_T$ defined by \eqref{J-def-doug}   can be extended 
   to a unital  completely contractive linear map
from the norm closure of $\cD^J_n+(\cD^J_n)^*$ to $B(\cH)$  by setting
\begin{equation}\label{J-def-ext}
\Phi^J_T(X):=\text{\rm WOT-}\lim_{r\to 1}K_{T,r}^* [X\otimes I)]
 K_{T,r}, \quad X\in \overline{\cD^J_n+(\cD^J_n)^*},
 \end{equation}
where   $\{K_{T,r}\}_{0<r<1}$  is the Poisson
kernel associated with $T$.
  \begin{enumerate}
  \item[(i)] If $\omega$ is $*$-SOT quasi-continuous and  admissible, then
  \begin{equation}\label{J-ine-omi}
 \omega\left([\Phi^J_{T,p}(F_1)]^*, \ldots, [\Phi^J_{T,p}(F_k)]^* \right)\leq 
 \omega\left(F_1^*, \ldots, F_k^* \right)
 \end{equation}
 for any $F_1,\ldots, F_k\in M_p\otimes W_J^\infty$, and $p\in \NN$.
   \item[(ii)]
   If $\omega$ is SOT quasi-continuous and strongly 
 admissible, then 
 \begin{equation}\label{J-ine-om5}
 \omega\left(\Phi^J_{T,p}(Y_1), \ldots, \Phi^J_{T,p}(Y_k) \right)\leq 
 \omega\left(Y_1, \ldots, Y_k \right)
 \end{equation}
 for any $Y_1,\ldots, Y_k\in M_p\otimes  \overline{\cD^J_n} $, and $p\in \NN$.
 \item[(iii)]
  If  $\omega$ is WOT quasi-continuous and strongly 
 admissible, then  inequality \eqref{J-ine-om5}
 holds for any $Y_1,\ldots, Y_k\in M_p\otimes  \overline{\cD^J_n+(\cD^J_n)^*}$, and 
 $p\in \NN$.
 \end{enumerate}
\end{theorem}
 
 \begin{proof}
 According to Lemma \ref{J-le-fo}, we have
 \begin{equation}\label{J-K*K}
 K_{T,r}^* (f(B_1,\ldots, B_n)p(B_1,\ldots, B_n)^* \otimes I_\cH)K_{T,r}=
 f_r(T_1,\ldots, T_n) p_r(T_1,\ldots, T_n)^*  
 \end{equation}
for any  $f(S_1,\ldots, S_n)\in F_n^\infty$ and any polynomial $p\in \cP$.
Since
$$
\|f_r(T_1,\ldots, T_n)\|\leq \|f(B_1,\ldots, B_n)\|\quad \text{ for any } 
\ 0<r<1,
$$
 and 
$p_r(T_1,\ldots, T_n)\to p(T_1,\ldots, T_n)$ in norm as $r\to 1$, 
one can use relations
\eqref{funct-cal} and \eqref{J-K*K} to deduce that
\begin{equation}\label{J-SOT1}
\text{\rm SOT-}\lim_{r\to 1} 
K_{T_r}^* [f(B_1,\ldots, B_n)p(B_1,\ldots, B_n)^* \otimes I_\cH]K_{T,r}=
 f(T_1,\ldots, T_n) p(T_1,\ldots, T_n)^*  
\end{equation}
 Hence, we  obtain
 \begin{equation}\label{J-SOT2}
 \begin{split}
\text{\rm WOT-}\lim_{r\to 1} 
K_{T_r}^* \bigl[\bigl(q(B_1,&\ldots, B_n)g(B_1,\ldots, B_n)^*\\
&+
f(B_1,\ldots, B_n)p(B_1,\ldots, B_n)^*\bigr) \otimes I_\cH\bigr]K_{T,r}\\
&
 =q(T_1,\ldots, T_n)g(T_1,\ldots, T_n)^*+
 f(T_1,\ldots, T_n) p(T_1,\ldots, T_n)^*  
 \end{split}
\end{equation}
for any $f,g\in F_n^\infty$ and $p,q\in \cP$. 
This shows that the linear map $\Phi_T:\cD^J_n+ (\cD^J_n)^*\to B(\cH)$ defined by
\eqref{J-def-doug} is  completely contractive and  the limit in \eqref{J-def-ext} exists
for any 
 $X\in \cD^J_n+ (\cD^J_n)^*$.
The rest of the proof is similar to that of Theorem \ref{von-ra-cnc}. We shall omit it.
 \end{proof}

 \begin{corollary} 
  If $X\in \overline{\cD_n^J}$, then 
$$
\Phi^J_T(X)=\text{\rm SOT-}\lim_{r\to 1}K_{T,r}^* [X\otimes I)]
 K_{T,r}.
 $$
 \end{corollary}

\begin{remark}\label{comm-ine}
Let $J_s$ be the WOT-closed two-sided ideal of $F_n^\infty$
generated  by the homogenous
polynomials
$$
\{S_iS_j-S_jS_i:\ i,j=1,\ldots, n\}.
$$
It is  easy to see that the subspace $\cN_{J_s}$ is equal to the symmetric
 Fock space $F_s^2(H_n)\subset F^2(H_n)$.
Let $T:=[T_1,\ldots, T_n]$, $T_i\in B(\cH)$, be a row contraction such that
$$
T_iT_j=T_jT_i,\quad i,j=1,\ldots, n.
$$
Then Theorem $\ref{J-von-ra}$ can be applied to the $n$-tuple $T$. As a   particular case,
 if $\omega$
is norm quasi-continuous and strongly 
admissible operator radius, 
we obtaind the following inequality
\begin{equation}\label{com-ineq}
\omega\left( \sum a_{\alpha\beta}^{(1)} T_\alpha T_\beta^*, \ldots 
 \sum a_{\alpha\beta}^{(k)} T_\alpha T_\beta^*\right)
 \leq 
 \omega\left( \sum a_{\alpha\beta}^{(1)} B_\alpha B_\beta^*, \ldots 
 \sum a_{\alpha\beta}^{(k)} B_\alpha B_\beta^*\right)
\end{equation}
for any polynomias $\sum a_{\alpha\beta}^{(j)} B_\alpha B_\beta^*$, 
$j=1,\ldots, k$, 
in $C^*(B_1,\ldots, B_n)$, where $B_i:=P_{F_s^2(H_n)} S_i|F_s^2(H_n)$, $i=1,\ldots, n$, are are the creation operators
on the symmetric Fock space $F_s^2(H_n)$.
We should remark that when $k=1$ and $\omega$ is the operator norm the inequality
\eqref{com-ineq} was obtained  by the author in \cite{Po-poisson} and
Arveson in  \cite{Arv1}.

We also remark that if $T$ is a $C_0$- (resp. c.n.c.) row contraction, 
one can apply Theorem $\ref{J-c0}$ (resp. Theorem $\ref{J-cnc}$) to $T$ and obtain 
commutative versions of these theorems.

\end{remark}

\smallskip

 \section{Multivariable  Haagerup--de la Harpe   inequalities
  }\label{HaHa}

In \cite{HD}, Haagerup and de la Harpe proved that any nilpotent 
contraction $T\in B(\cH)$ with  $T^m=0$, $m\geq 2$, satisfies the inequality 
$$
\sup_{\|h\|=1} |\left< Th, h\right>|=\omega(T)\leq \cos\frac {\pi} {m+1}.
$$

We apply the results of the previous sections
 to obtain several multivariable generalizations of the 
 Haagerup--de la Harpe  inequality.
 Using   a result of Boas and Kac
  \cite{BoK} as generalized by Janssen \cite{J}, we also obtain an
   epsilonized  version of the above inequality, when one gives up 
   the condition $T_\alpha=0$. This  extends a  recent result of
    Badea and Cassier \cite{BC}, to our setting.

First, we prove the following  nilpotent dilation result, which extends a result of Arveson
\cite{Ar}, to our multivariable (noncommutaive and commutative) setting.
\begin{theorem}\label{nilpotent}
 Let $m\geq 2$ be a nonnegative 
integer and let \ 
$[T_1,\ldots, T_n]$, \ $T_i\in B(\cH)$,  be an $n$-tuple of operators.
The following statements are equivalent:
\begin{enumerate}
\item[(i)] There exist a Hilbert space $\cK\supseteq \cH$ and a row contraction
$[N_1,\ldots, N_n]$, $N_i\in B(\cK)$, such that $N_\alpha=0$ for any 
$\alpha\in \FF_n^+$, $|\alpha|=m$, and 
\begin{equation}\label{nildil1}
p(T_1,\ldots, T_n)=P_\cH p(N_1,\ldots, N_n) |\cH\quad \text{ for any  }\ 
  p\in \cP_{m-1};
\end{equation}

\item[(ii)]
There is a Hilbert space $\cG$ such that $\cH$ can be identified with a
subspace of $ \cP_{m-1}\otimes\cG$
and 
\begin{equation}\label{nildil2}
 p(T_1,\ldots, T_n)=P_\cH p(S_1^{(m)}\otimes I_\cG,\ldots, 
 S_n^{(m)}\otimes I_\cG) |\cH\quad \text{ for any  }\ 
  p\in \cP_{m-1},
\end{equation}
  where $S_i^{(m)}:=P_{\cP_{m-1} }
 S_i |\cP_{m-1} $, \ $i=1,\ldots, n$;
\item[(iii)]   The multi-Toeplitz operator
 \begin{equation}
 \label{posi}
 \sum_{1\leq |\alpha|\leq m-1} 
  R_\alpha^* \otimes T_{\tilde \alpha}+ I\otimes I+
 \sum_{1\leq |\alpha|\leq m-1} 
 R_\alpha \otimes T_{\tilde \alpha}^* 
 \end{equation}
is positive.
 \end{enumerate}
 
 If $T_1,\ldots T_n $ are mutually commuting operators, then the above statements
 are  also equivalent 
 to 
 \begin{enumerate}
 
\item[(iv)]
There is a Hilbert space $\cG$ such that $\cH$ can be identified with a
subspace of $ \cP^s_{m-1}\otimes\cG$
and 
\begin{equation}\label{co-nildil2}
 p(T_1,\ldots, T_n)=P_\cH p(B_1^{(m)}\otimes I_\cG,\ldots, 
 B_n^{(m)}\otimes I_\cG) |\cH\quad \text{ for any  }\ 
  p\in \cP_{m-1},
\end{equation}
  where $B_i^{(m)}:=P_{\cP^s_{m-1} }
 S_i |\cP^s_{m-1} $, \ $i=1,\ldots, n$, and 
 $\cP^s_{m-1}:=P_{F^2_s(H_n)} (\cP_{m-1})$.
 \end{enumerate}
\end{theorem}

\begin{proof}
Assume that condition (i) holds. According to Theorem \ref{vN}, we have 
$$
\sum_{1\leq |\alpha|\leq m-1} 
  r^{|\alpha|} R_\alpha^* \otimes N_{\tilde \alpha}+ I+
 \sum_{1\leq |\alpha|\leq m-1} 
r^{|\alpha|} R_\alpha \otimes N_{\tilde \alpha}^* \geq 0
 $$
for any $0<r<1$. Taking $r\to 1$, we get (iii). Since the implication 
(ii)$\implies$(i) is clear, it remains to prove that 
(iii)$\implies$(ii). To this end, we assume that the multi-Toeplitz operator
given by \eqref{posi}  is positive.
Define the linear map
$$
\Phi: \text{\rm span}\{S^{(m)}_\alpha:\ |\alpha|\leq m-1\}\to 
 \text{\rm span}\{T_\alpha:\ |\alpha|\leq m-1\}
 $$
 by setting $\Phi(S^{(m)}_\alpha)=T_\alpha$, \ $|\alpha|\leq m-1$, where 
 $S^{(m)}_i:=P_{\cP_{m-1}} S_i |\cP_{m-1}$, $i=1,\ldots, n$.
 In what follows, we prove that $\Phi$ is completely contractive.
 Define the sequence of operators $A_{(\sigma)}:=T_\sigma$ if $0\leq |\sigma|\leq m-1$
 and $A_{(\sigma)}:=0$ if $|\sigma|\geq m$.
 Consider the multi-Toeplitz kernel $K:\FF_n^+\times \FF_n^+\to B(\cH)$ defined by
 $K(g_0, \sigma)=K(\sigma, g_0)^*:=A_{(\sigma)}$.
 As in the proof of Theorem \ref{vN}, one can show that condition (iii) implies 
 the positivity of the operator 
 $M_q:=[K(\omega, \sigma)]_{|\omega|, |\sigma|\leq q}$ for any $q=1,2,\ldots$, 
 i.e., the kernel $K$ is positive definite.
 According to Theorem 3.1 of \cite{Po-posdef}, there is a completely positive linear map
 $\mu:C^*(S_1,\ldots, S_n)\to B(\cH)$ such that
 $\mu(S_\sigma)=A_{(\sigma)}$ for any $\sigma\in \FF_n^+$.
 Since $\mu(I)=I$, $\mu$ is completely contractive.
 Therefore, $\mu_0:=\mu|\cA_n$ is completely contractive and 
 \begin{equation*}
 \mu_0(S_\sigma)=\begin{cases}T_\sigma \quad &\text{ if } |\sigma|\leq m-1\\
 0\quad &\text{ if } |\sigma|\geq m.
 \end{cases}
 \end{equation*}
 Note that $\mu_0$ is not multiplicative
  (in general $T_\sigma\neq 0$ if $|\sigma|\geq m$). However, since $\mu_0(J)=0$,
   where $J$ is the two-sided ideal generated by $\{S_\alpha:\ |\alpha|=m\}$, 
   $\mu_0$ induces a completely linear map 
   $\hat{\mu}_0:\cA_n/J\to B(\cH)$ such that 
   $$
   \hat{\mu}_0(S_\sigma+J)=T_\sigma, \qquad  |\sigma|\leq m-1.
   $$
   On the other hand, according to Theorem 4.1  of \cite{ArPo2}, 
   the linear map 
   $$
\Psi: \text{\rm span}\{S^{(m)}_\alpha:\ |\alpha|\leq m-1\}\to \cA_n/J
  $$
  defined by 
   $\Psi(S^{(m)}_\alpha)=S_\alpha+J$, \ $|\alpha|\leq m-1$, is completely isometric
    and onto the quotient $\cA_n/J$.
    Since $\Phi=\hat{\mu}_0\circ \Psi$, it follows 
    that $\Phi$ is completely contractive.
    Since $\Phi(I)=I$, Arveson's extension theorem shows that $\Phi$ has a 
    completely positive extension $\tilde\Phi$ to  $C^*(S^{(m)}_1,\ldots, S^{(m)}_n)$.
    By Stinespring's theorem, there exists a representation 
    $\pi:C^*(S^{(m)}_1,\ldots, S^{(m)}_n)\to B(\cK)$ such that 
    $\tilde\Phi(x)=P_\cH \pi(x)|\cH$ for any
     $x\in C^*(S^{(m)}_1,\ldots, S^{(m)}_n)$.
    Note that if $p,q$ are polynomials in $F^2(H_n)$, then
    \begin{equation}\label{pq}
    p(S_1,\ldots, S_n) P_{\CC} q(S_1,\ldots, S_n)^*\xi=
    \left< \xi, q(S_1,\ldots, S_n)(1)\right>
    p(S_1,\ldots, S_n)(1)
    \end{equation}
   for any $\xi\in F^2(H_n)$, where $ P_{\CC}=I-S_1S_1^*-\cdots -S_nS_n^*$.
    The operator $p(S_1,\ldots, S_n) P_{\CC} q(S_1,\ldots, S_n)^*$
     has rank one  and  takes $q$ into $\|q\|p$.
    If $p, q$ are polynomials in $\cP_{m-1}$ and since $\cP_{m-1}$ is an invariant subspace  under each operator
     $S_1^*, \ldots, S_n^*$, we  deduce that
     $$
     p(S^{(m)}_1,\ldots, S^{(m)}_n) (I-S^{(m)}_1{S^{(m)}_1}^*-\cdots -
     S^{(m)}_n{S^{(m)}_n}^*) q(S^{(m)}_1,\ldots, S^{(m)}_n)^*
     $$
    is a rank one operator acting on  $\cP_{m-1}$.
    Therefore,
    $C^*(S^{(m)}_1,\ldots, S^{(m)}_n)=B(\cP_{m-1})$  and, consequently,  
    $\pi$ must be unitarily equivalent to a multiple of the identity 
    representation. In particular, we have 
    $$
    T_\sigma=P_\cH \pi(S^{(m)}_\sigma)|\cH=P_\cH 
    ( S^{(m)}_\sigma \otimes I_\cG)|\cH
    $$
    for any $|\sigma|\leq m-1$, where $\cG$ is a separable
     Hilbert space.
    
     Now,   assume that $T_1,\ldots, T_n$ are commuting operators.
     As in the noncommutative case,  but now the ideal $J$ is replaced by the two
      sided ideal generated by $\{S_\alpha:\ |\alpha|=m\}$ and the
       commutators $S_iS_j-S_jS_i$, $i,j=1,\ldots, n$, 
       we find a $*$-representation 
      $\pi:C^*(B^{(m)}_1,\ldots, B^{(m)}_n)\to B(\cK)$ such that 
    $\tilde\Phi(x)=P_\cH \pi(x)|\cH$ for any
     $x\in C^*(B^{(m)}_1,\ldots, B^{(m)}_n)$.
  
     On the other hand,
     taking the compresion of \eqref{pq} to
      $\cP_{m-1}^s:=P_{F_s^2(H_n)}(\cP_{m-1})$ and taking into account that 
     the subspace  $\cP_{m-1}^s$ is invariant under each operator
      $S_1^*, \ldots, S_n^*$, we infer that
      $$C^*(B^{(m)}_1,\ldots, B^{(m)}_n)=B(\cP_{m-1}^s).
      $$
      Now the result follows,  as in  to the noncommutative case.
      This completes the proof.
\end{proof}

\begin{corollary}\label{ine-nilp}
  If the multi-Toeplitz operator defined by condition \eqref{posi} 
  is positive and $\omega$ is any  strongly admissible operator
 radius, then
 $$
 \omega\bigl(p_1(T_1,\ldots, T_n), \ldots, p_k(T_1,\ldots, T_n)\bigr)\leq
 \omega\bigl(p_1(S_1^{(m)},\ldots, S_n^{(m)}\bigr), \ldots, p_k(S_1^{(m)},\ldots, S_n^{(m)}))
$$
for any polynomials $p_1,\ldots, p_k\in \cP_{m-1}$.
If, in addition, $T_1,\ldots T_n $ are mutually commuting operators, then
$$
 \omega\bigl(p_1(T_1,\ldots, T_n), \ldots, p_k(T_1,\ldots, T_n)\bigr)\leq
 \omega\bigl(p_1(B_1^{(m)},\ldots, B_n^{(m)}\bigr), \ldots, p_k(B_1^{(m)},\ldots, B_n^{(m)})).
$$
\end{corollary}

We remark that in the particular case when $[T_1,\ldots, T_n]$ is a row contraction
such that $T_\alpha=0$ for any $\alpha\in \FF_n^+$ with $|\alpha|=m$, we already obtained more general inequalities in Example \ref{nilp}.

Let $S^{(m)}$ be the $m$-dimensional shift on $\CC^m$, given by the matrix
$$
\left(\begin{matrix}
0&0&0&\cdots& 0& 0\\
1&0&0&\cdots& 0& 0\\
0&1&0&\cdots& 0& 0\\
\vdots&\vdots&\vdots&\ddots&\vdots&\vdots\\
0&0&0&\cdots& 0& 0\\
0&0&0&\cdots& 1& 0\\
\end{matrix}
\right).
$$

It is well-known (see \cite{GR}, \cite{HD}) that  
\begin{equation}\label{nrf}
w(S^{(m)})=
 \cos\frac{\pi}{m+1}
 \end{equation}
and
$\left<S^{(m)}\xi,\xi\right>= \cos\frac{\pi}{m+1}$, where the vector
$\xi=(\xi_0,\ldots, \xi_{m-1})\in \CC^m$ is given
by
 \begin{equation}\label{coe}
 \xi_k=\sqrt{\frac{2}{m+1}}\sin\frac{km}{m+1}, \qquad k=0,1,\ldots, m-1.
 \end{equation}

In what follows, we obtain  a multivariable version of formula \eqref{nrf}.
\begin{theorem}\label{truncate}
Let $S_i^{(m)}:=P_{\cP_{m-1}}S_i|\cP_{m-1}$, \ $i=1,2,\ldots,  n$, be
 the truncated
left creation operators  acting on the set of all polynomials 
in $F^2(H_n)$ of degree $\leq m-1$.
Then the joint numerical range $W(S_1^{(m)},\ldots, S_n^{(m)})$ is a convex compact
subset of $\CC^n$ and
\begin{equation}\label{wwcos}
w_e(S_1^{(m)},\ldots, S_n^{(m)})=w(S_1^{(m)},\ldots, S_n^{(m)})
=\cos\frac{\pi}{m+1}.
\end{equation}
\end{theorem}
\begin{proof}
The fact that the joint numerical range 
$W(S_1^{(m)},\ldots, S_n^{(m)})$ is a convex compact
subset of $\CC^n$  is a consequence of Theorem \ref{spectru} part (i).
Now, notice that $S_\alpha^{(m)}=0$ for any 
$\alpha\in \FF_n^+$ with $|\alpha|=m$, and consequently
$$
(S_1\otimes {S_1^{(m)*}}+\cdots+ S_n\otimes S_n^{(m)*})^m=0.
$$
Applying Theorem \ref{nilpotent} (when $n=1$) to 
the operator $T:=S_1\otimes S_1^{(m)*}+\cdots+ S_n\otimes S_n^{(m)*}$ acting on
$\cH:=F^2(H_n)\otimes \cP_{m-1}$, we find a countable Hilbert space $\cG$
such that $\cH\subseteq \CC^m\otimes \cG$ and
$$
T^k=P_\cH(S^{(m)}\otimes I_\cG)^k|\cH,\quad  k=1,2,\ldots, m-1.
$$
Hence, using  Lemma \ref{prop2} and formula \eqref{nrf}, we get
$$
w(T)\leq w(S^{(m)}\otimes I_\cG)= w(S^{(m)})=\cos\frac{\pi}{m+1}.
$$
Since $w(S_1^{(m)},\ldots, S_n^{(m)})=w(T)$ (see   Corollary \ref{pro}),
  we infer that
\begin{equation}\label{ineleft}
w(S_1^{(m)},\ldots, S_n^{(m)})\leq \cos\frac{\pi}{m+1}.
\end{equation}

On the other hand, according to Theorem \ref{conv}, we have
\begin{equation}\label{wew2}
w_e(S_1^{(m)},\ldots, S_n^{(m)})\leq w(S_1^{(m)},\ldots, S_n^{(m)}).
\end{equation}
Now, define the vector $q\in \cP_{m-1}$ by setting
$$q:=\xi_0 1+\sum_{k=1}^{m-1} \xi_k e_1^k,
$$
where the coefficients $\xi_k$ are those  defined by \eqref{coe}.
Since $\|q\|_2=1$ and $\left<S_i^{(m)} q, q\right>=0$ for any $i=2,3,\ldots,n$, we
can use \eqref{nrf} and obtain
\begin{equation*}
\begin{split}
w_e(S_1^{(m)},\ldots, S_n^{(m)})&=
\sup\left\{ \left(\sum_{i=1}^n\left|\left<S_i^{(m)} p,p\right>\right|^2\right)^{1/2}
:\  p\in \cP_{m-1}, ~\|p\|=1\right\}\\
&\geq \left(\sum_{i=1}^n\left|\left<S_i^{(m)} q,q\right>\right|^2\right)^{1/2}
=\left|\left<S_1^{(m)} q, q\right>\right|\\
&=\left|\left<S^{(m)} \xi, \xi\right>\right|
=\cos\frac{\pi}{m+1}.
\end{split}
\end{equation*}
Hence, and using inequalities \eqref{ineleft} and \eqref{wew2}, we deduce  
\eqref{wwcos}. The proof is complete.
\end{proof} 

We should remark that $q$ is not unique with the above property.
Indeed, we can take for example $q=\xi_0 1+\sum_{k=1}^{m-1} \xi_k e_i^k$, 
where $i=2,\ldots, n$.
\begin{corollary}
Let $m\geq 2$ and 
let $B_i^{(m)}:=P_{\cP_{m-1}}B_i|\cP_{m-1}$, \ $i=1,2,\ldots,  n$, be
 the truncated
left creation operators  acting on the set of all polynomials 
in $F_s^2(H_n)$ of degree $\leq m-1$. Then
\begin{equation}\label{wwcos-B}
w_e(B_1^{(m)},\ldots, B_n^{(m)})=w(B_1^{(m)},\ldots, B_n^{(m)})
=\cos\frac{\pi}{m+1}.
\end{equation}
\end{corollary}

Now, we can prove the following multivariable generalization  of  Haagerup-de la Harpe inequality.

\begin{theorem}\label{gen-HH}  Let $m\geq 2$ be a nonnegative 
integer and let  $T_1,\ldots, T_n\in B(\cH)$ be such that
% $[T_1,\ldots, T_n]$ is a row contraction and
  the multi-Toeplitz operator
 \begin{equation}
 \label{posi1}
 \sum_{1\leq |\alpha|\leq m-1}  R_\alpha^* \otimes T_{\tilde \alpha}+ I+
 \sum_{1\leq |\alpha|\leq m-1}  R_\alpha \otimes T_{\tilde \alpha}^* 
 \end{equation}
is positive.
 Then,  for each $1\leq k\leq m-1$,
\begin{equation}
\label{HH1}
  w(T_\alpha:\ |\alpha|=k)
 \leq 
 \cos\frac {\pi} {\left[\frac{m-1} {k}\right]+2},
\end{equation}
where $[x]$ is the integer part of $[x]$.
\end{theorem}
\begin{proof}
According to Theorem \ref{nilpotent}, we have
$$
T_\beta= P_\cH (S_\beta^{(m)}\otimes I_\cG)|\cH,\quad |\beta|\leq m-1.
$$
Therefore, for each $1\leq k\leq m-1$, we have
\begin{equation}
\label{ww}
 w(T_\alpha:\ |\alpha|=k)\leq  w(S_\alpha^{(m)}\otimes I_\cG:\ |\alpha|=k).
 \end{equation}
 Let $X_1,\ldots, X_{n^k}$ be an  arrangement of the operators
 $S_\alpha^{(m)}\otimes I_\cG$, where $\alpha\in \FF_n^+$ and $|\alpha|=k$, 
 in a sequence.
 Notice that $X_\beta=0$ for any  
  $\beta\in \FF_n^+$  with 
  $$
  |\beta|=q:=\left[\frac{m-1} {k}\right] +1.
  $$ 
  Applying again Theorem \ref{nilpotent} to the row contraction
  $[X_1,\ldots, X_{n^k}]$, we find a Hilbert space $\cM$ such that 
  $$
  X_j= P_{\cP_{m-1}\otimes \cG}(S_j^{(q)}\otimes I_\cM)|\cP_{m-1}\otimes \cG, 
  \quad j=1,\ldots, n^k,
  $$
  where $S_1^{(q)},\ldots, S_{n^k}^{(q)}$ are the truncated left creation operators on 
  the full Fock space $F^2(H_{n^k})$.
  Hence, using Theorem \ref{propri} and Theorem \ref{truncate}, we obtain
  \begin{equation*}
  \begin{split}
  w(X_1,\ldots, X_{n^k})&\leq
  w(S_1^{(q)}\otimes I_\cM,\ldots, S_{n^k}^{(q)}\otimes I_\cM)\\
  &=w(S_1^{(q)},\ldots, S_{n^k}^{(q)})\\
  &\leq \cos \frac {\pi} {q+1}.
  \end{split}
  \end{equation*}
  Hence, and using inequality  \eqref{ww}, we obtain \eqref{HH1}. 
  The proof is complete.
\end{proof}

\begin{corollary}\label{w-vN}
 Let  $[T_1,\ldots, T_n]$ be  a row contraction
   such that the multi-Toeplitz operator 
 given by \eqref{posi1} is positive,  and let 
 $p :=\sum\limits_{|\alpha|\leq m-1}
 a_\alpha e_\alpha$ be a polynomial.
    Then
\begin{equation}\label{Ine}
w\left(p(T_1,\ldots, T_n) \right)
\leq  K\|p\|_2,
\end{equation}
where the constant $K$ is given by
$$
K:=\left(\sum_{k=0}^{m-1} \cos^2\frac {\pi} 
 {\left[\frac{m-1} {k}\right]+2}\right)^{1/2}.
 $$
\end{corollary}
\begin{proof}
Since,  according to 
Proposition \ref{2-ineq},
$$
w_e(T_\alpha:\ |\alpha|=k)\leq w(T_\alpha:\ |\alpha|=k),
$$
 Theorem \ref{gen-HH} implies
 \begin{equation}\label{weco}
 \sup_{\|h\|=1}\left(\sum_{|\alpha|=k}\left|
 \left<T_\alpha h,h\right>\right|^2\right)^{1/2}\leq
  \cos\frac {\pi} {\left[\frac{m-1} {k}\right]+2}.
 \end{equation}
 On the other hand, according to Theorem \ref{propri} the joint numerical radius is a norm.
 Applying Cauchy's inequality twice and using \eqref{weco}, we obtain
 \begin{equation*}
 \begin{split}
 w\left(p(T_1,\ldots, T_n) \right) &\leq \sum_{k=0}^{m-1}
 w\left( \sum_{|\alpha|=k} a_\alpha T_\alpha\right)\\
 &=\sum_{k=0}^{m-1}\sup_{\|h\|=1} \left|\sum_{|\alpha|=k}
  a_\alpha \left< T_\alpha h,h\right>\right|\\
 &\leq \sum_{k=0}^{m-1}\left[ \left(\sum_{|\alpha|=k} |a_\alpha|^2\right)^{1/2}
 \sup_{\|h\|=1} \left(\sum_{|\alpha|=k}
  \left|\left< T_\alpha h,h\right>\right|^2\right)^{1/2}\right]\\
  &\leq 
  \sum_{k=0}^{m-1}\left[\left(\sum_{|\alpha|=k} |a_\alpha|^2\right)^{1/2}
  \cos\frac {\pi} {\left[\frac{m-1} {k}\right]+2}\right]\\
  &\leq 
  \left(\sum_{|\alpha|\leq m-1} |a_\alpha|^2\right)^{1/2}
 \left( \sum_{k=0}^{m-1}\cos^2\frac {\pi} {\left[\frac{m-1} {k}\right]+2}
 \right)^{1/2},
 \end{split}
 \end{equation*}
  where we make the convention that the term corresponding to $k=0$ 
  in the latter summation is 1.
  The proof is complete.
\end{proof}

Notice that since $\|[T_1,\ldots, T_n]\|\leq 2 w(T_1,\ldots, T_n)$, inequality
\eqref{Ine} implies 
$$
\|p(T_1,\ldots, T_n)\|
\leq  2K\|p\|_2
$$
for any polynomial $p\in \cP_{m-1}$.
On the other hand, since
$$\sigma(p(T_1,\ldots, T_n))\subset 
\overline{\left(\CC^n\right)}_{w(p(T_1,\ldots, T_n))},
$$
one can use Corollary \ref{w-vN} to locate the spectrum
of the operator $p(T_1,\ldots, T_n)$.

 As   consequences of Theorem \ref{gen-HH}, we obtain the following 
  multivariable 
Haagerup--de la Harpe type inequalities.

\begin{corollary}\label{HA}
Let $T_1,\ldots, T_n\in B(\cH)$ be such that $T_\alpha=0$ for any 
$\alpha\in \FF_n^+$, $|\alpha|=m$
$(m\geq 2)$. Then 
\begin{equation}
\label{HH2}
  w(T_\alpha:\ |\alpha|=k)
\leq 
 \left\|\sum_{i=1}^n T_i T_i^*\right\|^{k/2} \cos\frac {\pi} 
 {\left[\frac{m-1} {k}\right]+2},
\end{equation}
for  $1\leq k\leq m-1$.
If, in addition,  $[T_1,\ldots, T_n]$ is a row contraction and
$p :=\sum\limits_{|\alpha|\leq m-1}
 a_\alpha e_\alpha$ is a polynomial,
    then
 inequality \eqref{Ine} holds.
\end{corollary}
\begin{proof}
Since $[T_1,\ldots, T_n]$ is a row contraction and 
$T_\alpha=0$ for any 
$\alpha\in \FF_n^+$, $|\alpha|=m$, we can use Theorem \ref{vN} to deduce that the 
operator given by \eqref{posi1} is positive. Applying Theorem \ref{gen-HH}, the result follows.
If $T_1,\ldots, T_n$ are arbitrary operators, one can use the homogeneity
of the joint numerical radius (see Theorem \ref{propri}) to deduce
the inequality \eqref{HH2}.
\end{proof}

\begin{corollary}\label{HA2}
Let $T_1,\ldots, T_n\in B(\cH)$ be such that $T_\alpha=0$ for any 
$\alpha\in \FF_n^+$, $|\alpha|=m$
$(m\geq 2)$. Then 
\begin{equation}
\label{HH3}
  w_e(T_\alpha:\ |\alpha|=k)
\leq \|(T_\alpha:\ |\alpha|=k)\|_e
 \cos\frac {\pi} 
 {\left[\frac{m-1} {k}\right]+2},
\end{equation}
for each $1\leq k\leq m-1$.
\end{corollary}
\begin{proof}
Since $(\lambda_1T_1+\cdots+\lambda_n T_n)^m=0$ for any
$(\lambda_1,\ldots, \lambda_n)\in \BB_n$, we can use Corollary
\ref{HA} in the particular case when $k=1$ and deduce that
$$
w(\lambda_1T_1+\cdots+\lambda_n T_n)\leq \|\lambda_1T_1+\cdots+\lambda_n T_n\|
\cos\frac {\pi} 
 { m +1}.
 $$
 Consequently, taking the supremum over $(\lambda_1,\ldots, \lambda_n)\in \BB_n$
and using Corollary  \ref{we-sup2}, we obtain
\begin{equation}\label{INe} 
w_e(T_1,\ldots, T_n)\leq \|(T_1\ldots, T_n)\|_e \cos\frac {\pi} 
 { m +1}.
\end{equation}
 Now, let us consider the general case when $ 1\leq k\leq m-1$.
  Let $Y_1,\ldots, Y_{n^k}$
 be an arrangement of $\{T_\alpha:\ |\alpha|=k\}$.
 Notice that $Y_\beta=0$ for any $\beta\in \FF_n^+$, $|\beta|=\left[\frac{m-1} {k}\right] +1$.
 Applying the first part of the proof to $(Y_1,\ldots, Y_{n^k})$,
 we deduce inequality \eqref{HH3}.
\end{proof}

Now, we present some results concerning the range of the joint numerical range and the 
euclidean operator radius  for $n$-tuples of operators.

  \begin{theorem}\label{range}  The joint numerical radius and the euclidean  operator  radius
  have the following properties:
    \begin{equation}\label{range1}
    \left\{w(T_1,\ldots, T_n):\ \|[T_1,\ldots, T_n]\|=1\right\}=
  \left[\frac {1} {2}, 1\right]
  \end{equation} 
  and 
   \begin{equation}\label{range2}
   \left\{w_e(T_1,\ldots, T_n):\ \|[T_1,\ldots, T_n]\|_e=1\right\}=
  \left[\frac {1} {2}, 1\right].
 \end{equation}
  Moreover, $\text{ \rm range } w=\text{ \rm range } w_e=[0,\infty)$.
  \end{theorem}
  \begin{proof}
  According to Theorem \ref{propri}, if $\|[T_1,\ldots, T_n]\|=1$, then
  $\frac {1} {2} \leq w(T_1,\ldots, T_n)\leq 1$.
  For each $i=1,\ldots,n$, and  $t\in [0,1]$,  define the operator 
  $T_i(t)$ on the Fock space $F^2(H_n)$ by
  $$
  T_i(t)\left(\sum_{\alpha\in \FF_n^+} a_\alpha e_\alpha\right):=
  a_0e_{g_i}+t\sum_{|\alpha|\geq 1} a_\alpha e_{g_i\alpha}.
  $$
  Note that $\|[T_1(t), \ldots, T_n(t)]\|=1$ for any $t\in [0,1]$.
  On the other hand,  using Theorem \ref{example} and Theorem \ref{truncate}, 
  we have
  $$
  w(T_1(1), \ldots, T_n(1))=w(S_1,\ldots, S_n)=1$$
  and 
  $$
  w(T_1(0), \ldots, T_n(0))=w(S_1^{(2)},\ldots, S_n^{(2)})=\frac{1} {2}
  $$
  Since the joint numerical radius is norm continuous  (see Theorem \ref{propri}),
   we deduce relation \eqref{range1}.
   Now, one can easily  show that 
   $\|[T_1(t), \ldots, T_n(t)]\|_e=1$ for any $t\in [0,1]$. As above, 
   using Theorem \ref{example}, Theorem \ref{truncate}, and the norm continuity
   of $w_e$ (see Theorem \ref{propri2}) we complete the proof.
  \end{proof}
  
  For each $m\geq 2$, define
  $$\cN_m:=\left\{ (T_1,\ldots, T_n)\in B(\cH)^{(n)}:\ 
  \|[T_1,\ldots, T_n]\|=1, \text{ and } T_\beta=0, \ |\beta|=m\right\}
  $$
  and
  $$\cN_m^e:=\left\{ (T_1,\ldots, T_n)\in B(\cH)^{(n)}:\ 
  \|[T_1,\ldots, T_n]\|_e=1, \text{ and } T_\beta=0, \ |\beta|=m\right\}
  $$
  
\begin{theorem} If $m\geq 2$, then
$$
\text{\rm range }w|\cN_m=\text{\rm range }w_e|\cN_m^e=\left[ \frac {1} {2}, ~
\cos\frac {\pi} 
 {m +1}\right].
 $$
\end{theorem}
\begin{proof}
According to Theorem \ref{propri} and Corollary \ref{HA}, if $\|[T_1,\ldots, T_n]\|=1$ 
and $T_\beta=0$ for any
$\beta\in \FF_n^+$ with $|\beta|=m$, then
  $$\frac {1} {2} \leq w(T_1,\ldots, T_n)\leq \cos\frac {\pi} 
 {m +1}.
 $$
 Similarly,  
 if $\|[T_1,\ldots, T_n]\|_e=1$ 
and $T_\beta=0$ for any
$\beta\in \FF_n^+$ with $|\beta|=m$, then, according to Theorem \ref{propri2}
 and Corollary \ref{HA2}, we have
  $$\frac {1} {2} \leq w_e(T_1,\ldots, T_n)\leq \cos\frac {\pi} 
 {m +1}.
 $$
  For each $i=1,\ldots,n$ and  $t\in [0,1]$,  define the operator 
  $T^{(m)}_i(t)$ on  $\cP_{m-1}$ by
  $$
  T^{(m)}_i(t)\left(\sum_{|\alpha|\leq m-1} a_\alpha e_\alpha\right):=
  a_0e_{g_i}+t\sum_{1\leq |\alpha|\leq m-2} a_\alpha e_{g_i\alpha}.
  $$
  Straightforward computations reveal  that 
  $$
  \|[T^{(m)}_1(t), \ldots, T^{(m)}_n(t)]\|=
  \|[T^{(m)}_1(t), \ldots, T^{(m)}_n(t)]\|_e=1
  $$
   for any $t\in [0,1]$.
  On the other hand,  using   Theorem \ref{truncate}, 
  we have
  \begin{equation*}\begin{split}
  w(T^{(m)}_1(1), \ldots, T^{(m)}_n(1))&=w_e(T^{(m)}_1(1), \ldots, T^{(m)}_n(1))\\
 &= w(S^{(m)}_1,\ldots, S^{(m)}_n)\\
 &=\cos\frac {\pi} 
 {m +1}
 \end{split}
 \end{equation*}
  and 
  \begin{equation*}
  \begin{split}
  w(T^{(m)}_1(0), \ldots, T^{(m)}_n(0))&=w_e(T^{(m)}_1(0), \ldots, T^{(m)}_n(0))\\
  &=
  w(S_1^{(2)},\ldots, S_n^{(2)})=\frac{1} {2}
  \end{split}
  \end{equation*}
  Since the joint numerical radius  and the euclidean numerical radius are 
   norm continuous  (see Theorem \ref{propri} and Theorem \ref{propri2}),
   we  can complete the proof.
\end{proof}

Another multivariable generalization of Haagerup-de la Harpe inequality is the following.

\begin{theorem}\label{Ha-la}
Let $m\geq 2$ and $T_1,\ldots, T_n\in B(\cH)$ be such that
\begin{equation} \label{lala}
\sum_{|\alpha|\leq m-1}\overline{\lambda}_\alpha T_\alpha^* +I +
\sum_{|\alpha|\leq m-1} {\lambda}_\alpha T_\alpha\geq 0
\end{equation}
for any $(\lambda_1,\ldots, \lambda_n)\in \BB_n$. Then the euclidean radius 
satisfies 
the inequality
\begin{equation}\label{we-cos}
w_e(T_1,\ldots, T_n)\leq \cos\frac {\pi} 
 {m +1}.
 \end{equation}
 Moreover, if  
$1\leq k\leq m-1$, then
\begin{equation}
\label{int-ine}
\sup_{|h|=1} \left[\sum_{|\alpha|=k}\left( \int_{\partial\BB_n}|
\xi_\alpha|^2 d\sigma(\xi)\right)\left|\left< T_\alpha h,h\right>\right|^2\right]^{1/2}
\leq \cos\frac {\pi}{\left[\frac { m-1} {k}\right] +2}.
\end{equation}
\end{theorem}
\begin{proof}
Let $\theta\in \RR$, $(\xi_1,\ldots, \xi_n)\in \BB_n$, and denote
 $\lambda_j:=e^{i\theta}\xi_j$, $j=1,\ldots, n$.
 Note that condition \eqref{lala} is equivalent to 
 \begin{equation*}
 \sum_{k=1}^{m-1} e^{-ik\theta} \left( \overline{\xi}_1T_1^*+\cdots +
  \overline{\xi}_n T_n^*\right)^k+
 I+\sum_{k=1}^{m-1} e^{ik\theta} \left( {\xi}_1 T_1+\cdots +
  {\xi}_n T_n\right)^k
 \end{equation*}
 for any  $\theta\in \RR$ and  $(\xi_1,\ldots, \xi_n)\in \BB_n$.
 The latter inequality is equivalent to
 $$
 \sum_{k=1}^{m-1} {S^*}^k \otimes \left( \overline{\xi}_1T_1^*+\cdots +
  \overline{\xi}_n T_n^*\right)^k+
 I+\sum_{k=1}^{m-1} {S}^k \otimes \left( {\xi}_1 T_1+\cdots +
  {\xi}_n T_n\right)^k\geq 0,
  $$
  where $S$ is the unilateral shift on the Hardy space $H^2$.
 Using Theorem \ref{gen-HH} (in the particular case $n=k=1$) in our setting,
 we obtain
 $$
 w_e(\xi_1T_1+\cdots+ \xi_n T_n)=
 w(\xi_1T_1+\cdots+ \xi_n T_n)\leq \cos\frac {\pi}{m+1}.
 $$
 Taking the supremum over $(\xi_1,\ldots, \xi_n)\in \BB_n$ and using
 Corollary \ref{we-sup2}, we obtain inequality \eqref{we-cos}.
 Similarly,  using Theorem \ref{gen-HH}, we deduce that
 \begin{equation}\label{xi-k}
 w((\xi_1T_1+\cdots+ \xi_nT_n)^k)\leq 
  \cos{\frac {\pi}{\left[\frac { m-1} {k}\right] +2}}.
 \end{equation}
 On the othe hand, we have
 \begin{equation*}\begin{split}
 \sup_{(\xi_1,\ldots, \xi_n)\in \BB_n}w((\xi_1T_1+\cdots+ \xi_nT_n)^k)&\geq 
 \sup_{(\xi_1,\ldots, \xi_n)\in \BB_n}\sup_{\|h\|=1}
 \left|\sum_{|\alpha|=k} \xi_\alpha \left< T_\alpha h, h\right>\right|\\
 &=\sup_{\|h\|=1}\left[\sup_{(\xi_1,\ldots, \xi_n)\in \BB_n}
 \left|\sum_{|\alpha|=k, |\beta|=k}\xi_\alpha \overline{\xi}_\beta \left< T_\alpha h, h\right>
\overline{ \left< T_\beta h, h\right>}\right|\right]^{1/2}\\
&\geq 
\sup_{\|h\|=1} 
\left|\sum_{|\alpha|=k, |\beta|=k} \int_{\partial \BB_n}
\xi_\alpha \overline{\xi}_\beta \left< T_\alpha h, h\right>
\overline{ \left< T_\beta h, h\right>} \right|^{1/2}\\
&= \sup_{\|h\|=1}\left[\sum_{|\alpha|=k}\left(\int_{\partial \BB_n}|\xi_\alpha|^2d 
\sigma(\xi)\right)\left|\left< T_\alpha h, h\right>\right|^2\right]^{1/2}.
 \end{split}
 \end{equation*}
 Hence, and using inequality \eqref{xi-k}, we deduce \eqref{int-ine}.
 The proof is complete.
\end{proof} 
 
 Let us point out two particular cases of Theorem \ref{Ha-la}.
 If $n=1$, we get a result of Badea-Cassier (see \cite{BC}).
 If $n=2$ and $T_1T_2=T_2T_1$, then
 $$
 w_e(T_1^2, T_2^2, T_1T_2)\leq \sqrt{3} 
 \cos\frac {\pi}{\left[\frac { m-1} {2}\right] +2}.
 $$
 
 We need to recall from \cite{BC} the following epsilonized
  version of Fej\' er inequality, which is a consequence of a result of Boas and Kac
  (\cite{BoK}) as generalized by Janssen (\cite{J}).
  
  If $f(e^{i\theta})=\sum_{k\in \ZZ} c_ke^{ik\theta}$ is positive such that
  $\sum_{k\in \ZZ} |c_k|<\infty$ with $c_0$=1 and $c_k\leq \epsilon$ for any 
  $k\geq m$, then
  $$
  |c_1|\leq \cos\left( \frac {\pi} {m+1}\right)+
  3\left[\pi \cos^4\frac{\pi} {2(m+1)}\right]^{1/3} 
  \left(\frac {\epsilon} {m+1}\right)^{2/3}.
  $$
 We denote by $K(\epsilon, m)$   the right side of this inequality.
 \begin{theorem}\label{epsi}
 Let $m\geq 2$ be a positive integer and let $[T_1,\ldots, T_n] $ be 
 a row contraction 
 such that
 $$
 \left\|\sum_{|\alpha|=m} T_\alpha T_\alpha^*\right\|^{1/2}\leq
  \epsilon\quad \text{ and }
 \quad
 \sum_{k=0}^\infty\left\|\sum_{|\alpha|=k} T_\alpha T_\alpha^*\right\|^{1/2}<\infty.
 $$
 Then
the joint numerical radius satisfies the inequality
 \begin{equation*}
 w(T_\alpha:\ |\alpha|=k)\leq
 \cos\left( \frac {\pi} {\left[\frac {m-1} {k}\right]+2}\right)+
  3\left[\pi \cos^4\frac{\pi} {2\left(\left[\frac {m-1} {k}\right]+2\right)}\right]^{1/3} 
  \left(\frac {\epsilon} {\left[\frac {m-1} {k}\right]+2}\right)^{2/3}.
  \end{equation*}
 \end{theorem}
 \begin{proof}
 Since $[e^{i\theta} T_1,\ldots, e^{i\theta}T_n]$ is a row contraction
  for any $\theta\in \RR$, Theorem \ref{vN} implies that the operator
 $$
 \sum_{k=1}^\infty \sum_{|\alpha|=k}r^k 
  e^{-ik\theta}R_\alpha^*\otimes T_{\tilde \alpha} +I+
  \sum_{k=1}^\infty \sum_{|\alpha|=k}r^k 
  e^{ik\theta}R_\alpha\otimes T_{\tilde \alpha}^* 
  $$
  is positive for any $0<r<1$ and $\theta\in \RR$,  
  where the convergence is in the operator norm.
   Hence, we deduce that
   \begin{equation}\label{X*}
   \sum_{k=1}^\infty  r^k e^{-ik\theta} {X^*}^k +I + 
   \sum_{k=1}^\infty  r^k e^{ik\theta} {X}^k\geq 0,
   \end{equation}
   where $X:= \sum_{i=1}^n R_i\otimes T_i^*$.
   Since
   $$
   \sum_{k=1}^\infty \|X^k\|= 
   \sum_{k=1}^\infty\left\|\sum_{|\alpha|=k}
    T_\alpha T_\alpha^*\right\|^{1/2}<\infty
    $$
    and taking $r\to 1$ in \eqref{X*}, we obtain
    $$
    \sum_{k=1}^\infty
    e^{-ik\theta}\left<{X^*}^k y,y\right>+1+
    \sum_{k=1}^\infty
    e^{ik\theta}\left<{X}^k y,y\right>\geq 0
    $$
    for any $\theta\in \RR$ and $y\in F^2(H_n)\otimes \cH$. 
    Fix $y\in F^2(H_n)\otimes \cH$ with $\|y\|=1$ and denote
    $$
    c_k:=\left<{X}^k y,y\right> \quad \text{ and } \quad 
     c_{-k}:=\left<{X^*}^k y,y\right>
     $$
     for any $k=0,1,\ldots$.
     Notice that 
     $\sum_{k\in \ZZ} |c_k|<\infty$ and, for any $q\geq m$, we have
     \begin{equation*}
     \begin{split}
     c_q&=\left| \left<{X}^q y,y\right>\right|\leq
     \left\|\sum_{|\alpha|=q}
    T_\alpha T_\alpha^*\right\|^{1/2}\\
    &\leq \left\|\sum_{|\alpha|=m}
    T_\alpha T_\alpha^*\right\|^{1/2}\left\|\sum_{|\alpha|=q-m}
    T_\alpha T_\alpha^*\right\|^{1/2}
    \leq \epsilon.
     \end{split}
     \end{equation*}
    Applying the result preceeding this theorem, we obtain
    $$|\left<Xy,y\right>|=|c_1|\leq K(\epsilon, m)
    $$
    for any  $y\in F^2(H_n)\otimes \cH$ with $\|y\|=1$.
    Using now Corollary \ref{pro} and the unitary invariance of 
    the joint numerical radius,
    we obtain inequality  
    \begin{equation*}
 w(T_1,\ldots, T_n)\leq
 \cos\left( \frac {\pi} { m+1}\right)+
  3\left[\pi \cos^4\frac{\pi} {2 (m+1)}\right]^{1/3} 
  \left(\frac {\epsilon} {m+1}\right)^{2/3},
  \end{equation*}
  which proves the inequality of the theorem, when $k=1$.
    Now, let us consider the general case when $ 1\leq k\leq m-1$.
  Let $Y_1,\ldots, Y_{n^k}$
 be an arrangement of $\{T_\alpha:\ |\alpha|=k\}$.
 Notice that   if 
  $q:==\left[\frac{m-1} {k}\right] +1$,
  then
  $$
 \left\|\sum_{|\beta|=q} Y_\beta Y_\beta^*\right\|^{1/2}\leq
  \epsilon\quad \text{ and }
 \quad
 \sum_{k=0}^\infty\left\|\sum_{|\beta|=k} Y_\beta Y_\beta^*\right\|^{1/2}<\infty.
 $$
 Applying the first part of the proof to $(Y_1,\ldots, Y_{n^k})$,
  we complete the proof of the theorem.
 \end{proof}

 The proof of the next result is similar to that of Corollary \ref{w-vN}. We shall omit it.
 \begin{corollary}
 Under the conditions of Theorem $\ref{epsi}$, the following inequalities hold:
 $$
 \frac{1} {2}\|p(T_1,\ldots, T_n)\|\leq w(p(T_1,\ldots, T_n))\leq
 \|p\|_2 \left(\sum_{k=0}^m K^2\left(\epsilon, 
 \left[\frac{m-1} {k}\right]+1\right)\right)^{1/2}
 $$
 for any polynomial $p:=\sum_{|\alpha|\leq m-1} a_\alpha e_\alpha$.
 \end{corollary}

\smallskip

 \section{Multivariable Fej\' er  inequalities  
  }\label{fejer}

 Haagerup and de la Harpe  proved (\cite{HD}) that their inequality 
 is equivalent to
Fej\' er's inequality for positive trigonometric polynomials
$\sum\limits_{k=-m+1}^{m-1} a_k e^{ik\theta}$, namely
$$
|a_1|\leq a_0 \cos\frac {\pi} {m+1}.
$$
 In what follows we  prove  multivariable  
  analogues of classical inequalities  (Fej\' er \cite{Fe},
   Egerv\' ary-Sz\' azs \cite{ES}) for the coefficients of positive
   trigonometric  polynomials, 
   and  of recent extensions to positive rational functions, 
   obtained by Badea and Cassier
   \cite{BC}.

The main result of this section is the following 
   
\begin{theorem}\label{gen-Fe}  Let $m\geq 2$ be a nonnegative 
integer and let  $\left\{A_{(\alpha)}\right\}_{|\alpha|\leq m-1}$,  be 
a sequence of operators in $B(\cH)$ such that
  the multi-Toeplitz operator
 \begin{equation}
 \label{posi2}
 \sum_{1\leq |\alpha|\leq m-1}  R_\alpha^* \otimes A_{(\alpha)}+ I\otimes 
 A_{0}+
 \sum_{1\leq |\alpha|\leq m-1}  R_\alpha \otimes A_{(\alpha)}^* 
 \end{equation}
is positive.
 Then,  for each $1\leq k\leq m-1$,
\begin{equation}
\label{fe1}
 w_e(A_{(\alpha)}:\ |\alpha|=k)\leq
w(A_{(\alpha)}:\ |\alpha|=k)
\leq 
 \|A_{0}\|\cos\frac {\pi} {\left[\frac{m-1} {k}\right]+2},
\end{equation}
where $[x]$ is the integer part of $[x]$.
\end{theorem}
\begin{proof}
Let $F$ denote the multi-Toeplitz operator defined by \eqref{posi2} and assume that
 $F\geq 0$ and $F\neq 0$. 
 According to the Fej\' er type factorization theorem for multi-Toeplitz
  polynomials \cite{Po-analytic}, there exists a multi-analytic operator
  $G:= \sum_{|\alpha|\leq m-1} R_\alpha\otimes C_{(\alpha)}$ such that $F=G^* G$.
  Hence, we have
  \begin{equation}
  \label{FGG}
  F=\sum_{|\alpha|, |\beta|\leq m-1} R_\alpha^* R_\beta \otimes C_{(\alpha)}^*
  C_{(\beta)}.
  \end{equation}
  This implies  $A_0=\sum_{|\alpha|\leq m-1} C_{(\alpha)}^*  C_{(\alpha)}\geq 0$.
  Assume for the moment that
  $A_0=I$.
  Given $\gamma\in \FF_n^+$ such that $|\gamma|\leq m-1$,  relation
  \eqref{FGG} implies 
  \begin{equation}\label{Agi}
  A_{(\gamma)}=\sum_{|\alpha|\leq m-|\gamma|-1} C_{(\alpha)}^* 
   C_{(\alpha \gamma)}.
  \end{equation}
  Define the operator $V:\cH\to \cP_{m-1}\otimes \cH$ be setting
  $$
  Vh:= \sum_{|\alpha|\leq m-1} C_{(\alpha)}^*h\otimes e_\alpha.
  $$
   Since $A_0=I$, we have
   $$
   \|Vh\|^2=\sum_{|\alpha|\leq m-1} \|C_{(\alpha)}^*h\|^2=\left<A_0h,h\right>
   =\|h\|^2,
   $$ 
   so $V$ is an isometry.
   On the other hand, we have
   \begin{equation*}
   \begin{split}
   \left<Vh, (I\otimes R_\gamma^{(m)}) Vk\right>&=
   \left<\sum_{|\beta|\leq m-1} C_{(\beta)}^*h\otimes e_\beta,
   \sum_{|\alpha|\leq m-\gamma-1} C_{(\alpha)}^*k\otimes e_{\alpha \tilde{\gamma}}
   \right>\\
   &=\sum_{|\alpha|\leq m-|\gamma|-1}\left< C_{(\alpha  
   \tilde{\gamma})} h, C_{(\alpha)}\right>\\
   &=\left< h, A_{( \tilde{\gamma})}\right>
   \end{split}
   \end{equation*}
  for any $h.k\in \cH$.
  Therefore,
  $$A_{( \tilde{\gamma})}= V^* (I\otimes R_\gamma^{(m)}) V
  $$
  for any  $\gamma\in \FF_n^+$ such that $|\gamma|\leq m-1$.
 Since $R_\alpha^{(m)}=0$ for any $|\alpha|=m$, we can apply Corollary
 \ref{HA}   and deduce that, for each $1\leq k\leq m-1$,
 \begin{equation*}
 \begin{split}
 w(A_{(\gamma)}:\ |\gamma|=k)
 &\leq w(V^* (I\otimes R^{(m)}_{\gamma}) V:\ |\gamma|=k)\\
 &\leq w(I\otimes R^{(m)}_{\gamma}:\ |\gamma|=k)\\
  &\leq w(  R^{(m)}_{\gamma}:\ |\gamma|=k)\\
  &\leq\cos\frac {\pi} {\left[\frac{m-1} {k}\right]+2}.
 \end{split}
 \end{equation*}
 
 Now let us consider the case when $A_0\geq 0$ but $A_0\neq I$.
 Let $\epsilon >0$ and  define the multi-Toeplitz operator 
 \begin{equation*}
 F_\epsilon:= 
 \sum_{1\leq |\alpha|\leq m-1}  R_\alpha^* \otimes C_{(\alpha)}+ I\otimes 
 C_{0}+
 \sum_{1\leq |\alpha|\leq m-1}  R_\alpha \otimes C_{(\alpha)}^*,
 \end{equation*} 
 where 
 \begin{equation}\label{calp}
 C_{(\alpha)}:=(A_0+\epsilon I)^{-1/2} A_{(\alpha)} (A_0+\epsilon I)^{-1/2} 
 \quad \text{ \rm if }\quad
 1\leq |\alpha|\leq m-1,
 \end{equation} 
 and $C_0=I$. 
 Since $F_\epsilon\geq 0$, we can apply the first part of the proof  and obtain
  \begin{equation}
\label{fe2}
w(C_{(\alpha)}:\ |\alpha|=k)
\leq 
  \cos\frac {\pi} {\left[\frac{m-1} {k}\right]+2}.
\end{equation}
 Now, using  Theorem \ref{propri} part (v)   and relation \eqref{calp},
  we obtain
 \begin{equation*}
 \begin{split}
 w(A_{(\alpha)}:\ |\alpha|=k)&=
 w(X_\epsilon ^*C_{(\alpha)}X_\epsilon:\ |\alpha|=k)\\
 &\leq
 \|X_\epsilon\|^2 ~w(C_{(\alpha)}:\ |\alpha|=k),
 \end{split}
 \end{equation*} 
 where $X_\epsilon:= A_0+\epsilon I$.
 Consequently, using inequality \eqref{fe2}  and taking $\epsilon \to 0$,
  we obtain  the second inequality in 
 \eqref{fe1}. The first inequality  of  \eqref{fe1} is a consequence of 
 Proposition \ref{2-ineq}. This completes the proof of the theorem.
\end{proof}

  In the particular case when $n=1$, we obtain the following operatorial version
  of the Fej\' er and Egerv\'ary-Sz\'azs inequalities.
  
  \begin{corollary}\label{oper-Fe}  Let $m\geq 2$ be a nonnegative 
integer and let  $\left\{A_{k}\right\}_{k=1}^{m-1}$  be 
a sequence of operators in $B(\cH)$ such that
  the  operator
 \begin{equation*}
 \sum_{1\leq k\leq m-1}  e^{-ik\theta} A_{k}^*+  
 A_{0}+
 \sum_{1\leq k\leq m-1} e^{ik\theta} A_k 
 \end{equation*}
is positive for any $e^{i\theta}\in \TT$.
 Then,  for each $1\leq k\leq m-1$,
\begin{equation}
\label{op-fe1}
w(A_k )
\leq 
 \|A_{0}\|\cos\frac {\pi} {\left[\frac{m-1} {k}\right]+2}.
\end{equation}
 \end{corollary}

Let  $\FF_n$ be the free group with generators $g_1,\ldots, g_n$, and let $\ell^2(\FF_n)$ the
Hilbert space defined by
$$
\ell^2(\FF_n):=\{f:\FF_n\to \CC:\sum\limits_{\sigma\in\FF_n}
|f(\sigma)|^2<\infty\}.
$$
The canonical basis of 
$\ell^2(\FF_n)$ is $\{\xi_ \sigma\}_{ \sigma\in\FF_n}$, where $\xi_ \sigma(t)=1$ if $t= \sigma$ and $\xi_ \sigma(t)=0$
otherwise.
For each $i=1,\dots,n$, let  $U_i\in B(\ell^2(\FF_n))$ be the
unitary operator defined by
$$
U_i\left(\sum\limits_{ \sigma\in\FF_n}a_ \sigma \xi_ \sigma\right):
=\sum\limits_{\sigma\in\FF_n} a_ \sigma \xi_{g_i \sigma}, \qquad
\left(\sum\limits_{\sigma\in\FF_n}| a_ \sigma|^2<\infty\right).
$$
 The reduced  group $C^*$-algebra  $C_{\text{\rm red}}^*(\FF_n)$
is the $C^*$-algebra generated by $U_1,\dots, U_n$.
The full group $C^*$-algebra $C^*(\FF_n)$ is generated by an $n$-tuple of universal
unitaries $\bf{U}_1,\ldots, \bf{U}_n$.

   \begin{theorem} \label{SUU} Let 
   $$
   f(X_1,\ldots, X_n):=\sum_{1\leq|\alpha|\leq m-1} \overline{a}_\alpha X_\alpha^* +a_0 I+
   \sum_{1\leq|\alpha|\leq m-1} a_\alpha X_\alpha,
   $$
   where $a_\alpha\in \CC$ and $X_1,\ldots, X_n$ are bounded linear operators on a Hilbert
    space.
   If $f(S_1,\ldots, S_n)$  \ (resp. $f(U_1,\ldots, U_n)$, \ 
   $f(\bf{U}_1,\ldots, \bf{U}_n))$
   is a positive element in $C^*(S_1,\ldots, S_n)$ 
   (resp. $C_{\text{\rm red}}^*(\FF_n)$, $C^*(\FF_n)$), then
   $$
   \left(\sum_{|\alpha|=k} |a_\alpha|^2\right)^{1/2}\leq a_0
   \cos\frac{\pi} {\left[\frac{m-1}{k}\right]+2}
   $$
    for each $1\leq k\leq m-1$, 
   where $[x]$ is the integer part of $x$.
   \end{theorem}
\begin{proof} If $f(S_1,\ldots, S_n)\geq 0$, then the result is a particular case
 of Theorem
\ref{gen-Fe}.
 To prove the second part of the theorem, notice that
the Hilbert space $\ell^2(\FF_n^+)$ can be seen as a subspace of 
$\ell^2(\FF_n)$ and  the full Fock space $F^2(H_n)$ can be naturally identified
to $\ell^2(\FF_n^+)$.
Under this identification, we have  
$U_i|_{F^2(H_n)}=S_i$, $i=1,\dots,n$, 
where $S_1,\dots,S_n$ are the left creation operators.  Consequently, we have
$$
f(S_1,\ldots, S_n)= P_{\ell^2(\FF_n^+)}f(U_1,\ldots, U_n)| \ell^2(\FF_n^+).
$$
Therefore, if  $f(U_1,\ldots, U_n)\geq 0$, then  $f(S_1,\ldots, S_n)\geq 0$ and
 the result  follows from the first part of the proof.
 Now assume that $f(\bf{U}_1,\ldots, \bf{U}_n)\geq 0$. 
 Due to the universal property of 
 $\bf{U}_1,\ldots, \bf{U}_n$, there is a $*$-representation $\pi$ of 
 $C^*(\bf{U}_1,\ldots, \bf{U}_n)$ onto $C^*(U_1,\ldots, U_n)$ such 
 that $\pi({\bf{U}}_i)=U_i$,
 \ $i=1,\ldots, n$. Therefore,
 $
 f(U_1,\ldots, U_n)=\pi(f(\bf{U}_1,\ldots, \bf{U}_n))\geq 0.
 $
 Using again the  first part of the proof, the result follows.
\end{proof}

We remark that, in a similar manner,  one can use Theorem \ref{gen-Fe} to obtain an operator-valued 
version of 
Theorem \ref{SUU} for the tensor products $C_{\text{\rm red}}^*(\FF_n)\otimes B(\cH)$ and $C^*(\FF_n)\otimes B(\cH)$.

The next result is a multivariable generalization of a recent result obtained by Badea and Cassier (see Theorem 6.3 of \cite{BC}).
Since the proof is similar to the proof of Theorem \ref{nilpotent}, we shall omit it.
\begin{theorem}\label{nilpotent2}
 Let $m\geq 2$ be a positive
integer 
and let  $\left\{A_{(\alpha)}\right\}_{ \alpha\in \FF_n^+}$  be 
a sequence of operators in $B(\cH)$ such that
 the multi-Toeplitz operator
 \begin{equation*}
 \sum_{1\leq |\alpha|\leq m-1} 
  R_\alpha^* \otimes A_{ (\alpha)}+ I\otimes I+
 \sum_{1\leq |\alpha|\leq m-1} 
 R_\alpha \otimes A_{(\alpha)}^* 
 \end{equation*}
is positive.  Then
there is a Hilbert space $\cG$ such that $\cH$ can be identified with 
a subspace of $ \cP_{m-1}\otimes\cG$
and 
\begin{equation*} 
  A_{(\alpha)}=P_\cH (S_\alpha^{(m)}\otimes I_\cG) |\cH\quad \text{ for any  }\ 
 |\alpha|\leq m-1,
\end{equation*}
  where $S_i^{(m)}:=P_{\cP_{m-1} }
 S_i |\cP_{m-1} $, \ $i=1,\ldots, n$.
 
 Moreover, if, in addition, $A_{(\alpha)}=A_{(\beta)}$ whenever $\lambda_\alpha=\lambda_\beta$ for any
 $(\lambda_1,\ldots, \lambda_n)\in \BB_n$, then
 \begin{equation*} 
  A_{(\alpha)}=P_\cH (B_\alpha^{(m)}\otimes I_\cG) |\cH\quad \text{ for any  }\ 
 |\alpha|\leq m-1,
\end{equation*}
where $B_i^{(m)}:=P_{\cP_{m-1}^s} S_i|
\cP_{m-1}^s$, $i=1,\ldots, n$.
\end{theorem}

Finnaly, we mention that Theorem \ref{nilpotent2} and Theorem \ref{truncate} can be used to provide another proof of 
Theorem \ref{gen-Fe}.

 \bigskip
 \bigskip

\end{document}